\documentclass[11pt,final,leqno,onefignum,onetabnum]{siamltex}
\usepackage{fullpage}
\usepackage{mathtools}  
\mathtoolsset{showonlyrefs}  
\newtheorem{remarkk}{Remark}
\let\olddefinition\remarkk
\renewcommand{\remarkk}{\olddefinition\normalfont}  
\usepackage[algo2e,linesnumbered,ruled]{algorithm2e}
\usepackage{amsmath,amssymb}
\usepackage{amsopn}
\usepackage{color}
\usepackage{xr}
\usepackage{mathtools}
\usepackage{amssymb}
\usepackage{enumerate}
\usepackage{hyperref}
\usepackage{color, colortbl}
\usepackage{graphicx}
\usepackage{caption}
\usepackage{subcaption}
\usepackage{verbatim}
\usepackage{algorithmic,algorithm}

\definecolor{Blue}{rgb}{0,0,1}
\definecolor{Red}{rgb}{1,0,0}
\definecolor{Purple}{rgb}{0.5, 0, 0.5}
\definecolor{grey}{gray}{0.6}

\newcommand{\comp}{b}
\newcommand{\ideal}{a}

\newcommand{\bm}[1]{\ensuremath{{\mathbf #1}}}

\newcommand{\nat}[1]{\ensuremath{\mathbb N(#1)}} 
\newcommand{\innat}[1]{\ensuremath{\in\nat{#1}}}
\newcommand{\innatseq}[1]{\ensuremath{= 1,\ldots,#1}}

\newcommand{\RR}[1]{\mathbb{R}^{#1}}
\newcommand{\RRstar}[1]{\mathbb{R}_\star^{#1}}

\newcommand{\grassman}[2]{\mathcal G\left(#1,#2\right)}
\newcommand{\innerProd}[3]{\left(#1,#2\right)_{#3}}
\newcommand{\innerProdAj}[2]{\innerProd{#1}{#2}{\Aj}}
\newcommand{\normAj}[1]{\left\|#1\right\|_{\Aj}}
\newcommand{\normA}[1]{\left\|#1\right\|_{\A}}
\newcommand{\normOutput}[1]{\left\|#1\right\|_{\outputMat^T\outputMat}}
\newcommand{\A}{{\mathbf A}}
\newcommand{\Aj}{\A_j}
\newcommand{\Ajm}{\A_{j-1}}
\newcommand{\Ajmi}{\A_{j-i}}

\newcommand{\M}{{\mathbf M}}
\renewcommand{\prec}{\M}
\newcommand{\precj}{\M_j}
\newcommand{\precInvj}{\precj^{-1}}
\newcommand{\precInv}{\prec^{-1}}
\newcommand{\precNonlinlNo}{\bar \M^{l}}
\newcommand{\precNonlinl}[1]{\bar \M^{l}(#1)}
\newcommand{\precNonlinNo}{\bar \M}
\newcommand{\precNonlin}[1]{\bar \M(#1)}
\newcommand{\objectiveNo}{g}
\newcommand{\objective}[3]{\objectiveNo_{#1}^{#2}(#3)}
\newcommand{\objectiveMap}[2]{\objectiveNo_{#1}^{#2}}
\newcommand{\objectiveibeam}[2]{\objectiveNo_{#1}(#2)}
\newcommand{\objectiveibeamMap}[1]{\objectiveNo_{#1}}

\newcommand{\x}{{\mathbf x}}
\newcommand{\dummyvec}{{\mathbf z}}
\newcommand{\y}{{\mathbf y}}
\newcommand{\snapshot}{{\mathbf s}}
\newcommand{\snapshotMat}{{\mathbf S}}
\newcommand{\weight}{\gamma}
\newcommand{\weightVec}{\boldsymbol \gamma}
\newcommand{\metric}{{ \boldsymbol \Theta}}
\newcommand{\xj}{\x_j}
\newcommand{\xjm}{\x_{j-1}}

\newcommand{\xjArg}[1]{\xj^{(#1)}}
\newcommand{\xArg}[1]{\x^{(#1)}}
\newcommand{\xProj}{\tilde \x}
\newcommand{\xArgHat}[1]{\hat \x^{(#1)}}
\newcommand{\xjk}{\xjArg{k}}

\newcommand{\xjkhat}{\hat\x_j^{(k)}}
\newcommand{\xhat}{\hat\x}
\newcommand{\kitNo}{k}
\newcommand{\kit}[1]{k_{#1}}
\newcommand{\xjkit}[1]{\xj^{(\kit{#1})}}
\newcommand{\memory}{\omega}
\newcommand{\bandwidth}{\tau}
\newcommand{\bandwidthj}{\bandwidth_j}
\newcommand{\xExact}{\x^\star}
\newcommand{\xExactArg}[1]{\xExact_{#1}}
\newcommand{\xExactj}{\xj^\star}
\newcommand{\xExactjm}{\xjm^\star}
\newcommand{\xExactjIn}{\xj^{\parallel}}
\newcommand{\xExactjOrth}{\xj^{\perp}}
\newcommand{\xExactjCOrth}{\xj^{\outputMat^\perp}}
\newcommand{\xExactjCIn}{\xj^{\outputMat}}
\newcommand{\xExactjCInWin}{\underline \xj^{\parallel}}
\newcommand{\xExactjCInWout}{\underline \xj^{\perp}}
\renewcommand{\b}{{\mathbf b}}
\newcommand{\bj}{\b_j}
\newcommand{\bjm}{\b_{j-1}}
\newcommand{\z}{{\mathbf z}}

\newcommand{\zArg}[1]{\z^{(#1)}}
\newcommand{\p}{{\mathbf p}}

\newcommand{\pArg}[1]{\p^{(#1)}}
\newcommand{\vecLength}{\gamma}
\newcommand{\vecLengthArg}[1]{\vecLength^{(#1)}}

\newcommand{\vecLengthAll}{{\bf \Gamma}}
\newcommand{\vecLengthAllj}{\vecLengthAll_j}
\newcommand{\alphaVec}{\hat {\bf v}}

\newcommand{\alphaArg}[1]{\alpha^{(#1)}}

\newcommand{\augComponent}{{\boldsymbol \mu}}

\newcommand{\augComponentArg}[1]{\augComponent^{(#1)}}
\renewcommand{\r}{{\mathbf r}}

\newcommand{\rArg}[1]{\r^{(#1)}}
\newcommand{\outputfunNo}{{\mathbf q}}
\newcommand{\outputfun}[1]{\outputfunNo\left(#1\right)}
\newcommand{\outputMat}{{\mathbf C}}
\newcommand{\Deltax}{\delta \mathbf x}

\newcommand{\Deltaxj}{\Deltax_j}
\newcommand{\DeltaxExactj}{\Deltaxj^\star}
\newcommand{\krylov}{{ K}}
\newcommand{\krylovjkNo}{\krylov^{(\kitNo)}}

\newcommand{\krylovjk}[2]{\krylovjkNo\left(#1,#2\right)}
\newcommand{\krylovj}[3]{\krylov^{(#3)}\left(#1,#2\right)}

\newcommand{\krylovjkUsed}{\krylovjkNo_j}
\newcommand{\krylovjkArbUsed}[1]{\krylov_j^{(#1)}}

\newcommand{\projOrth}[3]{\mathbf P^{#3}_{#1}\left({#2}\right)}
\newcommand{\projOrthNo}[2]{\mathbf P^{#2}_{#1}}
\newcommand{\projError}[3]{\left(\I-\mathbf P^{#3}_{#1}\right)\left({#2}\right)}
\newcommand{\projOrthAjmZjmSingVal}{\sigma_1}

\newcommand{\I}{\mathbf I}
\newcommand{\zero}{\mathbf 0}

\newcommand{\stageonebasis}{\mathbf W}
\newcommand{\stageonebasisArg}[1]{\stageonebasis_{#1}}
\newcommand{\stageonebasisj}{\stageonebasisArg{j}}
\newcommand{\rangestageonebasis}{\mathcal W}
\newcommand{\rangestageonebasisArg}[1]{\rangestageonebasis_{#1}}
\newcommand{\rangestageonebasisj}{\rangestageonebasisArg{j}}
\newcommand{\cholRed}{\hat {\bf R}}
\newcommand{\cholRedj}{\cholRed_j}
\newcommand{\innercholRedjArg}[1]{\bar {\bf R}_j^{(#1)}}
\newcommand{\innercholRedjk}{\innercholRedjArg{k}}
\newcommand{\chol}{{\bf R}}
\newcommand{\cholj}{\chol_j}
\newcommand{\cholEnforce}{{\bf L}^T}
\newcommand{\cholEnforceArg}[1]{\cholEnforce_{#1}}
\newcommand{\cholEnforcej}{\cholEnforceArg{j}}
\newcommand{\cholEnforcejinv}{{\bf L}^{-T}}
\newcommand{\cholEnforcejT}{{\bf L}_j}
\newcommand{\stageonedim}{w}
\newcommand{\stageonedimArg}[1]{\stageonedim_{#1}}
\newcommand{\stageonedimj}{\stageonedimArg{j}}
\newcommand{\stageonesolred}{\hat{\mathbf w}}
\newcommand{\stageonesolredj}{\stageonesolred_j}
\newcommand{\stageonesubspace}{{\mathcal W}}
\newcommand{\stageonesubspacej}{{\stageonesubspace_j}}

\newcommand{\stageoneaugj}{{\hat\stageonebasis_j}}
\newcommand{\rangestageoneaugj}{{\hat\stageonesubspacej}}

\newcommand{\stagetwobasis}{\mathbf X}
\newcommand{\stagetwobasisArg}[1]{\stagetwobasis_{#1}}
\newcommand{\stagetwobasisj}{\stagetwobasisArg{j}}
\newcommand{\rangestagetwobasis}{\mathcal X}
\newcommand{\rangestagetwobasisArg}[1]{\rangestagetwobasis_{#1}}
\newcommand{\rangestagetwobasisj}{\rangestagetwobasisArg{j}}
\newcommand{\stagetwosolred}{\hat{\mathbf x}}
\newcommand{\stagetwosolredj}{\stagetwosolred_j}
\newcommand{\innerstagetwosolredjArg}[1]{\stagetwosolred_j^{(#1)}}
\newcommand{\innerstagetwosolredjk}{\innerstagetwosolredjArg{k}}

\newcommand{\stagetwoaug}{{\hat\stagetwobasis}}
\newcommand{\stagetwoaugj}{\stagetwoaug_j}
\newcommand{\rangestagetwoaugj}{\hat\rangestagetwobasisj}
\newcommand{\innerstagetwoaugjArg}[1]{\bar \stagetwobasis_j^{(#1)}}
\newcommand{\innerstagetwoaugjk}{\innerstagetwoaugjArg{k}}
\newcommand{\stagetwoit}{\hat k}
\newcommand{\stagetwoitj}{\hat k_j}
\newcommand{\innerstagetwoitjArg}[1]{\bar k_j^{(#1)}}
\newcommand{\innerstagetwoitjk}{\innerstagetwoitjArg{k}}
\newcommand{\stagetwodim}{\stagetwoit}
\newcommand{\stagetwodimj}{\stagetwodim_j}

\newcommand{\augSpace}{\mathcal Y}
\newcommand{\augSpaceArg}[1]{\augSpace_{#1}}
\newcommand{\augSpacej}{\augSpaceArg{j}}
\newcommand{\augSpacejIdeal}{\augSpacej^\star}
\newcommand{\augSpacejIdealOut}{\bar \augSpace_j^\star}

\newcommand{\augSpaceWeightsVec}{\boldsymbol\eta}
\newcommand{\augSpaceWeightsVecPrev}{\augSpaceWeightsVec^\text{prev}_j}
\newcommand{\augSpaceWeightsEntryPrev}{[\augSpaceWeightsVecPrev]}
\newcommand{\augSpaceWeightsVecPrevArg}[1]{\augSpaceWeightsVec^\text{prev}_{#1}}
\newcommand{\augSpaceWeightsVecRBF}{\augSpaceWeightsVec^\text{RBF}_j}
\newcommand{\augSpaceWeightsVecRBFArg}[1]{\augSpaceWeightsVec^\text{RBF}_{#1}}
\newcommand{\augSpaceWeightsAjVec}{\boldsymbol\eta^{\star}_j}
\newcommand{\augSpaceWeightsAj}{[\augSpaceWeightsAjVec]}
\newcommand{\augSpaceWeightsAVecArg}[1]{\boldsymbol\eta^{\star}_{#1}}

\newcommand{\augSpaceWeightsCtC}{\bar{\boldsymbol\eta}^{\star}_j}
\newcommand{\augSpaceWeightsCtCVec}{\bar{\boldsymbol\eta}^{\star}_j}
\newcommand{\augSpaceBasis}{\mathbf Y}
\newcommand{\augSpaceBasisArg}[1]{\augSpaceBasis_{#1}}
\newcommand{\augSpaceBasisj}{\augSpaceBasisArg{j}}
\newcommand{\augSpaceBasisjp}{\augSpaceBasisArg{j+1}}
\newcommand{\augSpaceBasisjm}{\augSpaceBasisArg{j-1}}
\newcommand{\rangeaugSpaceBasis}{\mathcal Y}
\newcommand{\rangeaugSpaceBasisArg}[1]{\rangeaugSpaceBasis_{#1}}
\newcommand{\rangeaugSpaceBasisj}{\rangeaugSpaceBasisArg{j}}

\newcommand{\augSpaceSolredExactj}{\hat{\bm y}_j^\star}
\newcommand{\augSpaceSolredj}{\hat{\bm y}_j}

\newcommand{\currentData}{\mathbf Z}
\newcommand{\currentDataArg}[1]{\currentData_{#1}}
\newcommand{\currentDataj}{\mathbf Z_{j+1}}
\newcommand{\currentDatajm}{\mathbf Z_{j}}
\newcommand{\currentDatajmi}{\mathbf Z_{j+1-i}}

\newcommand{\rangecurrentDatajm}{\mathcal Z_{j}}

\newcommand{\currentDatavec}{\mathbf z}
\newcommand{\currentDataveci}[1]{\currentDatavec_{#1}}
\newcommand{\currentDataSpace}{\mathcal Z}
\newcommand{\currentDatajmSingVal}{\sigma_\text{min}}
\newcommand{\currentSpaceDim}{z}

\newcommand{\currentSpaceDimjm}{{\currentSpaceDim_{j}}}
\newcommand{\currentSpaceDimjmhalf}{{\currentSpaceDim_{j}^{1/2}}}
\newcommand{\currentDataSpaceArg}[1]{\currentDataSpace_{#1}}

\newcommand{\currentDataSpacejm}{\currentDataSpaceArg{j}}
\newcommand{\spaceBound}{{\alpha}}

\newcommand{\kryVec}{{\mathbf V}}
\newcommand{\rangekryVec}{{\mathcal V}}
\newcommand{\basis}{{\kryVec}}
\newcommand{\rangebasis}{{\mathcal V}}
\newcommand{\kryVecArg}[1]{\kryVec_{#1}}
\newcommand{\rangekryVecArg}[1]{\rangekryVec_{#1}}
\newcommand{\kryVecj}{\kryVecArg{j}}
\newcommand{\kryVecjm}{\kryVecArg{j-1}}
\newcommand{\rangekryVecj}{\rangekryVecArg{j}}
\newcommand{\kryVecjk}{\kryVecj^{(k)}}

\newcommand{\kryVecjkArg}[1]{\kryVecj^{(#1)}}

\newcommand{\kryVecjkT}{\kryVecj^{(k)T}}
\newcommand{\redKryj}{\hat{\mathbf v}_j}

\newcommand{\rj}{\r_j}

\newcommand{\xo}{\x^{(0)}}
\newcommand{\xjo}{\xo_j}
\newcommand{\xGuess}{\bar \x}
\newcommand{\xGuessArg}[1]{\xGuess_{#1}}
\newcommand{\xGuessj}{\xGuessArg{j}}
\newcommand{\xGuessjm}{\xGuessArg{j-1}}
\newcommand{\xjAugSpacej}{\xj^{\augSpace_j}}

\newcommand{\podSubspaceNo}{\mathcal U}
\newcommand{\podSubspace}[4]{\podSubspaceNo^{#4}_{#1}\left(#2,#3\right)}
\newcommand{\metricIdeal}{\metric^\ideal}
\newcommand{\metricPractical}{\metric^\text{\comp }}
\newcommand{\weightVecIdeal}{\weightVec^\ideal}
\newcommand{\weightEntryIdeal}{\weight^{\ideal }}
\newcommand{\weightVecPractical}{\weightVec^\text{\comp }}
\newcommand{\weightEntryPractical}{\weight^\text{\comp }}
\newcommand{\podSubspaceIdeal}{\podSubspaceNo^{\metricIdeal}_{\augSpaceDim}\left(\snapshotMat,\weightVecIdeal\right)}
\newcommand{\podSubspaceIdealShort}{\podSubspaceNo^{\ideal}}
\newcommand{\podSubspacePractical}{\podSubspaceNo^{\metricPractical}_{\augSpaceDim}\left(\snapshotMat,\weightVecPractical\right)}
\newcommand{\podSubspacePracticalShort}{\podSubspaceNo^{\comp}}
\newcommand{\podBasisNo}{{\mathbf U}}
\newcommand{\podBasis}[4]{\podBasisNo^{#4}_{#1}\left(#2,#3\right)}
\newcommand{\podBasisIdeal}{\podBasisNo^{\metricIdeal}_{\augSpaceDim}\left(\snapshotMat,\weightVecIdeal\right)}
\newcommand{\podBasisPractical}{\podBasisNo^{\metricPractical}_{\augSpaceDim}\left(\snapshotMat,\weightVecPractical\right)}
\newcommand{\podBasisSimple}[2]{\podBasisNo^{#2}_{#1}}

\newcommand{\feasibleSet}{A}
\newcommand{\energyCrit}{\varepsilon}
\newcommand{\energyCritStageTwo}{\energyCrit_y}
\newcommand{\energyCritStageOne}{\energyCrit_w}

\newcommand{\podVecNo}{\mathbf u}

\newcommand{\podVecSimple}[2]{\podVecNo^{#2}_{#1}}

\newcommand{\kronecker}[2]{\delta_{#1#2}}

\newcommand{\deflationEigenvectors}{\mathbf G}
\newcommand{\deflationEigenvector}{\mathbf g}
\newcommand{\deflationEigenvalues}{\boldsymbol \Lambda}
\newcommand{\deflationEigenvalue}{ \lambda}
\newcommand{\harmonicRitzVec}{\y}
\newcommand{\harmonicRitzVeci}{\harmonicRitzVec_i}
\newcommand{\harmonicRitzVecfirsti}{\bar\harmonicRitzVec_i}
\newcommand{\harmonicRitzValinv}{ \deflationEigenvalue}
\newcommand{\harmonicRitzValinvi}{\harmonicRitzValinv_i}
\newcommand{\harmonicRitzVal}{ \theta}
\newcommand{\harmonicRitzVali}{\harmonicRitzVal_i}
\newcommand{\harmonicRitzVals}{\boldsymbol\Lambda^{-1}}
\newcommand{\testHarmonicRitzVec}{\boldsymbol w}

\newcommand{\defeq}{\vcentcolon=}
\newcommand{\range}[1]{\text{range}(#1)}
\newcommand{\myspan}[1]{\text{span}\{#1\}}
\newcommand{\T}{\boldsymbol\Gamma}
\newcommand{\TArg}[1]{\T_{#1}}
\newcommand{\Tj}{\TArg{j}}
\newcommand{\redTj}{\hat{\boldsymbol\Gamma}_j}
\newcommand{\innerredTjArg}[1]{\bar{\boldsymbol\Gamma}_j^{(#1)}}
\newcommand{\innerredTjk}{\innerredTjArg{k}}

\newcommand{\tol}{\epsilon}
\newcommand{\tolstagetwo}{\hat \tol}
\newcommand{\tolstagetwoArg}[1]{\tolstagetwo_{#1}}
\newcommand{\tolstagetwoj}{\hat \tol_j}
\newcommand{\innertolstagetwo}{\bar \tol}
\newcommand{\innertolstagetwoArg}[1]{\innertolstagetwo_{#1}}
\newcommand{\innertolstagetwoj}{\innertolstagetwoArg{j}}
\newcommand{\tolj}{\tol_j}
\newcommand{\tolOutput}{\tau}
\newcommand{\nlhs}{p}
\newcommand{\ndof}{n}
\newcommand{\nout}{q}
\newcommand{\augSpaceDim}{y}
\newcommand{\augSpaceDimj}{{\augSpaceDim_j}}
\newcommand{\augSpaceDimjp}{{\augSpaceDim_{j+1}}}
\newcommand{\augSpaceDimjm}{{\augSpaceDim_{j-1}}}

\newcommand{\stageoneset}{\mathcal K}
\newcommand{\stageonesetj}{\stageoneset_j}
\newcommand{\nsnapshots}{s}
\newcommand{\spd}[1]{\text{SPD}(#1)}
\newcommand{\spsd}[1]{\text{SPSD}(#1)}

\renewcommand{\diag}{\text{diag}}

\newcommand{\fext}{\ensuremath{\boldsymbol f^\text{ext}}}
\newcommand{\fintNo}{\ensuremath{\boldsymbol f^\text{int}}}
\newcommand{\fint}[1]{\ensuremath{\fintNo\left(#1\right)}}
\newcommand{\fcontactNo}{\ensuremath{\boldsymbol f^\text{contact}}}
\newcommand{\fcontact}[2]{\ensuremath{\fcontactNo\left(#1,#2\right)}}
\newcommand{\penaltyNo}{\ensuremath{\lambda}}
\newcommand{\penaltyArg}[1]{\penaltyNo^{#1}}
\newcommand{\penaltyl}{\penaltyArg{l}}
\newcommand{\displac}{\ensuremath{\boldsymbol u}}
\newcommand{\displacDummy}{\bar\displac}
\newcommand{\displaci}{\displac_i}
\newcommand{\displacilk}{\displac_i^{l(k)}}
\newcommand{\displacik}{\displac_i^{(k)}}
\newcommand{\timei}{t_i}
\newcommand{\rbfWeight}[2]{\rho^{\text{RBF}}(#1,#2)}

\newcommand{\rbfIDWWeightNo}{\rho^{\text{IDW}}}
\newcommand{\rbfIDWWeight}[1]{\rbfIDWWeightNo(#1)}

\newcommand{\ntimeSteps}{T}
\newcommand{\npenalty}{L}
\newcommand{\temperature}{\ensuremath{\Delta T}}
\newcommand{\strainThermal}{\ensuremath{\epsilon_\text{thermal}}}

\newcommand{\vectomat}[2]{\left[#1_1\ \cdots\ #1_{#2}\right]}
\newcommand{\totuple}[2]{\left(#1_1, \ldots, #1_{#2}\right)}

\newenvironment{proofofthing}[1][Proof]{\begin{trivlist}
\item[\hskip \labelsep {\bfseries Proof of #1}]}{$\square$ \end{trivlist}}

\usepackage{verbatim}
\newcommand{\stageonealgorithmname}{\texttt{direct\_reduced\_solve}}
\newcommand{\stagetwoalgorithmname}{\texttt{augmented\_pcg}}
\newcommand{\overallalgorithmname}{\texttt{three\_stage\_algorithm}}
\newcommand{\compressionalgorithmname}{\texttt{compression}}
\newcommand{\podsvdalgorithmname}{\texttt{pod\_svd}}
\newcommand{\podevdalgorithmname}{\texttt{pod\_evd}}
\newcommand{\stagethreeoptionone}{\varphi}
\newcommand{\true}{1\ }
\newcommand{\newvectorsdirect}{\varrho}
\newcommand{\augSpaceThreshold}{\bar y}

\newcommand{\jcompress}{\bar j}

\newcommand{\stiefel}[2]{\mathcal S_{#1}({\RR{#2}})}
\newcommand{\distanceMetric}[2]{d\left(#1,#2\right)}

\newcommand{\myA}[1]{\A_{#1}}
\newcommand{\myZ}[1]{\snapshotMat}
\newcommand{\myIdealY}[1]{\podBasisNo^\ideal}
\newcommand{\myPracticalY}[1]{\podBasisNo^{\comp }}
\newcommand{\myPracticalMetric}{\metric^\comp }
\newcommand{\myIdealMetric}{\metric^\ideal}
\newcommand{\myMetricPerturb}{\boldsymbol \Delta} 
\newcommand{\myIdealX}[1]{\bar{\mathbf V}^\ideal}
\newcommand{\myPracticalX}[1]{\bar{\mathbf V}^{\comp }}
\newcommand{\myIdealTruncatedX}[1]{{\mathbf V}^\ideal}
\newcommand{\myPracticalTruncatedX}[1]{{\mathbf V}^{\comp }}

\newcommand{\myIdealEtaMatrix}[1]{\boldsymbol \Gamma^\ideal}
	\newcommand{\myPracticalEtaMatrix}[1]{\boldsymbol \Gamma^{\comp }}
\newcommand{\myIdealEta}[1]{\weightVec^\ideal}
\newcommand{\myPracticalEta}[1]{\weightVec^{\comp }_{#1}}
\newcommand{\myIdealLambda}[1]{\bar{\boldsymbol \Lambda}^\ideal}
\newcommand{\myPracticalLambda}[1]{\bar{\boldsymbol \Lambda}^{\comp }}
\newcommand{\myIdealTruncatedLambda}[1]{{\boldsymbol \Lambda}^\ideal}
\newcommand{\myPracticalTruncatedLambda}[1]{{\boldsymbol \Lambda}^\comp}
\newcommand{\myGeneralPerpIdealTruncatedLambda}{\boldsymbol\Lambda^\ideal_{B,\perp}}
\newcommand{\myPerpIdealTruncatedLambda}[1]{\boldsymbol\Lambda^\ideal_\perp}
\newcommand{\myGeneralPracticalTruncatedLambda}{\boldsymbol\Lambda_B^{\comp }}
\newcommand{\E}{\mathbf E} 
\newcommand{\myGeneralPractical}{\mathbf B^{\comp }}
\newcommand{\myGeneralPracticalSubspace}{\mathcal S^{\comp }}
\newcommand{\myGeneralIdeal}{\mathbf B^\ideal}
\newcommand{\myGeneralIdealSubspace}{\mathcal S^\ideal}

\newcommand{\myRR}[1]{\mathbb{R}^{#1}}
\newcommand{\myTruncDimension}{y}
\newcommand{\mySnapDimension}{s}

\newcommand{\myRelativeWeightBound}{\hat{\delta}}
\newcommand{\myIdealCondNumber}{\kappa^\ideal}
\newcommand{\myPracticalCondNumber}{\kappa^{\comp }}
\newcommand{\myCombined}[1]{\mathbf D}
\newcommand{\myCC}{\outputMat^T \outputMat}

\newcommand{\myZComplement}{\mathbf Z_\perp}

\usepackage[algo2e,linesnumbered,ruled]{algorithm2e}
\DontPrintSemicolon
\SetKwComment{tcm}{ }{}
\title{Krylov-subspace recycling via\\ the POD-augmented conjugate-gradient method}
\author{Kevin Carlberg\thanks{Sandia National Laboratories
(\href{mailto:ktcarlb@sandia.gov}{ktcarlb@sandia.gov},
\href{mailto:rstumin@sandia.gov}{rstumin@sandia.gov}).}
\and Virginia Forstall\thanks{University of Maryland (\href{mailto:vhfors@gmail.com}{vhfors@gmail.com}).}
\and Ray Tuminaro\footnotemark[1]
}%

\begin{document}
\maketitle

\begin{abstract} 
This work presents a new Krylov-subspace-recycling method for efficiently solving sequences of linear systems of equations characterized by varying right-hand sides and symmetric-positive-definite matrices. As opposed to typical truncation strategies used in recycling such as deflation, we propose a truncation method inspired by goal-oriented proper orthogonal decomposition (POD) from model reduction.  This idea is based on the observation that model reduction aims to compute a low-dimensional subspace that contains an \textit{accurate} solution; as such, we expect the proposed method to generate a low-dimensional subspace that is well suited for computing solutions that can satisfy \textit{inexact} tolerances. In particular, we propose specific goal-oriented POD `ingredients' that align the optimality properties of POD with the objective of Krylov-subspace recycling.  To compute solutions in the resulting `augmented' POD subspace, we propose a hybrid direct/iterative three-stage method that leverages 1) the optimal ordering of POD basis vectors, and 2) well-conditioned reduced matrices. Numerical experiments performed on solid-mechanics problems highlight the benefits of the proposed method over existing approaches for Krylov-subspace recycling.
\end{abstract}


\section{Introduction}

This work considers solving a sequence of linear systems of equations
characterized by varying right-hand sides and symmetric-positive-definite
matrices. Such problems arise in a variety of engineering and science
applications, including structural optimization, nonlinear structural
dynamics, unconstrained numerical optimization, and nonlinear
electromagnetics.  In particular, we consider solving these linear systems to
\textit{inexact tolerances} using the (preconditioned) conjugate-gradient
method; further, we allow for a solution-dependent output to serve as a
quantity of interest.

While each linear-system solve can be executed independently of previous
solves, reusing data generated during these solves can lead to improved
convergence; this observation has led to the emergence of
\textit{Krylov-subspace recycling} methods.  These recycling methods can also
be considered `augmented Krylov subspace methods' \cite{saad1997aak} because
they `augment' the typical Krylov subspace with a subspace computed from
previous data, and subsequently compute solutions via projection with this augmented
subspace. Initially, researchers developed methods that employ the space spanned
by \textit{all} Krylov vectors generated during the solution of previous
linear systems as the augmenting subspace.
In this framework, researchers developed 
methods to treat the case of an invariant matrix and multiple right-hand sides
available either simultaneously (i.e., block Krylov methods) \cite{oleary1980bcg}
or sequentially
\cite{saad1987lms,farhatmultiRHS,erhel2000acg}; the
sequential case occurs, for example, when restarting
Krylov-subspace methods.  These ideas were also extended to solve multiple
systems with varying matrices; approximate orthogonalization techniques
\cite{rey1994dda,roux1995}, projection methods
\cite{farhat1995optimizing,farhat2000sgf}, and an efficient full
orthogonalization method \cite{risler2000iaa} were developed for this purpose.

However, retaining the accumulation of all previous Krylov vectors can be
computationally expensive and memory intensive, particularly when convergence
is slow, the number of previous linear systems is large, or the length of the Krylov vectors is not large (e.g., preconditioners
based on domain decomposition).  This has led to the development of
\textit{truncation methods} that retain a subspace of the span of previous Krylov vectors as
the augmenting subspace. First, deflation techniques for sequences
of systems with invariant \cite{chapman1997daa,saad2000deflated} and varying
\cite{rey1998rrp,parks2007rks} matrices were developed.  These methods employ
approximated eigenvectors associated with the smallest eigenvalues of the
governing matrices as the augmenting subspace. As such, they are
effective primarily in cases where convergence is hampered by a few small
eigenvalues.  An alternative approach computes the augmenting subspace
as the subspace that most accurately represents the Krylov
subspace in the orthogonalization step of the generalized conjugate residual
(GCR) method; this was also developed for both the invariant-matrix
\cite{desturler1999truncation} and varying-matrices \cite{parks2007rks} cases.

These truncated Krylov-subspace-recycling techniques do not target the
efficient solution of \textit{inexact solutions}, which is the focus of this
work. To this end, we propose a new proper orthogonal decomposition
(POD)-augmented conjugate-gradient method
for Krylov-subspace recycling. This approach---which 
employs goal-oriented model reduction to truncate previous Krylov vectors---is inspired by
the observation that model-reduction techniques aim to generate
low-dimensional approximations that preserve high levels of accuracy (i.e.,
satisfy inexact tolerances); it is based on preliminary work presented in
Ref.~\cite{carlberg-adaptive}.  The paper consists of the following new
contributions.
\begin{itemize} 
  \item We perform analyses that expose the close relationship between goal-oriented POD
	and Krylov-subspace recycling (Theorems \ref{Anorm}--\ref{weightsDifference} and Corollary
	\ref{cor:computableWeightsIdeal}).
	\item 
		We propose goal-oriented POD ingredients for truncating previous Krylov vectors,
		including
    \begin{itemize}
		\item \textit{Snapshots} comprising all previous Krylov vectors, 
		\item \textit{Metrics} induced by 1) the system matrix and 2) the
		output quantity of interest, and
		\item \textit{Snapshot weights} arising from 1) the linear system before
		truncation and 2) a radial-basis-function approximation of the solution.
    \end{itemize}
Further, we bound the distance between these goal-oriented POD subspaces and
the corresponding `ideal' subspaces (Theorem \ref{thm:podDistanceGen}).
  \item
		We develop a novel `three-stage' algorithm, which accelerates the solution over the
		augmenting subspace using a hybrid direct/iterative approach. The algorithm leverages
		1) the optimal ordering of POD basis vectors and 2) well-conditioned
		reduced matrices. The algorithm comprises
		\begin{itemize} 
		\item \textit{Stage 1}: Direct solution over the first few (high-energy) POD
		basis vectors,
		\item \textit{Stage 2}: Iterative solution over the full augmenting space
		using the augmented CG algorithm, and
		\item \textit{Stage 3}: Iterative solution over the full space
		using the augmented preconditioned CG algorithm. This stage is equipped
		with new strategies for efficiently orthogonalizing against the entire
		augmenting subspace.
		\end{itemize}
\end{itemize}

Section \ref{sec:prob} provides the problem formulation, Section
\ref{sec:podaugcg} describes the proposed POD-augmented conjugate-gradient
algorithm, Section \ref{sec:threestage} describes the three-stage algorithm,
Section \ref{sec:experiments} provides numerical experiments, and Section
\ref{sec:conclusions} concludes the paper. Note that Appendix
\ref{app:proofs} provides proofs for theoretical results.
In the remainder of this manuscript, we
denote matrices by capitalized bold letters, vectors by lowercase bold
letters, and scalars by lowercase letters.  We denote the columns of a matrix
$\bm{A}\in\RR{m\times n}$ by $\bm{a}_i\in\RR{m}$, $i\in\nat{n}$ with
$\nat{a}\defeq\{1,\ldots, a\}$ such that $\bm{A}\defeq\left[\bm{a}_1\ \cdots\
\bm{a}_{n}\right]$. We denote the scalar-valued matrix elements by
$a_{ij}\in\RR{}$ such that  $\bm{a}_j\equiv\left[a_{1j}\ \cdots\
a_{mj}\right]^T$, $j\in\nat{n}$.  In addition, we denote the range of a matrix
by its calligraphic counterpart, i.e.,  $\range{\A} \equiv \mathcal A$; we
sometimes refer to $\A$ as the `basis' for $\mathcal A$, although---more
precisely---it is the `basis in matrix form' for $\mathcal A$.

\section{Problem formulation}\label{sec:prob}
This work considers solving a sequence of linear systems with a varying
matrix
 \begin{equation} \label{eq:sequence}
 \Aj\xExactj = \bj, \quad j\innatseq{\nlhs}.
  \end{equation} 
Here, $\Aj\in\spd{\ndof}$ and $\bj\in\RR{\ndof}$ denote
the $j$th sparse system matrix and right-hand side,
respectively, with $\spd{\ndof}$ denoting the set of
symmetric-positive-definite (SPD) $\ndof\times\ndof$ matrices. The quantity
$\xExactj\in\RR{\ndof}$ is implicitly defined as the
(exact) solutions to Eqs.~\eqref{eq:sequence}. 
Further, we assume that the primary objective 
is to compute an output quantity of interest $\outputfun{\xj}$,
$j\innatseq{\nlhs}$ with
\begin{align} \label{eq:outputs}
\outputfunNo:&\RR{\ndof}\rightarrow\RR{\nout}\\
& \x\mapsto\outputMat\x,
\end{align} 
where 
$\outputMat\in\RR{\nout\times\ndof}$.
We consider computing a sequence of \textit{inexact} solutions $\xj$, $j\innatseq{\nlhs}$ to
Eqs.~\eqref{eq:sequence} that satisfy
 \begin{equation} \label{eq:sequenceInexact}
 \|\bj-\Aj\xj\|_2\leq \tolj ,  \quad j\innatseq{\nlhs}.
  \end{equation} 
	where $\tolj\geq 0$, $j\innatseq{\nlhs}$ denotes the forcing sequence
	\cite{eisenstat1996choosing}. 

\subsection{Conjugate-gradient method}

To compute each inexact solution $\xj$, we consider
	applying the preconditioned conjugate gradient (PCG) algorithm, which
	computes a sequence of solutions that minimize the energy norm of the
	error over the (current) Krylov subspace. For the $j$th linear system, these solutions satisfy
\begin{equation}\label{eq:KrylovSolve} 
\xjk= \arg\min_{\x\in\xjo + \krylovjkUsed(\xjo)}\normAj{\xExactj - \x},\quad
k\innatseq{\kit{j}},
\end{equation} 
where
$\krylovjkUsed:\x\mapsto
\krylovjk{\precInvj\Aj}{\precInvj\left(\bj - \Aj\x\right)}\subseteq\RR{\ndof}$ is the affine
search subspace,
$\krylovjkNo:(\A,\b)\mapsto\myspan{\b,\A\b,\ldots,\A^{k-1}\b}\subseteq\RR{\ndof}$ is the
Krylov subspace at iteration $k$, 
$\precj\in\spd{\ndof}$ is a preconditioner, 
 $\xjo\in\RR{\ndof}$ is the initial approximation,
$\innerProdAj{\x}{\y}\defeq \x^T\Aj\y$ and 
$\normAj{\x}\defeq\sqrt{\innerProdAj{\x}{\x}}$ denote the $\Aj$-weighted inner
product and norm, and $\kit{j}$ denotes the number of
iterations required such that $\xj \defeq \xjkit{j}$ satisfies 
inequality \eqref{eq:sequenceInexact}.

Optimality property \eqref{eq:KrylovSolve} can be interpreted in a number of
ways. In particular, given a basis $\basis\in\RRstar{\ndof\times\augSpaceDim}$
and centering point $\xGuess\in\RR{\ndof}$ for an affine subspace 
$\xGuess + {\rangebasis}$ (recall $\rangebasis \equiv \range{\basis}$), then
the solution $\xProj$ that minimizes the $\A$-norm of the error 
with respect to the exact solution $\xExact$---which is implicitly defined by $\A\xExact = \b$
with $\A\in\spd{\ndof}$ and $\b\in\RR{\ndof}$---over this subspace can be written as
\begin{equation}\label{eq:KrylovSolveGen} 
\xProj = \arg\min_{\x\in\xGuess + \rangebasis}\normA{\xExact - \x}.
\end{equation} 
Here, $\RRstar{m \times n}$ denotes the
noncompact Stiefel manifold: the set of full-column-rank $m\times n$ matrices.
By the definition of an orthogonal projection \cite{saad2003iterative}, 
Eq.~\eqref{eq:KrylovSolveGen} can be equivalently expressed as 
\begin{gather}\label{eq:KrylovSolveProjGen} 
\xProj = \projOrth{\A}{\xExact}{\xGuess +\rangebasis},
\end{gather} 
where $\projOrth{\A}{\x}{\krylov}$ denotes the $\A$-orthogonal projection of
the vector $\x\in\RR{\ndof}$ onto the (possibly affine) subspace $\krylov\subseteq\RR{\ndof}$, i.e.,
$\projOrth{\A}{\x}{\krylov}\in \krylov,\ \forall\x\in\RR{\ndof}$; $\left(\x -
\projOrth{\A}{\x}{\krylov},\dummyvec\right)_\A= 0,\ \forall\x\in\RR{\ndof},\
\forall \dummyvec\in\krylov$; and
$\projOrth{\A}{\projOrth{\A}{\x}{\krylov}}{\krylov} =
\projOrth{\A}{\x}{\krylov}$ (idempotency).
This orthogonal projector can be defined algebraically as
\begin{equation} \label{eq:algebraicProjectorGen}
\projOrthNo{\A}{\xGuess + {\rangebasis}}:\x\mapsto \xGuess +
\basis\left(\basis^T\A\basis\right)^{-1}\basis^T\A(\x- \xGuess).
\end{equation} 
If a symmetric factorization $\A = (\A^{1/2})^T\A^{1/2}$ is available (where $\A^{1/2}$
need not be upper triangular), then an equivalent definition is
\begin{equation} \label{eq:algebraicProjectorPseudoGen}
\projOrthNo{\A}{\xGuess + {\rangebasis}}:\x\mapsto \xGuess +
\basis\left(\A^{1/2}\basis\right)^{+}\A^{1/2}(\x- \xGuess).
\end{equation} 
where a superscript $+$ denotes the Moore--Penrose pseudoinverse.

Substituting definition \eqref{eq:algebraicProjectorGen} in
Eq.~\eqref{eq:KrylovSolveProjGen} reveals that the solution $\xProj$ can be
computed via Galerkin
projection
\begin{gather}\label{eq:KrylovSolveGalGen} 
\basis^T\A\basis\xhat =
\basis^T\r(\xGuess),\quad
 \xProj = \xGuess + \basis\xhat,
\end{gather}
where $\r:\x\mapsto \b-\A\x$ denotes the residual
operator and we have used $\A\xExact = \b$.
If $\basis^T\A\basis$ is invertible, then Eqs.\ \eqref{eq:KrylovSolveProjGen}
and \eqref{eq:KrylovSolveGalGen} imply that 
a vector $\xProj$ is the $\A$-orthogonal projection of
the exact solution $\xExact$ to $\A\xExact = \b$ onto an affine subspace $\xGuess + {\rangebasis}$ if and only if its associated
residual is orthogonal to $\rangebasis$ in the Euclidean inner product, i.e.,
\begin{equation}\label{eq:equivalenceResidualOrth}
\x = \projOrth{\A}{\xExact}{\xGuess + {\rangebasis}}\Leftrightarrow
\basis^T\r(
\x
) = 0.
\end{equation}

Thus, equivalent interpretations of approximate solutions $\xjk$ are 1)
as the solution to minimization problem \eqref{eq:KrylovSolve}; 2)
as the orthogonal projection (from Eq.~\eqref{eq:KrylovSolveGen})
\begin{gather}
\xjk = \projOrth{\Aj}{\xExactj}{\xjo +\krylovjkUsed(\xjo)},\quad
k\innatseq{\kit{j}};
\end{gather}
3) as the Galerkin solution (from Eq.~\eqref{eq:KrylovSolveGalGen})
\begin{gather}
\kryVecjkT\Aj\kryVecjk\xjkhat =
\kryVecjkT\rj(\xjo),
\\
\label{eq:KrylovSolveGal2} \xjk = \xjo + \kryVecjk\xjkhat,\quad
k\innatseq{\kit{j}},
\end{gather} 
where $\kryVecjk\in\RRstar{\ndof\times k}$ constitutes a basis for the
subspace $\krylovjkUsed(\xjo)$ and 
$\rj:\x\mapsto\bj-\Aj\x$ denotes the residual for linear system $j$;
and 
4)
as the solution yielding residual orthogonality
\begin{equation} 
\kryVecjkT\rj( \xjk) = 0.
\end{equation}


\subsection{Augmented conjugate-gradient method}\label{sec:augCG}

Krylov-subspace recycling aims to reduce the computational burden of solving
linear system $j$ of 
Eqs.~\eqref{eq:sequence} 
by exploiting a previously computed `augmenting' subspace
$\augSpacej\subseteq\RR{\ndof}$ of dimension $\augSpaceDimj\leq \ndof$ spanned by a basis
$\augSpaceBasisj\in\RRstar{\ndof\times \augSpaceDimj}$. 
To achieve this, recycling strategies employ augmented Krylov-subspace
methods as a primary tool; rather than perform the typical sequence of Krylov iterations
(e.g., Eqs.~\eqref{eq:KrylovSolve}), these methods first compute a solution in
the augmenting subspace, and then compute an increment in a newly generated
Krylov subspace. Critically, these methods ensure that the final solution
minimizes the error over the sum of augmenting and Krylov subspaces; this
generally requires maintaining orthogonality of new Krylov vectors (or
$\Aj$-orthogonality of new search directions) to the augmenting subspace.

In the context of the conjugate-gradient method, these methods first solve a minimization
problem in the affine subspace $\xGuessj + \augSpacej$, where $\xGuessj$ is an
initial guess:
\begin{equation}\label{eq:KrylovSolveAug} 
\xjAugSpacej= \arg\min_{\x\in\xGuessj + \augSpacej}\normAj{\xExactj - \x} = \projOrth{\Aj}{\xExactj}{\xGuessj+\augSpacej}. 
\end{equation} 
As before, 
the solution can be computed via Galerkin projection
\begin{gather}
\label{eq:augGal}\augSpaceBasisj^T\Aj\augSpaceBasisj\augSpaceSolredj =
\augSpaceBasisj^T\left(\bj-\Aj\xGuessj\right),\quad \xjAugSpacej=\xGuessj + \augSpaceBasisj\augSpaceSolredj.
\end{gather}
Subsequently, these methods solve a minimization problem in the `augmented'
Krylov subspace 
\begin{equation}\label{eq:KrylovSolveAugment} 
\xjk= \arg\min_{\x\in\xGuessj + \augSpacej + \krylovjkUsed(\xjAugSpacej)}\normAj{\xExactj -
\x},\quad k\innatseq{\kit{j}}.
\end{equation} 
The final solution can then be expressed as 
\begin{equation}\label{eq:xjAug} 
\kryVecj^T\Aj\kryVecj\redKryj =
\kryVecj^T\left(\bj-\Aj\xGuessj\right),\quad 
\xj=\xGuessj + \augSpaceBasisj\augSpaceSolredj +
\kryVecj\redKryj, 
\end{equation}
where $\kryVecj\defeq\kryVecjkArg{\kit{j}}$ denotes a basis that satisfies
 \begin{gather} 
\label{eq:Ajorthogonality}\kryVecj^T\Aj\kryVecj = \vecLengthAllj,
\quad\augSpaceBasisj^T\Aj\kryVecj=\zero,\\
\label{eq:directSum}{\rangeaugSpaceBasisj} \oplus
{\rangekryVecj} = \augSpacej
+\krylovjkUsed(\xjAugSpacej),
  \end{gather} 
	where $\vecLengthAllj$ is diagonal.
Critically, $\Aj$-orthogonality  of the augmented CG search directions to the
augmenting basis \eqref{eq:Ajorthogonality} ensures the `direct
sum' property \eqref{eq:directSum}. This allows the 
solution to be
expressed as the (unique) sum of components in $\rangeaugSpaceBasisj$ and
$\rangekryVecj$ by decoupling
the Galerkin
problems \eqref{eq:augGal} and \eqref{eq:xjAug}.
The solution increment $\kryVecj\redKryj$ in Eq.~\eqref{eq:xjAug} can be
computed via the augmented PCG algorithm
\ref{alg:augmentedCGgen}
as 
$$
(\kit{j}, \redKryj, \kryVecj,
\vecLengthAllj) = \stagetwoalgorithmname(
\Aj,  \bj - \Aj\xGuessj,
\augSpaceSolredj,
\augSpaceBasisj, \precj, \tolj).
$$
Here, steps \ref{step:aug2}--\ref{step:aug2part2}, and
\ref{step:aug3}--\ref{step:aug3part2} enforce $\Aj$-orthogonality property
\eqref{eq:Ajorthogonality}.

\begin{algorithm}[htbp]
\caption{ \stagetwoalgorithmname }
\begin{algorithmic}[1]\label{alg:augmentedCGgen}
\small
\REQUIRE $\A$, $\b$, $\xArgHat{0}$, $\augSpaceBasis$, $\prec$, 
 $\tol$
\ENSURE $\kitNo$, $\alphaVec$, $\kryVec$, $\vecLengthAll$  
\STATE $\xArg{0} = \augSpaceBasis\xArgHat{0}$
\STATE $\rArg{0} = \b - \A\xArg{0}$
\STATE $\zArg{0} = \precInv \rArg{0}$
\STATE \label{step:aug2}Solve $\augSpaceBasis^T \A   \augSpaceBasis\augComponentArg{0} =
\augSpaceBasis^T \A  \zArg{0}$. 
\STATE\label{step:aug2part2} $\pArg{0} = \zArg{0} - \augSpaceBasis\augComponentArg{0} $
\FOR{$k = 0, 1, \ldots$}
\STATE $\vecLengthArg{k} = (\A\pArg{k},\pArg{k})$
\STATE  $\alphaArg{k} = (\rArg{k},\zArg{k})/\vecLengthArg{k}$
\STATE  $\xArg{k+1} = \xArg{k} + \alphaArg{k} \pArg{k} $
\STATE\label{step:res}  $\rArg{k+1} = \rArg{k} - \alphaArg{k} \A  \pArg{k}$
\STATE  $\zArg{k+1} = \precInv\rArg{k+1}$
\STATE  $\beta^{(k+1)} = \frac{(\rArg{k+1},\zArg{k+1})}{(\rArg{k},\zArg{k})}$
\STATE  \label{step:aug3}Solve $\augSpaceBasis^T \A  \augSpaceBasis
\augComponentArg{k+1} =
\augSpaceBasis^T \A  \zArg{k+1} $.
\STATE\label{step:aug3part2}  $ \pArg{k+1} = \zArg{k+1} + \beta^{(k+1)} \pArg{k} -
\augSpaceBasis\augComponentArg{k+1}$
\IF{$\|\rArg{k+1}\|\leq \tol$  } 
\STATE Exit.
\ENDIF
 \ENDFOR
 \STATE $\kitNo = k+1$, $\alphaVec =
 \left[\alphaArg{0}\ \cdots\ \alphaArg{k-1}\right]^T$, $\kryVec =
  \left[\pArg{0}\ \cdots\
 \pArg{k-1}\right]$, $\vecLengthAll =
 \diag{(\vecLengthArg{0},\ldots,\vecLengthArg{k-1})}$ 
\end{algorithmic}
\end{algorithm}

Several strategies exist for selecting the augmenting basis $\augSpaceBasisj$. Typically,
the columns of this matrix consist of \textit{all} Krylov vectors generated
during the solution of
previous linear systems; this provides the interpretation of `recycling'
Krylov subspaces. 
In this case, we have
\begin{equation} 
\augSpaceBasisj = \vectomat{\kryVec}{j-1}.
\end{equation} 
However, after a modest number
of linear systems has been solved, it becomes memory- and
computation-intensive to retain and orthogonalize against this complete set of
Krylov vectors.  Therefore, truncation techniques have been devised to 
retain only the subspace that is
`most important' in some sense. In particular, after solving linear system $j-1$, these methods compute
$\augSpaceBasisArg{j}$ such that
$\augSpaceArg{j}\subseteq\currentDataSpacejm$, where 
$\currentDatajm \defeq
\left[\augSpaceBasisjm,\ \kryVecjm\right]\in\RR{\ndof\times \currentSpaceDimjm}$ denotes the matrix of all preserved
vectors accumulated over the first $j-1$ linear solves,
with
$\currentSpaceDimjm = \augSpaceDimjm + \kit{j-1}$. Note that
$\currentDatajm = \vectomat{\kryVec}{j-1}$ before truncation first
occurs.

\subsection{Deflation with harmonic Ritz vectors}

Deflation techniques aim to retain the subspace associated with eigenvectors
that tend to hamper convergence, i.e., those with eigenvalues close to zero.
Such techniques have been developed for multiple linear systems with an
invariant (e.g., restarting) \cite{chapman1997daa,saad2000deflated} and
varying \cite{rey1998rrp,parks2007rks} matrix.

To accomplish this, these techniques compute the harmonic Ritz vectors---which
approximate the eigenvectors of $\Aj$ with the smallest eigenvalues---by
solving the following problem: Find
$\harmonicRitzVecfirsti\in\range{\Ajm\currentDatajm}$ and
$\harmonicRitzValinvi\in\RR{}$, $i = 1,\ldots,\augSpaceDimj$ with $\augSpaceDimj \leq
\currentSpaceDimjm$ such that 
\begin{equation} 
\left(\testHarmonicRitzVec,\Ajm^{-1}\harmonicRitzVecfirsti -
\harmonicRitzVecfirsti\harmonicRitzValinvi\right) =
0,\quad\forall\testHarmonicRitzVec\in\range{\Ajm\currentDatajm}.
\end{equation} 
An equivalent problem statement is the following: Find
$\harmonicRitzVeci\in{\rangecurrentDatajm}$ and
$\harmonicRitzVali\in\RR{}$, $i = 1,\ldots,\augSpaceDimj$ with $\augSpaceDimj \leq
\currentSpaceDimjm$ such that 
\begin{equation} 
\left(\testHarmonicRitzVec,\Ajm\harmonicRitzVeci -
\harmonicRitzVeci\harmonicRitzVali\right) =
0,\quad\forall\testHarmonicRitzVec\in\range{\Ajm\currentDatajm},
\end{equation} 
where $\harmonicRitzVali=\harmonicRitzValinvi^{-1}$ and
$\harmonicRitzVecfirsti = \Ajm\harmonicRitzVeci$,
$i\innatseq{\augSpaceDimjp}$. Because vector $\harmonicRitzVeci$ is
equivalent to eigenvector $\harmonicRitzVecfirsti$ with one step of inverse
iteration, vectors $\harmonicRitzVeci$ are typically employed as eigenvector
approximations
\cite{morgan1991computing}. 
	Algebraically, this corresponds to solving the generalized eigenvalue
	problem
\begin{equation}\label{eq:deflationEVP} 
\currentDatajm^T\Ajm^T\Ajm\currentDatajm\deflationEigenvectors=\currentDatajm^T\Ajm\currentDatajm\deflationEigenvectors\harmonicRitzVals
\end{equation} 
with
$\deflationEigenvalues=\diag\totuple{\harmonicRitzValinv}{\currentSpaceDimjm}$ and subsequently setting
$\augSpaceBasisj \leftarrow
\vectomat{\currentDatajm\deflationEigenvector}{\augSpaceDimj}$, where the
eigenvalues $\harmonicRitzValinvi$ are ordered in decreasing magnitude.
 
 Deflation has been pursued in the context of GMRES with deflated
 restarting (GMRES-DR \cite{chapman1997daa,rey1998rrp,risler2000iaa,morgan2003gmres}), GCR with orthogonalization and deflated
 restarting (GCRO-DR) \cite{parks2007rks}, and the deflated conjugate gradient method
 \cite{chapman1997daa}.
This approach is effective primarily in cases where the matrix is characterized by
a small number of eigenvalues close to zero that hamper convergence. Further,
because this approach aims to promote convergence to the `exact' solution,
it is not tailored for the efficient computation of \textit{inexact} solutions, which
is the focus of this work. To this end, we propose a novel truncation strategy
inspired by model reduction.

\section{POD-augmented CG}\label{sec:podaugcg}


Because this work focuses on efficiently computing \textit{inexact} solutions,
we aim to compute an `ideal' low-dimensional subspace $\augSpacejIdeal$
that directly minimizes the error in the augmenting-subspace
solution appearing in Eq.~\eqref{eq:KrylovSolveAug}, i.e.,
\begin{equation} \label{eq:grassmanOpt}
\augSpacejIdeal =
\arg\min_{\augSpace\in\grassman{\augSpaceDimj}{\ndof}}\normAj{\xExactj- \projOrth{\Aj}{\xExactj}{\xGuessj+\augSpacej}},
\end{equation} 
where $\grassman{m}{n}$ with $m\leq n$  denotes the set of $m$-dimensional subspaces of
$\RR{n}$  (the Grassmannian). 
Alternatively, we may be interested in computing inexact solutions that most accurately
represent output quantities of interest. In this case, we can also consider an
ideal output-oriented subspace
\begin{equation} \label{eq:grassmanOptOut}
\augSpacejIdealOut =
\arg\min_{\augSpace\in\grassman{\augSpaceDimj}{\ndof}}
\|\outputfun{\xExactj}-
\outputfun{\projOrth{\Aj}{\xExactj}{\xGuessj+\augSpacej}}\|_2.
\end{equation}

Clearly, these ideal subspaces are computable if
the exact solution $\xExactj$ is known; in this case, we can enforce
$\xExactj-\xGuessj\in\augSpacejIdeal$ (resp.\
$\xExactj-\xGuessj\in\augSpacejIdealOut$)---which yields a
zero objective-function value in Eq.~\eqref{eq:grassmanOpt} (resp.\
Eq.~\eqref{eq:grassmanOptOut})---for any $\augSpaceDimj\geq 1$; if $\augSpaceDimj=1$, then
$\augSpacejIdeal=\myspan{\xExactj-\xGuessj}$. However, the exact solution $\xExactj$ is
not known before the $j$th linear system is solved. In this case, we aim to
compute a subspace that \textit{approximately} solves minimization
problem \eqref{eq:grassmanOpt} (resp.\ \eqref{eq:grassmanOptOut}).
For this purpose, we employ goal-oriented POD with carefully chosen snapshots,
weights, and metrics
\cite{CarlbergCPOD,carlbergCpodJour}.

\subsection{POD}

The POD method \cite{POD,rathinam:newlook} generates a basis that optimally represents a set of vectors
(or `snapshots') in a certain sense. The technique was developed in the
context of model reduction, and has been applied to reduced-order modeling 
in structural dynamics \cite{lall2003structure,camReview,carlberg2012spd} and fluid
dynamics
\cite{sirovich1987tad3,sirisup2004spectral,noack2005need,wang2012proper,san2013proper,carlbergJCP} among others; it is also closely
related to the singular value decomposition, principal component analysis, and
the Karhunen--Lo\`{e}ve expansion. We aim to select goal-oriented POD
ingredients (i.e., snapshots, weights, and metric) \cite{CarlbergCPOD,carlbergCpodJour}
that align the POD optimality property with the objective function in
Eq.~\eqref{eq:grassmanOpt}.

Given a matrix $\snapshotMat\in\RR{\ndof\times \nsnapshots}$ whose columns
represent $\nsnapshots$ snapshots,
a vector of weights
$\weightVec\in\RR{\nsnapshots}$, and pseudometric
associated with matrix $\metric\in\spsd{\ndof}$, where $\spsd{\ndof}$ denotes
the set of $\ndof\times\ndof$ symmetric-positive-semidefinite matrices, the
POD method computes a subspace of dimension 
$\augSpaceDim\leq\nsnapshots$ that minimizes the sum of squared projection errors,
i.e.,
\begin{equation} \label{eq:podOpt}
\podSubspace{\augSpaceDim}{\snapshotMat}{\weightVec}{\metric}
= \arg\min_{\mathcal
A\in\grassman{\augSpaceDim}{\ndof}}\sum_{i=1}^\nsnapshots\|\projError{\metric}{\weight_i\snapshot_i}{\mathcal
A}\|_{\metric}^2,\quad \augSpaceDim\innatseq{\nsnapshots}.
\end{equation} 
A basis
$\podBasis{\augSpaceDim}{\snapshotMat}{\weightVec}{\metric}\in\RRstar{\ndof\times\augSpaceDim}$
for this POD subspace can be computed by solving the eigenvalue problem
\begin{equation}
\diag\totuple{\weight}{\nsnapshots}\snapshotMat^T\metric\snapshotMat\diag\totuple{\weight}{\nsnapshots}\mathbf V=\mathbf V\boldsymbol \Sigma^2
\end{equation} 
with $\boldsymbol\Sigma\defeq\diag\totuple{\sigma}{\nsnapshots}$ and $\mathbf
V\in\stiefel{\nsnapshots}{\nsnapshots}$  and setting 
\begin{equation}\label{eq:podAlgebraicDefinition}
\podBasis{\augSpaceDim}{\snapshotMat}{\weightVec}{\metric}
=\snapshotMat\diag\totuple{\weight}{\nsnapshots}\left[\frac{1}{\sigma_1}\mathbf v_1\ \ldots\
\frac{1}{\sigma_{\augSpaceDim}}\mathbf v_{\augSpaceDim}\right],
\end{equation}
where the eigenvalues $\sigma^2_i$ are ordered in decreasing magnitude.
Here, $\stiefel{m}{n}$ denotes the
Stiefel manifold: the set of orthonormal $m$-frames in $\RR{n}$.
Alternatively, the POD basis can be computed via the singular value decomposition;
Appendix \ref{app:podcomputation} describes both algorithms.
The resulting basis is nested
 \begin{equation} \label{eq:podordering}
\podBasisSimple{i+1}{\metric} = 
[\podBasisSimple{i}{\metric}\ 
\podVecSimple{i+1}{\metric}] 
,\quad i\innatseq{\nsnapshots-1}
  \end{equation}
with $\podVecSimple{i}{\metric}\in\RR{\ndof}$,
$i\innatseq{s}$ and $\podBasisSimple{1}{\metric} = 
\podVecSimple{1}{\metric}$ and exhibits $\metric$-orthogonality, i.e., 
\begin{equation}\label{eq:podOrth}
\innerProd{\podVecSimple{i}{\metric}}
{\podVecSimple{j}{\metric}}{\metric}
= \kronecker{i}{j},
\end{equation}
where $\kronecker{i}{j}$ denotes the Kronecker delta. We note that the
algebraic definition \eqref{eq:algebraicProjectorGen} is not valid when
considering pseudometrics, as $\basis^T\metric\basis$ may not be invertible if
$\metric$ is semidefinite. In this case, definition
\eqref{eq:algebraicProjectorPseudoGen} is appropriate, as $\metric^{1/2}\basis$
always has a pseudoinverse.
\begin{remarkk}\label{rem:podvecorder}
Note that Eqs.~\eqref{eq:podOpt}--\eqref{eq:podordering} imply that the POD
basis vectors are optimally ordered, i.e., the first $\augSpaceDim$ POD basis
vectors span the optimal $\augSpaceDim$-dimensional subspace in the sense of
Eq.~\eqref{eq:podOpt}.
\end{remarkk}

\subsection{POD and the augmented conjugate
gradient method}

We propose employing goal-oriented POD to define the augmenting subspace, i.e.,
\begin{equation}\label{eq:augSpaceIsPOD}
\augSpacej =
\podSubspace{\augSpaceDim}{\snapshotMat}{\weightVec}{\metric}
\end{equation}
for specific choices of the snapshots $\snapshotMat$, weights $\weightVec$,
	and metric $\metric$. The following
	results provide guidance toward this end. We first show that POD is related
	to minimizing $\normAj{\xExactj- \projOrth{\Aj}{\xExactj}{\xGuessj+\augSpacej}}$, which is the error
	minimized by ideal subspace $\augSpacejIdeal$ in problem \eqref{eq:grassmanOpt}.
\begin{theorem}\label{Anorm}
The POD subspace 
$\podSubspace{\augSpaceDim}{\currentDatajm}{\augSpaceWeightsAjVec}{\Aj}$ with
 \begin{equation}\label{eq:weightsIdealAj}
\augSpaceWeightsAjVec \defeq
\left(\currentDatajm^T\Aj\currentDatajm\right)^{-1}\currentDatajm^T\left(\bj-\Aj\xGuessj\right)
  \end{equation} 
minimizes an
upper bound for $\normAj{\xExactj- \projOrth{\Aj}{\xExactj}{\xGuessj+\augSpacej}}$ over all $\augSpaceDim$-dimensional
subspaces of $\currentDataSpacejm\subseteq\RR{\ndof}$, where $\augSpaceDim\leq\currentSpaceDimjm$.
\end{theorem}
See Appendix \ref{app:proofs} for the proof (and all subsequent proofs).

We now show that other POD ingredients can be selected to align POD 
	with minimizing $\|\outputfun{\xExactj}-
\outputfun{\projOrth{\Aj}{\xExactj}{\xGuessj+\augSpacej}}\|_2$,
	which is the error
	minimized by ideal subspace $\augSpacejIdealOut$ in problem
	\eqref{eq:grassmanOptOut}.
\begin{theorem}\label{outNorm}
The POD subspace 
$\podSubspace{\augSpaceDim}{\currentDatajm}{\augSpaceWeightsCtCVec}{\outputMat^T\outputMat}$
 with
 \begin{equation} \label{eq:weightsIdealOut}
\augSpaceWeightsCtCVec \defeq
 \left(\outputMat\currentDatajm\right)^{+}\outputMat
 \left(\xExactj-\xGuessj\right)
  \end{equation} 
minimizes an
upper bound for the output error $\|\outputfun{\xExactj}-
\outputfun{\projOrth{\Aj}{\xExactj}{\xGuessj+\augSpacej}}\|_2$ over all $\augSpaceDim$-dimensional
subspaces of $\currentDataSpacejm\subseteq\RR{\ndof}$, where $\augSpaceDim\leq\currentSpaceDimjm\leq \ndof$.

\end{theorem}

\subsection{Goal-oriented POD ingredients}\label{sec:podingredients}

In light of these theoretical results, we now propose several
practical choices for POD ingredients that align goal-oriented POD with augmented CG.

\subsubsection{Snapshots}
Theorems \ref{Anorm} and \ref{outNorm} show that POD minimizes an
upper bound for errors of interest if snapshots are set to vectors accumulated
over the first $j-1$ linear solves.  We therefore employ $\snapshotMat =
\currentDatajm$ in Eq.~\eqref{eq:augSpaceIsPOD}. 

\subsubsection{Weights}\label{sec:podweights}

In practice, the augmenting subspace $\augSpacej$ is computed after linear system
$j-1$ is solved. As such, we cannot compute the `ideal weights'
$\augSpaceWeightsCtCVec$ defined in \eqref{eq:weightsIdealOut}, which requires
knowledge of $\xExactj$. While $\augSpaceWeightsAjVec$ in
Eq.~\eqref{eq:weightsIdealAj} can be computed by solving
$\currentDatajm^T\Aj\currentDatajm\augSpaceWeightsAjVec =
\currentDatajm^T\left(\bj-\Aj\xGuessj\right)$, this is not practical: the
computational cost of doing so is equivalent to solving Eq.~\eqref{eq:augGal}
with $\augSpaceBasisj  = \currentDatajm$, i.e., employing an augmenting subspace
that has not
been truncated. As we aim to employ a truncated augmented space for the $j$th
linear-system solve, we consider two approximations to these
ideal weights.
\begin{enumerate} 
\item \textbf{Previous weights}. This approach employs 
weights of
 \begin{equation} \label{eq:weightsPrevious}
\augSpaceWeightsVecPrev =
\left(\currentDatajm^T\Ajm\currentDatajm\right)^{-1}\currentDatajm^T\left(\bjm-\Ajm\xGuessjm\right).
  \end{equation} 
	Comparing Eqs.~\eqref{eq:weightsIdealAj} and \eqref{eq:weightsPrevious}
	reveals that $\augSpaceWeightsVecPrev$ is `close' to
	$\augSpaceWeightsAjVec$, but employs readily available data, as these
	weights are equal to the coefficient in the expansion of the solution at the
	previous time step, i.e., $\xjm = \xGuessjm + \currentDatajm
	\augSpaceWeightsVecPrev$. We now provide bounds for the difference between
	these computable weights and the ideal weights.
\begin{theorem}\label{weightsDifference}
The difference between the previous weights and ideal weights can be bounded
as
\begin{equation} \label{eq:weightsDiffAj}
\|\augSpaceWeightsAjVec - \augSpaceWeightsVecPrev\| \leq\frac{1}{\currentDatajmSingVal}\left(
\|\projOrthNo{\Aj}{\currentDataSpacejm} -
\projOrthNo{\Ajm}{\currentDataSpacejm}\|\|{\xExactj - \xGuessj}\| +
\projOrthAjmZjmSingVal\| (\xExactj - \xGuessj) - (\xExactjm - \xGuessjm)\|
\right)
\end{equation} 
and
\begin{equation} \label{eq:weightsDiffCtC}
\|\augSpaceWeightsCtCVec - \augSpaceWeightsVecPrev\| \leq\frac{1}{\currentDatajmSingVal}\left(
\|\projOrthNo{\outputMat^T\outputMat}{\currentDataSpacejm} -
\projOrthNo{\Ajm}{\currentDataSpacejm}\|\|{\xExactj - \xGuessj}\| +
\projOrthAjmZjmSingVal\| (\xExactj - \xGuessj) - (\xExactjm - \xGuessjm)\|
\right),
\end{equation} 
where $\currentDatajmSingVal$ denotes the smallest singular value of
$\currentDatajm$ and $\projOrthAjmZjmSingVal$ denotes the largest singular
value of $\projOrthNo{\Ajm}{\currentDataSpacejm}$.
\end{theorem}
We now provide conditions under which the computable weights are equal to the
ideal weights.
\begin{corollary} \label{cor:computableWeightsIdeal}
If $\Ajm = \Aj$ and $\xExactj-\xGuessj = \xExactjm-\xGuessjm$, then
$\augSpaceWeightsAjVec = \augSpaceWeightsVecPrev$. Alternatively, if
$\outputMat^T\outputMat = \Aj$
and $\xExactj-\xGuessj = \xExactjm-\xGuessjm$, then $\augSpaceWeightsCtCVec -
\augSpaceWeightsVecPrev$.
\begin{proof}
The result follows trivially from Eqs.~\eqref{eq:weightsDiffAj} and \eqref{eq:weightsDiffCtC}.
\end{proof}
	\end{corollary}

\item \textbf{Radial-basis-function weights}. Noting that
$\currentDatajm\augSpaceWeightsAjVec = \projOrth{\Aj}{\xExactj -
\xGuessj}{\currentDataSpacejm}$, we can approximate ideal weights
$\augSpaceWeightsAjVec$ by approximating the component of $\xExactj -
\xGuessj$ in ${\rangecurrentDatajm}$. To achieve this, we assume the $j$th
solution can be approximated as a linear combination of previous solutions,
i.e., 
\begin{equation}
\projOrth{\Aj}{\xExactj - \xGuessj}{\currentDataSpacejm} \approx
\sum_{i=1}^\memory\rbfWeight{j}{j-i}(\projOrth{\A_{j-i}}{\xExactArg{j-i} -
\xGuessArg{j-i}}{\currentDataSpaceArg{j-i}}),
\end{equation}
where $\memory\innat{j-1}$ denotes the number of previous solutions
to include and $\rbfWeight{i}{j}$ denotes a radial basis function. This
implies weights of
 \begin{align} \label{eq:weightsRBF}
 \begin{split}
\augSpaceWeightsVecRBF(\memory) &=
\sum_{i=1}^\memory\rbfWeight{j}{j-i}
\left(\currentDatajmi^T\Ajmi\currentDatajmi\right)^{-1}\currentDatajmi^T\Ajmi\left(\xExactArg{j-i}-\xGuessArg{j-i}\right)\\
&=
\sum_{i=1}^\memory\rbfWeight{j}{j-i}
\augSpaceWeightsVecPrevArg{j+1-i}
.
 \end{split}
  \end{align} 
In practice, we employ a inverse-distance-weight radial basis function of the form
$\rbfWeight{i}{j}= \rbfIDWWeight{|i-j|}$ with $\rbfIDWWeightNo:r\mapsto
1/2^{r-1}$; we set $\memory$ to be the number of linear systems since the
most recent truncation.
\end{enumerate}


\subsubsection{Metric}\label{sec:podmetric}
Theorems \ref{Anorm} and \ref{outNorm} demonstrate that POD minimizes
an upper bound for the errors $\normAj{\xExactj- \projOrth{\Aj}{\xExactj}{\xGuessj+\augSpacej}}$ and
$\|\outputfun{\xExactj}- \outputfun{\projOrth{\Aj}{\xExactj}{\xGuessj+\augSpacej}}\|_2$ only if POD metrics of
$\metric = \Aj$ and $\metric = \outputMat^T\outputMat$ are used, respectively.
However, as previously discussed, we aim to avoid using $\Aj$ to truncate the
basis $\currentDatajm$, as it entails reduced computations with the $j$th
linear system (i.e., solving Eqs.~\eqref{eq:augGal} with $\augSpaceBasisj  =
\currentDatajm$). Therefore, we employ two practical choices for metrics in
Eq.~\eqref{eq:augSpaceIsPOD}:
\begin{enumerate} 
\item $\metric = \Ajm$. This approach aligns the truncation with
minimizing $\normAj{\xExactj- \projOrth{\Aj}{\xExactj}{\xGuessj+\augSpacej}}$. From
Eq.~\eqref{eq:podOrth}, the resulting basis is
$\Ajm$-orthogonal. Here, computing the POD basis via Algorithm \ref{alg:PODEVD} is
appropriate, as a symmetric factorization of $\Ajm$ is not readily available.
\item $\metric = \outputMat^T\outputMat$. This is an output-oriented approach
associated with minimizing the output error. From Eq.~\eqref{eq:podOrth}, the
resulting basis is $\outputMat^T\outputMat$-orthogonal. Computing this POD
basis via Algorithm \ref{alg:PODSVD} is appropriate, as a symmetric
factor of the (pseudo)metric is readily available as $\outputMat$.
\end{enumerate}
\subsubsection{Computable goal-oriented POD methods}\label{sec:podingredientssummary}

\noindent In summary, the four (computable) goal-oriented POD truncation methods we propose are
 \begin{enumerate} 
  \item \label{option:compressAjprev}
	$
	\augSpaceBasisj  =
	\podBasis{\augSpaceDimj}{\currentDatajm}{\augSpaceWeightsVecPrev}{\Ajm}
	$, computable by 
$(\augSpaceBasisj) = \podevdalgorithmname(\currentDatajm,
\augSpaceWeightsVecPrev, \Ajm,\energyCritStageTwo)$,
\item\label{option:compressAjRBF} 
	$
	\augSpaceBasisj  =
	\podBasis{\augSpaceDimj}{\currentDatajm}{\augSpaceWeightsVecRBF(\memory)}{\Ajm}
	$, computable by 
$(\augSpaceBasisj) = \podevdalgorithmname(\currentDatajm,
\augSpaceWeightsVecRBF(\memory), \Ajm,\energyCritStageTwo)$,
\item \label{option:compressCprev}
	$
	\augSpaceBasisj  =
	\podBasis{\augSpaceDimj}{\currentDatajm}{\augSpaceWeightsVecPrev}{\outputMat^T\outputMat}
	$, computable by 
$(\augSpaceBasisj) = \podsvdalgorithmname(\currentDatajm,
\augSpaceWeightsVecPrev, \outputMat,\energyCritStageTwo)$,
\item \label{option:compressCRBF}
	$
	\augSpaceBasisj  =
	\podBasis{\augSpaceDimj}{\currentDatajm}{\augSpaceWeightsVecRBF(\memory)}{\outputMat^T\outputMat}
	$, computable by 
$(\augSpaceBasisj) = \podsvdalgorithmname(\currentDatajm,
\augSpaceWeightsVecRBF(\memory), \outputMat,\energyCritStageTwo)$,
	 \end{enumerate}
	 where $\energyCritStageTwo\in\left[0,1\right]$ is a statistical `energy
	 criterion' used to truncate the POD basis, and \podevdalgorithmname\ and
	 \podsvdalgorithmname\ are described in Algorithms \ref{alg:PODEVD} and
	 \ref{alg:PODSVD}, respectively.
We now provide a result that allows us to bound the distance between the
four computable goal-oriented POD subspaces above and their ideal counterparts,
which correspond to $\podSubspace{\augSpaceDimj}{\currentDatajm}{\augSpaceWeightsAjVec}{\Aj}$ for
\ref{option:compressAjprev}--\ref{option:compressAjRBF} and 
$\podSubspace{\augSpaceDimj}{\currentDatajm}{\augSpaceWeightsCtCVec}{\outputMat^T\outputMat}$
for
\ref{option:compressCprev}--\ref{option:compressCRBF} above.

\begin{theorem}[Distance between goal-oriented POD subspaces]\label{thm:podDistanceGen}
The distance between two goal-oriented POD subspaces
$\podSubspaceIdealShort\defeq\podSubspaceIdeal$ and
$\podSubspacePracticalShort\defeq\podSubspacePractical$---which are characterized by different
metrics 
$\myIdealMetric\in\spsd{\ndof}$ 
and 
$\myPracticalMetric\in\spsd{\ndof}$ 
with $\myIdealMetric
= \myPracticalMetric + \myMetricPerturb$
and different snapshot weights 
$\myIdealEtaMatrix{j}\defeq\diag\totuple{\weightEntryIdeal}{\nsnapshots}$
and
$\myPracticalEtaMatrix{j}\defeq\diag\totuple{\weightEntryPractical}{\nsnapshots}$ 
associated with snapshots
$\myZ{j-1}\in\RR{\ndof\times\nsnapshots}$---can be bounded as
\begin{align} \label{eq:subspaceDistanceBound}
 \distanceMetric{\podSubspaceIdealShort}{\podSubspacePracticalShort}
\le\kappa(\myZ{j-1} \myIdealEtaMatrix{j}) \kappa([\myIdealEtaMatrix{j}]^{-1} \myPracticalEtaMatrix{j}
\myPracticalTruncatedX{j}) \Bigl(&\| 
 \left(\myPracticalEtaMatrix{j} +
\myIdealEtaMatrix{j}\right)\left(\myPracticalEtaMatrix{j} -
\myIdealEtaMatrix{j}\right)[\myIdealEtaMatrix{j}]^{-1} \myZ{j-1}^T
\myPracticalMetric \myZ{j-1} \myIdealEtaMatrix{j} \|_2 \\
& + \| \myMetricPerturb
\|_2 \| \myZ{j-1} \myIdealEtaMatrix{j} \|^2_2\Bigr)/
\text{abssep}(\myPerpIdealTruncatedLambda{j},\myPracticalTruncatedLambda{j})
\end{align}
if the matrix $\myZ{j-1}
\myIdealEtaMatrix{j}$ has full column rank.
Here, $\distanceMetric{\mathcal U}{\mathcal V}\defeq \max_{u\in\mathcal
U, \|u\|=1}{\min_{v\in\mathcal V}}\|u - v\| = \sin(\theta_\text{max})$, where
$\theta_\text{max}$ is the largest principal angle between $\mathcal U$ and $\mathcal
V$, $\kappa(\A)$ denotes the condition number of matrix $\A$, and $
\text{abssep}(\boldsymbol\Lambda_1,\boldsymbol\Lambda_2)\defeq
\min_{\|\boldsymbol Z\|_2 = 1}\|\boldsymbol \Lambda_1 \boldsymbol Z -
\boldsymbol Z\boldsymbol\Lambda_2\|_2
$ denotes the absolute separation between two spectra.
In addition, 
$\myPracticalTruncatedX{j}\in\stiefel{\myTruncDimension}{\mySnapDimension}$
and $\myPracticalTruncatedLambda{j}\in\RR{\myTruncDimension\times\myTruncDimension}$
denote the first $\augSpaceDim$ eigenvectors and (diagonal matrix of) eigenvalues, respectively,
of the matrix
$\myPracticalEtaMatrix{j} \myZ{j-1}^T \myPracticalMetric \myZ{j-1}
\myPracticalEtaMatrix{j}$, assuming eigenvalues are ordered in descending
magnitude.
Similarly, $\myPerpIdealTruncatedLambda{j}$ denotes
the diagonal matrix of eigenvalues associated with 
the orthogonal complement to the
invariant subspace spanned by the first $\augSpaceDim$ eigenvectors of $\myIdealEtaMatrix{j} \myZ{j-1}^T  
\myIdealMetric \myZ{j-1} \myIdealEtaMatrix{j}$.
\end{theorem}

This theorem allows us to bound the distances between the four computable POD
subspaces and their ideal counterparts by substituting the appropriate
quantities into Theorem \ref{thm:podDistanceGen}.
That is, the theorem allows us to bound the following quantities:
\begin{enumerate} 
\item\label{distanceone}
$\distanceMetric{\podSubspace{\augSpaceDim}{\currentDatajm}{\augSpaceWeightsVecPrev}{\Ajm}}{\podSubspace{\augSpaceDim}{\currentDatajm}{\augSpaceWeightsAjVec}{\Aj}}$,
\item\label{distancetwo}
$\distanceMetric{\podSubspace{\augSpaceDim}{\currentDatajm}{\augSpaceWeightsVecRBF(\memory)}{\Ajm}}{\podSubspace{\augSpaceDim}{\currentDatajm}{\augSpaceWeightsAjVec}{\Aj}}$,
\item\label{distancethree}
$\distanceMetric{\podSubspace{\augSpaceDim}{\currentDatajm}{\augSpaceWeightsVecPrev}{\outputMat^T\outputMat}}{\podSubspace{\augSpaceDim}{\currentDatajm}{\augSpaceWeightsAjVec}{\outputMat^T\outputMat}}$,
and
\item\label{distancefour}
$\distanceMetric{\podSubspace{\augSpaceDim}{\currentDatajm}{\augSpaceWeightsVecRBF(\memory)}{\outputMat^T\outputMat}}{\podSubspace{\augSpaceDim}{\currentDatajm}{\augSpaceWeightsAjVec}{\outputMat^T\outputMat}}$.
\end{enumerate}
\begin{remarkk}
\label{rem:subspaceDistance}
The subspace-distance bound \eqref{eq:subspaceDistanceBound}
simplifies under certain conditions. \\
\textit{Fixed snapshot weights.}
First, note that when the snapshot weights are identical but the metric
differes, i.e., 
$\myPracticalEtaMatrix{j} = \myIdealEtaMatrix{j}$, the first term within parentheses vanishes and 
$\kappa([\myIdealEtaMatrix{j}]^{-1} \myPracticalEtaMatrix{j}
\myPracticalTruncatedX{j})=1
$ due to symmetry of $\myGeneralPractical$, which yields
\begin{align} \label{eq:subspaceDistanceBoundMetricOnly}
 \distanceMetric{\podSubspaceIdealShort}{\podSubspacePracticalShort}
\le\kappa(\myZ{j-1} \myIdealEtaMatrix{j}) 
 \| \myMetricPerturb
\|_2 \| \myZ{j-1} \myIdealEtaMatrix{j} \|^2_2/
\text{abssep}(\myPerpIdealTruncatedLambda{j},\myPracticalTruncatedLambda{j}).
\end{align}
Thus, the subspace distance depends linearly on the norm of the metric
perturbation in this case.

\noindent \textit{Fixed metric.} If instead the snapshot weights differ while $\myMetricPerturb = {\bf 0}$, then 
the second term within parentheses vanishes and we have
\begin{align} \label{eq:subspaceDistanceBoundWeightsOnly}
 \distanceMetric{\podSubspaceIdealShort}{\podSubspacePracticalShort}
\le\kappa(\myZ{j-1} \myIdealEtaMatrix{j}) \kappa([\myIdealEtaMatrix{j}]^{-1} \myPracticalEtaMatrix{j}
\myPracticalTruncatedX{j}) &\| 
\left(\myPracticalEtaMatrix{j} +
\myIdealEtaMatrix{j}\right)\left(\myPracticalEtaMatrix{j} -
\myIdealEtaMatrix{j}\right)[\myIdealEtaMatrix{j}]^{-1} \myZ{j-1}^T
\myPracticalMetric \myZ{j-1} \myIdealEtaMatrix{j} \|_2 
/
\text{abssep}(\myPerpIdealTruncatedLambda{j},\myPracticalTruncatedLambda{j})
\end{align}
In practice, this situation corresponds to cases
\ref{distancethree}--\ref{distancefour} above (where 
$\myPracticalMetric =\myIdealMetric= \myCC$) or to cases
\ref{distanceone}--\ref{distancetwo} when 
$\myA{j} = \myA{j-1}$.
Notice that the relative weight perturbation $(\myPracticalEtaMatrix{j} -
\myIdealEtaMatrix{j})[\myIdealEtaMatrix{j}]^{-1}$
appears within the remaining term. When the relative weight perturbation can
be bounded as
$$
\max_k{(\myPracticalEta{j} - \myIdealEta{j})_k / (\myIdealEta{j})_k} \le \myRelativeWeightBound ,
$$
then it can be shown that the bound simplifies to 
\begin{align} \label{eq:subspaceDistanceBoundWeightsOnlyRelFix}
 \distanceMetric{\podSubspaceIdealShort}{\podSubspacePracticalShort}
\le\kappa(\myZ{j-1} \myIdealEtaMatrix{j}) \kappa([\myIdealEtaMatrix{j}]^{-1} \myPracticalEtaMatrix{j}
\myPracticalTruncatedX{j}) &\myRelativeWeightBound (2 +
\myRelativeWeightBound) \| \myIdealEtaMatrix{j} \myZ{j-1}^T \myIdealMetric
\myZ{j-1} \myIdealEtaMatrix{j} \|_2
/
\text{abssep}(\myPerpIdealTruncatedLambda{j},\myPracticalTruncatedLambda{j}).
\end{align}
Further, when 
$\myIdealEtaMatrix{j}$ commutes with $\myZ{j-1}^T \myPracticalMetric \myZ{j-1}$, then the bound simplifies to 
\begin{align} \label{eq:subspaceDistanceBoundWeightsOnlyCommutes}
 \distanceMetric{\podSubspaceIdealShort}{\podSubspacePracticalShort}
\le\kappa(\myZ{j-1} \myIdealEtaMatrix{j}) \kappa([\myIdealEtaMatrix{j}]^{-1} \myPracticalEtaMatrix{j}
\myPracticalTruncatedX{j}) &\| 
\left(\myPracticalEtaMatrix{j} +
\myIdealEtaMatrix{j}\right)\left(\myPracticalEtaMatrix{j} -
\myIdealEtaMatrix{j}\right)  \myZ{j-1}^T
\myPracticalMetric \myZ{j-1}  \|_2 
/
\text{abssep}(\myPerpIdealTruncatedLambda{j},\myPracticalTruncatedLambda{j}),
\end{align}
which now depends on the absolute weight perturbation,
$\myPracticalEtaMatrix{j} - \myIdealEtaMatrix{j}$, rather than
a relative perturbation. This occurs for cases
\ref{distanceone}--\ref{distancetwo} above when $\A_1=\cdots=\Aj=\A$
(i.e., an invariant matrix),
as  $\myZ{j-1}^T \myPracticalMetric \myZ{j-1}$ is diagonal due to enforced $\A$-orthogonality of $\myZ{j-1}$.

\noindent \textit{Strongly separated eigenvalues.} If the eigenvalues associated with 
$\myPracticalTruncatedLambda{j}$ and $\myPerpIdealTruncatedLambda{j}$ in Theorem \ref{thm:podDistanceGen} are strongly separated, i.e., 
$$0<\delta_a\defeq\max\{\min_{\lambda\in\myPracticalTruncatedLambda{j}}|\lambda|-\max_{\hat\lambda\in\myPerpIdealTruncatedLambda{j}}|\hat
\lambda|, \min_{\hat\lambda \in\myPerpIdealTruncatedLambda{j}}|\hat
\lambda| - \max_{\lambda\in\myPracticalTruncatedLambda{j}}|\lambda|\}$$
then the expression for the subspace-distance bound becomes \cite[Theorem 5.3]{ipsen2000absolute}
\begin{align} \label{eq:subspaceDistanceBoundNoAbsSep}
 \distanceMetric{\podSubspaceIdealShort}{\podSubspacePracticalShort}
\le\kappa(\myZ{j-1} \myIdealEtaMatrix{j}) \kappa([\myIdealEtaMatrix{j}]^{-1} \myPracticalEtaMatrix{j}
\myPracticalTruncatedX{j}) \Bigl(&\| 
 \left(\myPracticalEtaMatrix{j} +
\myIdealEtaMatrix{j}\right)\left(\myPracticalEtaMatrix{j} -
\myIdealEtaMatrix{j}\right)[\myIdealEtaMatrix{j}]^{-1} \myZ{j-1}^T
\myPracticalMetric \myZ{j-1} \myIdealEtaMatrix{j} \|_2 \\
& + \| \myMetricPerturb
\|_2 \| \myZ{j-1} \myIdealEtaMatrix{j} \|^2_2\Bigr)/
\delta_a.
\end{align}

Finally, we note that we can apply the results in Eq.~\eqref{eq:weightsDiffAj}
from Theorem \ref{weightsDifference} to those of
Theorem \ref{thm:podDistanceGen} via the triangle inequality to obtain a more
explicit subspace-distance bound for case \ref{distanceone} above, where 
$\myIdealEtaMatrix{j} = \diag\totuple{[\augSpaceWeightsAjVec]}{\nsnapshots}$
and $\myPracticalEtaMatrix{j} = 
\diag\totuple{\augSpaceWeightsEntryPrev}{\nsnapshots}$.
\end{remarkk}

\section{Three-stage algorithm}\label{sec:threestage}

From Section \ref{sec:augCG}, we know that augmented PCG relies on first computing a solution 
over the augmenting subspace
$\xGuessj+\augSpacej$. While this can be accomplished by solving the Galerkin
system~\eqref{eq:augGal} directly, this requires assembling
the reduced matrix $\augSpaceBasisj^T\Aj\augSpaceBasisj$, which incurs  $\augSpaceDimj$
matrix--vector products and 
$\mathcal
	O((\bandwidthj + \augSpaceDimj)\augSpaceDimj\ndof)$ flops, where 
	$\bandwidthj$
	denotes the average number of nonzeros per row of $\Aj$
	This cost can be significant when the dimension of the augmenting
	space $\augSpaceDimj$ becomes large.

To mitigate this effect, we propose a novel three-stage algorithm for augmented
PCG. The algorithm leverages an augmenting basis that is optimally ordered and
yields a well-conditioned reduced matrix. In particular, for linear system $j$,
we assume that an augmenting basis
$\augSpaceBasisj\in\RRstar{\ndof\times\augSpaceDimj}$ is available that 1)
contains a low-dimensional basis
$\stageonebasisj\in\RRstar{\ndof\times\stageonedimj}$ with
${\rangestageonebasisj}\subseteq {\rangeaugSpaceBasisj}\subseteq \RR{\ndof}$
and $\stageonedimj\leq\augSpaceDimj\leq \ndof$ that can capture an accurate
solution, and 2) yields a well-conditioned reduced matrix
$\augSpaceBasisj^T\Aj\augSpaceBasisj$. Given this basis, stage 1 first 
computes an accurate solution over the subspace
$\xGuessj+\rangestageonebasisj$; a direct approach is suitable due to its assumed
low dimensionality. Stage 2 iteratively computes a
solution over $\xGuessj+\augSpacej$. This should be
efficient, as $\augSpaceBasisj^T\Aj\augSpaceBasisj$ need never be formed:
the well-conditioned assumption implies that a preconditioner is not required
for fast convergence and matrix--vector products of the form
	$\augSpaceBasis^T\A\augSpaceBasis\p$ can be computed using a single
	matrix--vector product as $\augSpaceBasis^T(\A(\augSpaceBasis\p))$. Finally,
	stage 3 applies augmented PCG in the full space to
	compute a solution satisfying the specified tolerance $\tolj$.
	
	Goal-oriented POD as proposed in Section \ref{sec:podingredients} fits
	naturally into this framework. First, the basis vectors are
	optimally ordered (Remark \ref{rem:podvecorder}): the first few POD basis
	vectors span an optimal subspace in the sense of minimizing the objective
	function in Eq.~\eqref{eq:podOpt}, which is an upper bound for the
	$\Aj$-norm and $\outputMat^T\outputMat$-norm of the error, respectively
	(Theorems \ref{Anorm} and \ref{outNorm}). This implies that the
	first few POD vectors could be employed as $\stageonebasisj$. Second, it
	automatically yields a well-conditioned reduced matrix if $\metric = \Aj$;
	this will be further discussed in Section \ref{sec:integratePOD}. However,
	other truncation methods that satisfy these properties can also be
	considered within the proposed three-stage algorithm.

\subsection{Stage 1}

The objective of Stage 1 is to `jump start' the algorithm by computing an accurate solution at very low cost.
To achieve this, this stage computes $\xGuessj+\stageonebasisj\stageonesolredj
= \projOrth{\Aj}{\xExactj}{\xGuessj+{\rangestageonebasisj}}$ by solving
 \begin{equation} 
\stageonebasisj^T\Aj\stageonebasisj\stageonesolredj = \stageonebasisj^T(\bj -
\Aj\xGuessj)
  \end{equation} 
 directly. The assumptions placed on $\stageonebasisj$ imply that 
 the solution
$\xGuessj+\stageonebasisj\stageonesolredj$ will be accurate and
 a direct
 solve will be inexpensive. Computing this
solution can be executed via Algorithm
\ref{alg:stageone} 
as
$$
( \stageonesolredj, \cholRedj) = \stageonealgorithmname ( \Aj,  \bj -
\Aj\xGuessj, \stageonebasisj ).
$$
\begin{algorithm}[H]
\caption{{\stageonealgorithmname}}
\begin{algorithmic}[1]\label{alg:stageone}
\small
\REQUIRE $\A$, $\b$, $\stageonebasis$
\ENSURE $\stageonesolred$, $\cholRed$
\STATE Compute $\hat\A = \stageonebasis^T\A\stageonebasis$ and $\hat \b =
\stageonebasis^T\b$
\STATE Solve $\hat\A\stageonesolred =
\hat \b$ by Cholesky factorization $\hat\A = \cholRed^T\cholRed$
\end{algorithmic}
\end{algorithm}

\subsection{Stage 2}\label{sec:stage2}

The objective of Stage 2 is to improve upon the stage-1 solution $\xGuessj +
\stageonebasisj\stageonesolredj$ by efficiently solving over the entire
augmenting subspace ${\rangeaugSpaceBasisj}$, whose dimension $\augSpaceDimj$
may be large, while avoiding the cost of explicitly forming the reduced matrix $\augSpaceBasisj^T\Aj\augSpaceBasisj$.
To achieve this, this
stage solves  
\begin{equation}
\augSpaceBasisj^T\Aj \augSpaceBasisj \augSpaceSolredExactj=
\augSpaceBasisj^T\left(\bj - \Aj\xGuessj-
\Aj\stageonebasisj\stageonesolredj\right)
\end{equation}
iteratively to tolerance $\tolstagetwoj$ via augmented CG with augmenting
basis $\stageoneaugj\in\RRstar{\augSpaceDimj\times \stageonedimj}$, which
represents the stage-1 basis in augmenting-subspace coordinates, i.e.,
$\stageonebasisj = \augSpaceBasisj\stageoneaugj$.
This approach is
promising for two reasons. First, it precludes the need to explicitly compute
the reduced matrix, as
	matrix--vector products of the form $\augSpaceBasisj^T\Aj\augSpaceBasisj\p$
	can be computed via $\augSpaceBasisj^T(\Aj(\augSpaceBasisj\p))$ using a single
	matrix--vector product and $\mathcal
	O((2\augSpaceDimj + \bandwidthj)\ndof)$ flops. Second, a small number of iterations will be needed for
	convergence if the reduced matrix $\augSpaceBasisj^T\Aj
	\augSpaceBasisj $ is well conditioned, as has been assumed.

 Therefore, this stage amounts to executing Algorithm
\ref{alg:augmentedCGgen} 
as
$$
(
\stagetwoitj, \stagetwosolredj,
\stagetwoaugj, \redTj
) = \stagetwoalgorithmname(
\augSpaceBasisj^T\Aj\augSpaceBasisj, 
\augSpaceBasisj^T(\bj - \Aj\xGuessj), 
\stageonesolredj,
\stageoneaugj,  \I,   \tolstagetwoj)
$$
with an implementation optimization:
 the solves in steps \ref{step:aug2}
and \ref{step:aug3} of
Algorithm
\ref{alg:augmentedCGgen}
can be performed
directly in $\mathcal O(\stageonedimj^2)$ flops, as the Cholesky factorization of
$\stageonebasisj^T \Aj \stageonebasisj =
\cholRedj^T\cholRedj$ can be
reused from stage 1.

	The (inexact) solution increment $
	\augSpaceBasisj \augSpaceSolredj = 
	\stagetwobasisj\stagetwosolredj$ computed
	by stage 2 lies in the
	range of a (reduced) Krylov basis $\stagetwobasisj =
	\augSpaceBasisj\stagetwoaugj \in \RR{\ndof\times\stagetwoitj}$ with
	${\rangestagetwobasisj}\subseteq {\rangeaugSpaceBasisj}\subseteq
	\RR{\ndof}$ and $\stagetwoaugj\in\RR{\augSpaceDimj\times \stagetwoitj}$,
	where $\stagetwoitj\leq\augSpaceDimj$ denotes the number of stage-2
	iterations.  
This basis satisfies
 \begin{gather} 
\stagetwobasisj^T\Aj\stagetwobasisj = \redTj,
\quad\stagetwobasisj^T\Aj\stageonebasisj = 0,\\
{\rangestageoneaugj} \oplus
{\rangestagetwoaugj} = \rangestageoneaugj
+\krylovj{\augSpaceBasisj^T\Aj\augSpaceBasisj}{\augSpaceBasisj^T(\bj -
\Aj\xGuessj-\Aj\stageonebasisj\stageonesolredj)}{\stagetwoitj},
  \end{gather} 
	while the resulting solution satisfies
	\begin{align}
	\xGuessj+\stageonebasisj\stageonesolredj  + \stagetwobasisj\stagetwosolredj&=
	\projOrth{\Aj}{\xExactj}{\xGuessj+{\rangestageonebasisj}\oplus{\rangestagetwobasisj}}
	\end{align}
Ref.\ \cite{risler2000iaa} proposed a similar idea referred to as the
iterative reuse of Krylov subspaces (IRKS). This approach did not employ
stage-1 direct solve and supported using
either
$\augSpaceBasisj = \vectomat{\kryVec}{j-1}$ or $\augSpaceBasisj =
\kryVecArg{j-1}$ (i.e., no truncation).

\subsection{Stage 3}\label{sec:stage3}

The objective of stage 3 is to continue iterating in the full space until the
specified tolerance $\tolj$ is satisfied. This stage therefore solves
 \begin{equation} 
\Aj\DeltaxExactj = \bj - \Aj\xGuessj-
	\Aj\stageonebasisj\stageonesolredj- \Aj\stagetwobasisj\stagetwosolredj
  \end{equation} 
via augmented PCG with augmenting basis $[\stageonebasisj,\ \stagetwobasisj]$.
Critically, note that the augmenting subspace used in stage 3 satisfies $\range{[\stageonebasisj,\
\stagetwobasisj]}\subseteq\augSpacej$. This choice is made because
$\Aj[\stageonebasisj,\
\stagetwobasisj]$ and 
$[\stageonebasisj,\ \stagetwobasisj]^T\Aj[\stageonebasisj,\ \stagetwobasisj]$---which are required for orthogonalization in step \ref{step:aug3} of
Algorithm \ref{alg:augmentedCGgen}---are available from stages 1 and 2. Treating 
$\augSpacej$ as the augmenting subspace in stage 3 would necessitate computing
$\Aj\augSpaceBasisj$ and $\augSpaceBasisj^T\Aj\augSpaceBasisj$, which would
eliminate any cost savings realized during stages 1 and 2, as computing and
factorizing the reduced matrix $\augSpaceBasisj^T\Aj\augSpaceBasisj$ incurs $\mathcal
O((\bandwidthj + \augSpaceDimj^2)\ndof+\augSpaceDimj^3)$ flops: the same cost as
executing stage 1 with $\stageonebasisj = \augSpaceBasisj$.
Thus, stage 3 executes Algorithm \ref{alg:augmentedCGgen} as 
$$
(
\kit{j},
\redKryj, \kryVecj,
\Tj
) = \stagetwoalgorithmname ( 
 \Aj,  \bj - \Aj\xGuessj,
\left[\stageonesolredj^T,\ \stagetwosolredj^T\right]^T,
[\stageonebasisj,\ \stagetwobasisj], 
\precj,  \tolj
)
$$
with the following optimization:
the solves in steps \ref{step:aug2}
and \ref{step:aug3} can be performed
directly in $\mathcal O(\stageonedimj^2 + \stagetwodimj)$ flops, as the
Cholesky factorization 
$ [\stageonebasisj,\ \stagetwobasisj]^T \A
 [\stageonebasisj,\ \stagetwobasisj] =
\cholj^T\cholj$ is readily available from stages 1 and 2 as 
$$
\cholj = \left[\begin{array}{cc}
\cholRedj & \zero\\
\zero & \sqrt{\redTj}\\
\end{array}\right].
$$


The solution increment $\Deltaxj = \kryVecj\redKryj$ computed by stage 3
lies in the range of a Krylov basis $\kryVecj\in\RRstar{\ndof\times\kit{j}}$,
where $\kit{j}\leq\ndof$ denotes the number of stage-3 iterations.
This basis satisfies
 \begin{gather} 
\kryVecj^T\Aj\kryVecj = \Tj,
\quad\kryVecj^T\Aj[\stageonebasisj,\ \stagetwobasisj]  = \zero,\\
\rangestageonebasisj \oplus\rangestagetwobasisj
\oplus \rangekryVecj= \rangestageonebasisj
+\augSpaceBasisj\krylovj{\augSpaceBasisj^T\Aj\augSpaceBasisj}{\augSpaceBasisj^T(\bj -
\Aj\xGuessj-\Aj\stageonebasisj\stageonesolredj)}{\stagetwoitj} +
\krylovjkArbUsed{\kit{j}}(\xGuessj+\stageonebasisj\stageonesolredj  + \stagetwobasisj\stagetwosolredj)
  \end{gather} 
with $\Tj\in\RR{\kit{j}\times\kit{j}}$ diagonal.
	The resulting solution satisfies
	\begin{equation}
	\xGuessj+\stageonebasisj\stageonesolredj  + \stagetwobasisj\stagetwosolredj
	+ \kryVecj\redKryj = 
	\projOrth{\Aj}{\xExactj}{\xGuessj+{\rangestageonebasisj}\oplus{\rangestagetwobasisj}\oplus{\rangekryVecj}}.
	\end{equation}

\subsubsection{Accounting for entire augmenting
subspace}\label{sec:stage3special}

To preserve computational-cost savings, the stage-3 implementation described above augments with the subspace
$\range{[\stageonebasisj,\ \stagetwobasisj]}\subseteq{\rangeaugSpaceBasisj}$.
However, we know that the stage-2 solution $\xGuessj+\stageonebasisj\stageonesolredj  + \stagetwobasisj\stagetwosolredj$ satisfies
 \begin{equation} 
 \|\augSpaceBasisj^T\rj(\xGuessj+\stageonebasisj\stageonesolredj  + \stagetwobasisj\stagetwosolredj)\|\leq\tolstagetwoj ,
  \end{equation} 
where  $\tolstagetwoj$ is the stage-2 tolerance. If this tolerance is small,
then orthogonality condition \eqref{eq:equivalenceResidualOrth} implies that this solution is 
\textit{nearly} optimal over the entire augmenting subspace, i.e.,
 \begin{equation*} 
	\xGuessj+\stageonebasisj\stageonesolredj  + \stagetwobasisj\stagetwosolredj = 
	\projOrth{\Aj}{\xExactj}{\xGuessj+{\rangestageonebasisj}\oplus{\rangestagetwobasisj}}
\approx
\projOrth{\Aj}{\xExactj}{\xGuessj+{\rangeaugSpaceBasisj}}.
  \end{equation*}
Therefore, employing
$\augSpaceBasis\leftarrow\augSpaceBasisj$ in stage 3
may reduce the number of iterations to convergence, as this
ensures that new search directions remain $\Aj$-orthogonal to
the full augmenting subspace ${\rangeaugSpaceBasisj}$ over which the solution
has already been computed to tolerance $\tolstagetwoj$. 
However, for this approach to be practical, the reduced solves in 
steps \ref{step:aug2} and \ref{step:aug3} must be performed efficiently, i.e.,
without
assembling and factorizing the matrix
$\augSpaceBasisj^T\Aj\augSpaceBasisj$. We therefore
propose performing these solves \textit{iteratively}, i.e, by executing stage 2
within stage 3. In particular, this stage-3 variant executes Algorithm
\ref{alg:augmentedCGgen} 
as
$$
( \kit{j}, \redKryj, \kryVecj, \Tj) = \stagetwoalgorithmname ( \Aj,  \bj -
\Aj\xGuessj, \stageoneaugj\stageonesolredj+\stagetwoaugj\stagetwosolredj,
\augSpaceBasisj, \precj, \tolj)
$$
with the following optimizations:
\begin{itemize}
\item
At iteration $k$, the solves in steps \ref{step:aug2} and \ref{step:aug3} can
be performed iteratively by executing Algorithm \ref{alg:augmentedCGgen} as
$$
(
\innerstagetwoitjk\, \innerstagetwosolredjk\,
\innerstagetwoaugjk\, \innerredTjk\
)
=
\stagetwoalgorithmname
(
\augSpaceBasisj^T\Aj\augSpaceBasisj, 
\augSpaceBasisj^T(
\Aj  \zArg{k+1}_j
), 
\zero,
[\stageoneaugj,\ \stagetwoaugj,\
\innerstagetwoaugjArg{0},\ \cdots,\ \innerstagetwoaugjArg{k-1}],  \I, 
\innertolstagetwoj
)
$$
with the following optimizations:
 \begin{itemize} 
 \item The solves in steps \ref{step:aug2} and \ref{step:aug3} can be
 performed directly in $\mathcal O(\stageonedimj^2 + \stagetwodimj +
 \sum_{\ell=0}^{k-1}\innerstagetwoitjArg{\ell})$ flops, as the
 Cholesky factorization $\augSpaceBasis^T \A
 \augSpaceBasis = \innercholRedjk\innercholRedjk$ with
$$
\innercholRedjk = \left[\begin{array}{ccccc}
\cholRedj & \zero & \zero& \cdots & \zero\\
\zero & \sqrt{\redTj} & \zero& \cdots & \zero\\
\zero & \zero & \sqrt{\innerredTjArg{0}} & \cdots& \zero\\
\vdots & \vdots & \vdots & \ddots& \vdots\\
\zero & \zero & \zero & \zero& \sqrt{\innerredTjArg{k-1}}
\end{array}\right]
$$
 can be reused from stage 1, stage 2, and earlier stage-3 iterations.
	 \end{itemize}
\end{itemize}
This solution increment $\Deltaxj = \kryVecj\redKryj$ also
lies in the range of a Krylov basis $\kryVecj\in\RRstar{\ndof\times\kit{j}}$;
however, the basis now satisfies
 \begin{gather} 
\kryVecj^T\Aj\kryVecj = \Tj,
\quad\kryVecj^T\Aj[\stageonebasisj,\ \stagetwobasisj,\ \augSpaceBasisj
\innerstagetwoaugjArg{0},\ \cdots,\
\augSpaceBasisj\innerstagetwoaugjArg{\kit{j}-1}
]  = \zero,\\
\augSpacej
\oplus \rangekryVecj\approx \augSpacej
 +
\krylovjkArbUsed{\kit{j}}(\xGuessj+\stageonebasisj\stageonesolredj  +
\stagetwobasisj\stagetwosolredj),
  \end{gather} 
while
	the resulting solution satisfies
	\begin{equation}
	\xGuessj+\stageonebasisj\stageonesolredj  + \stagetwobasisj\stagetwosolredj
	+ \kryVecj\redKryj \approx
\projOrth{\Aj}{\xExactj}{\xGuessj+{\rangeaugSpaceBasisj}}.
	\end{equation}
Note that the solution does not associate with an exact orthogonal projection.


	\subsection{Overall algorithm}\label{sec:overall}
Algorithm \ref{alg:overall} reports the proposed three-stage algorithm. Here,
$\augSpaceThreshold$ denotes the maximum number of accumulated vectors to
preserve before truncation. 
The variable $\newvectorsdirect\in[0,1]$ is a
threshold for determining which new vectors should be included in the stage 1
basis.
Finally, $\stagethreeoptionone$ is a boolean variable whose value is
\true if stage 3 enforcing orthogonality with the entire augmenting subspace (not just stage 1 and 2 directions)
as discussed in Section
\ref{sec:stage3special}. 

\begin{algorithm}[htbp]
\caption{\overallalgorithmname}
\begin{algorithmic}[1]\label{alg:overall}
\small
\REQUIRE $\{\Aj\}_{j=1}^{\nlhs}$,
$\{\bj\}_{j=1}^{\nlhs}$,
$\{\xGuessj\}_{j=1}^{\nlhs}$,
forcing sequence $\{\tolj\}_{j=1}^{\nlhs}$, storage threshold
$\augSpaceThreshold$,
stage-1 option $\newvectorsdirect$, stage-3 option $\stagethreeoptionone$
\ENSURE 
\STATE $\stageonebasisArg{1} \leftarrow \emptyset$,
$\stagetwobasisArg{1}\leftarrow \emptyset$, $\jcompress
\leftarrow 0$
\FOR{$j = 1, \ldots, \nlhs$}
\IF{$\stageonebasisj\neq\emptyset$}
\STATE \textit{Stage 1}: ($\stageonesolredj$, $\cholRedj$) =
\stageonealgorithmname($\Aj$, $\bj - \Aj\xGuessj$,
$\stageonebasisj$)
\ENDIF
\IF{$\stagetwobasisj\neq\emptyset$}
\STATE \textit{Stage 2}: ($\stagetwoitj$, $\stagetwosolredj$,
$\stagetwoaugj$, $\redTj$) = \stagetwoalgorithmname($
\augSpaceBasisj^T\Aj\augSpaceBasisj$, $
\augSpaceBasisj^T(\bj - \Aj\xGuessj)$, 
$\stageonesolredj$,
$\stageoneaugj$, $\I$,  $\tolstagetwoj$).\\
Note optimizations discussed in Section
\ref{sec:stage2}
\ENDIF
\IF[orthogonalize against entire $\augSpaceBasisj$]{$\stagethreeoptionone=$\true}
\STATE \textit{Stage 3}: 
($\kit{j}$, $\redKryj$, $\kryVecj$,
$\Tj$) = \stagetwoalgorithmname(
$\Aj$, $ \bj - \Aj\xGuessj$,
$\stageoneaugj\stageonesolredj+\stagetwoaugj\stagetwosolredj$,
$\augSpaceBasisj$, $\precj$, and $\tolj$).\\ Note optimizations discussed in Section \ref{sec:stage3special}
\ELSE
\STATE \textit{Stage 3}: 
($\kit{j}$,
$\redKryj$, $\kryVecj$,
$\Tj$) = \stagetwoalgorithmname(
$ \Aj$, $ \bj - \Aj\xGuessj$,
$\left[\stageonesolredj^T\ \stagetwosolredj^T\right]^T$,
$[\stageonebasisj\ \stagetwobasisj]$, $
\precj$, $ \tolj$).\\ Note optimizations discussed in Section \ref{sec:stage3}
\ENDIF
\STATE$\stageonebasisArg{j+1} \leftarrow\stageonebasisArg{j}$
\STATE $\stageonesetj \leftarrow \{i\ |\
[\Tj]_{ii}/\sum_{k=1}^{\kit{j}}[\Tj]_{kk} >
\newvectorsdirect \}$
\FOR{$k\in\stageonesetj$}
\STATE $\stageonebasisArg{j+1} \leftarrow[\stageonebasisArg{j+1},\
[\kryVecj]_k\sqrt{[\Tj]_{kk}}]$
\ENDFOR
\IF[include new vectors in stage-1 basis]{$\newvectorsdirect=$\true}
\STATE $\stageonebasisArg{j+1} \leftarrow[\stageonebasisArg{j},\ \kryVecj\sqrt{\Tj}]$
\ELSE
\STATE$\stageonebasisArg{j+1} \leftarrow\stageonebasisArg{j}$
\ENDIF
\IF[truncate]{$\augSpaceDimjp > \augSpaceThreshold$}
\STATE\label{step:compress1} $\augSpaceBasisjp =
\compressionalgorithmname([\augSpaceBasisj,\kryVecj])$ (note that $\currentDataj =[\augSpaceBasisj,\kryVecj]$)
\STATE\label{step:enforceOrthog} Enforce $\Aj$-orthogonality:
$\augSpaceBasisjp^T\Aj\augSpaceBasisjp = \cholEnforcejT\cholEnforcej$
(Cholesky factorization),
$\augSpaceBasisjp\leftarrow \augSpaceBasisjp\cholEnforcejinv$
\STATE\label{step:compresslast} Determine $\stageonebasisArg{j+1}$ with
$\range{\stageonebasisArg{j+1}}\subseteq \range{\augSpaceBasisArg{j+1}}$
\STATE Update lastest-truncation index $\jcompress \leftarrow j$
\ELSE
\STATE $\augSpaceBasisjp \leftarrow \left[\augSpaceBasisj,\ \kryVecj\sqrt{\Tj}\right]$
\ENDIF
\ENDFOR
\end{algorithmic}
\end{algorithm}

Step \ref{step:enforceOrthog} of Algorithm \ref{alg:overall} ensures
that
$\augSpaceBasisArg{\jcompress+1}^T\A_{\jcompress}\augSpaceBasisArg{\jcompress+1}
= \I$. We now show that this step is critical for enabling fast stage-2
convergence: it leads to well-conditioned reduced matrices for slowly varying
matrix sequences.

\begin{theorem}\label{wellConditioned}
In Algorithm \ref{alg:overall}, if either 1) $\newvectorsdirect=\true$ and $\stageonebasisArg{\jcompress +
1} = \augSpaceBasisArg{\jcompress + 1}$ in step \ref{step:compresslast}, or
2) $\stagethreeoptionone=\true$ and $\tolstagetwoArg{k} = \innertolstagetwoArg{k} =
0$, $k=\jcompress + 1,\ldots,j-1$,
then 
\begin{equation} \label{eq:wellConditioned}
\|\augSpaceBasisj\Aj\augSpaceBasisj - \I\| \leq \sum_{k = \jcompress +
1}^j\|\augSpaceBasisArg{k}\|^2\|\A_k - \A_{k-1}\|.
\end{equation} 
 
\end{theorem}


	\subsubsection{Integration with goal-oriented POD}\label{sec:integratePOD}

Section \ref{sec:podingredientssummary} described four sets of goal-oriented
POD ingredients (with $\memory = j - \jcompress$) that could be employed as \compressionalgorithmname\ in step
\ref{step:compress1} of Algorithm \ref{alg:overall}.
In all cases, step
\ref{step:compresslast} of Algorithm \ref{alg:overall} can be executed by
setting $\stageonebasisArg{j+1}$ equal to the first $\stageonedimArg{j}$
vectors of $\augSpaceBasisArg{j+1}$, where  $\stageonedimArg{j}$ is determined
from the appropriate POD algorithm (Algorithm \ref{alg:PODEVD} or \ref{alg:PODSVD}) with a
modest statistical energy criterion $\energyCrit\leftarrow
\energyCritStageOne$ with $\energyCritStageOne\leq \energyCritStageTwo$.

The first two options 
$\augSpaceBasisjp  =
\podBasis{\augSpaceDimjp}{\currentDataj}{\augSpaceWeightsVecPrev}{\Aj}$ and 
	$
	\augSpaceBasisjp  = \podBasis{\augSpaceDimjp}{\currentDataj}{\augSpaceWeightsVecRBF(\memory)}{\Aj}
	$ expose two implementation
optimizations.
First, during any truncation iteration $\jcompress$, the algorithm should
employ $\stageonebasisArg{\jcompress} = \augSpaceBasisArg{\jcompress}$. The
reason is that the truncation step requires computing
$\augSpaceBasisArg{\jcompress}^T\A_{\jcompress}\augSpaceBasisArg{\jcompress}$ explicitly;
because this is the same matrix used in stage one with
$\stageonebasisArg{\jcompress} =
\augSpaceBasisArg{\jcompress}$, employing this choice can lead to faster stage-3 convergence
at no additional computational cost.
Second, orthogonalization step \ref{step:enforceOrthog} in Algorithm \ref{alg:overall}
requires \textit{no operations}, as the basis is automatically
$\A_{\jcompress}$-orthogonal (see Eq.~\eqref{eq:podOrth}). 


\section{Numerical experiments}\label{sec:experiments}

\subsection{Problem description}\label{sec:problemDescription}
We now assess the proposed methodology using model problems from
Sierra/SolidMechanics \cite{mitchell2001adagio}, which is a Lagrangian,
three-dimensional code for finite element analysis of solids and structures. 

\subsubsection{Problem 1: Pancake problem}\label{sec:pancakeDescription}
\begin{figure}[htb] 
\centering 
         \begin{subfigure}[b]{0.4\textwidth}
                \includegraphics[width=\textwidth]{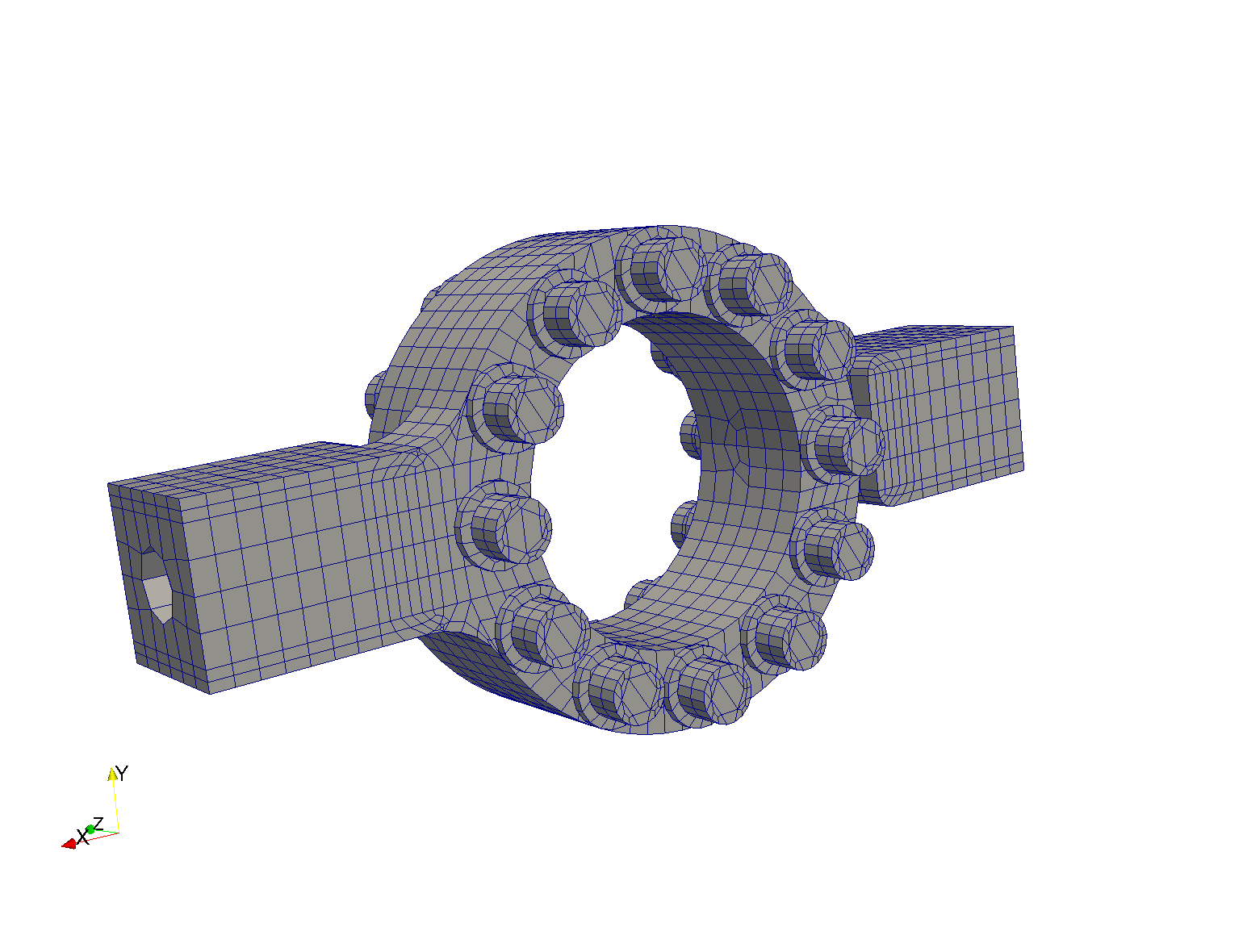} 
                \caption{Finite-element mesh}
        \end{subfigure}
         \begin{subfigure}[b]{0.4\textwidth}
                \includegraphics[width=\textwidth]{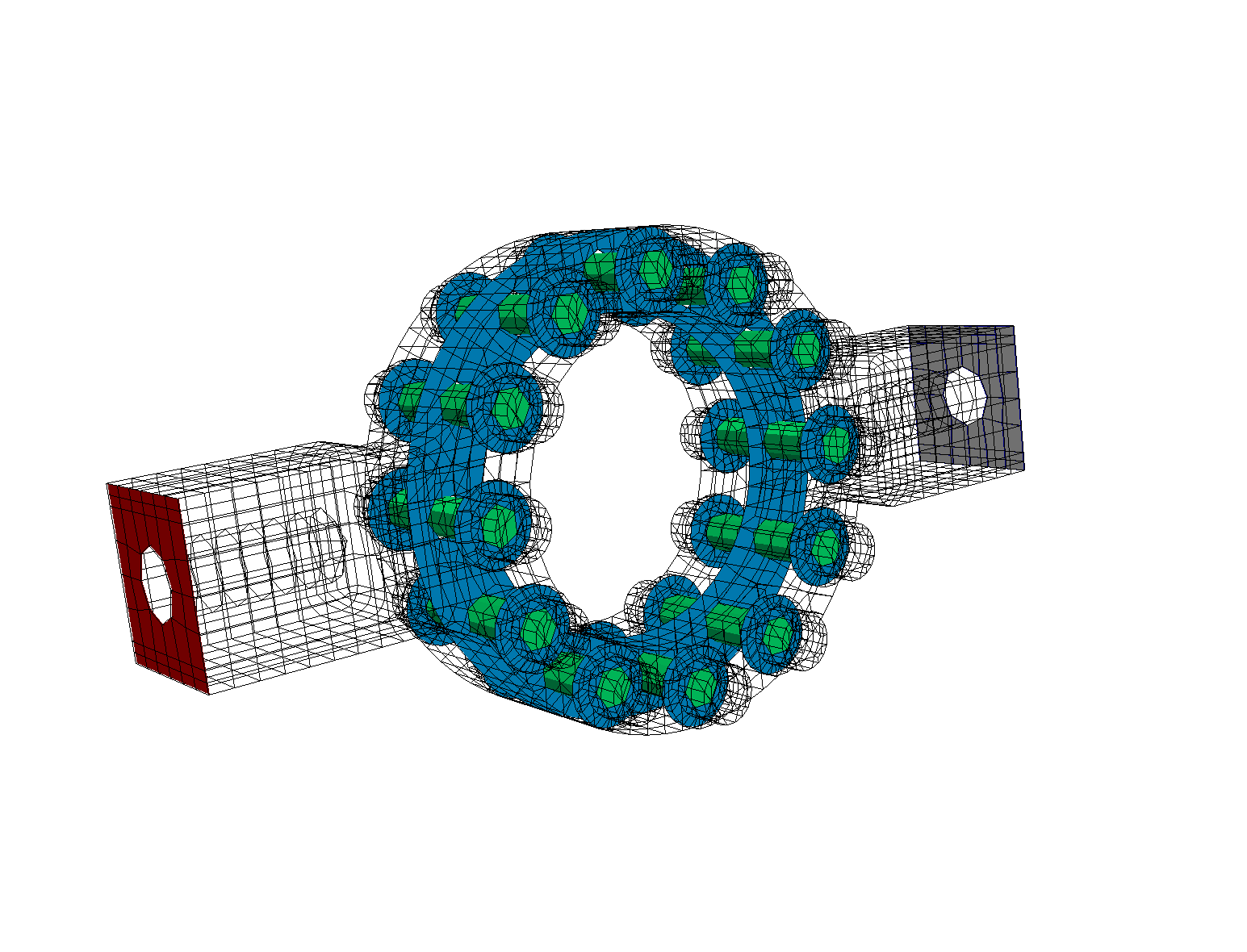} 
                \caption{Pressure-loaded surface (red), contact surfaces (blue), prescribed
								temperature (green), Dirichlet boundary condition (gray).}
  
        \end{subfigure}
				\caption{Pancake problem.}
\label{fig:pancakeMesh} 
\end{figure}

We first consider computing the quasistatic response of the
`pancake' domain pictured in Figure \ref{fig:pancakeMesh}. The material is
steel, which is characterized by a Young's modulus of $E = 2.0\times 10^{8}\
\frac{\text{N}}{\text{mm}\ \cdot\ \text{s}^2}$, Poisson's ratio of
$\nu = 0.3$, and density of 
$\rho = 7.86\times 10^{-6}\ \text{kg}/\text{mm}^3$. The 
logarithmic thermal
strain of steel is linearly dependent on temperature $\strainThermal = (11.7\times
10^{-6})\temperature\ \frac{1}{\text{K}}$, where $\temperature$ denotes the
change in 
temperature (in Kelvin) from the reference temperature.

\begin{figure}[htb] 
\centering 
         \begin{subfigure}[b]{0.40\textwidth}
                \includegraphics[width=\textwidth]{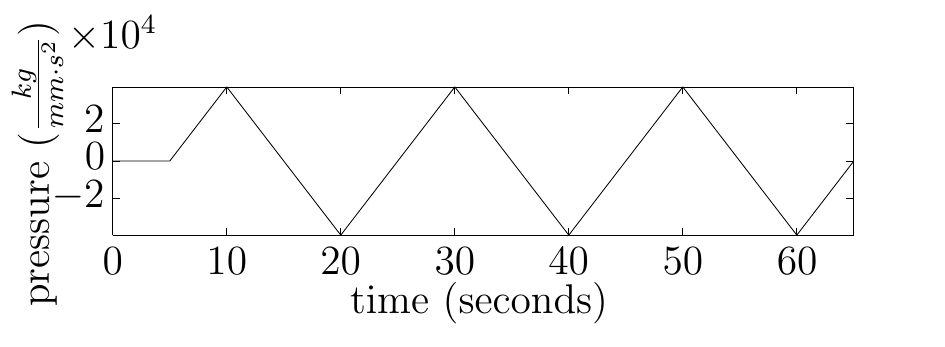} 
                \caption{Time-dependent pressure load applied to leftmost surface (extrema are
$\pm 3.94\times 10^4\ \frac{\text{kg}}{\text{mm}\cdot\text{s}^2}$). }
	\label{fig:pressureLoad}
        \end{subfigure}
         \begin{subfigure}[b]{0.40\textwidth}
                \includegraphics[width=\textwidth]{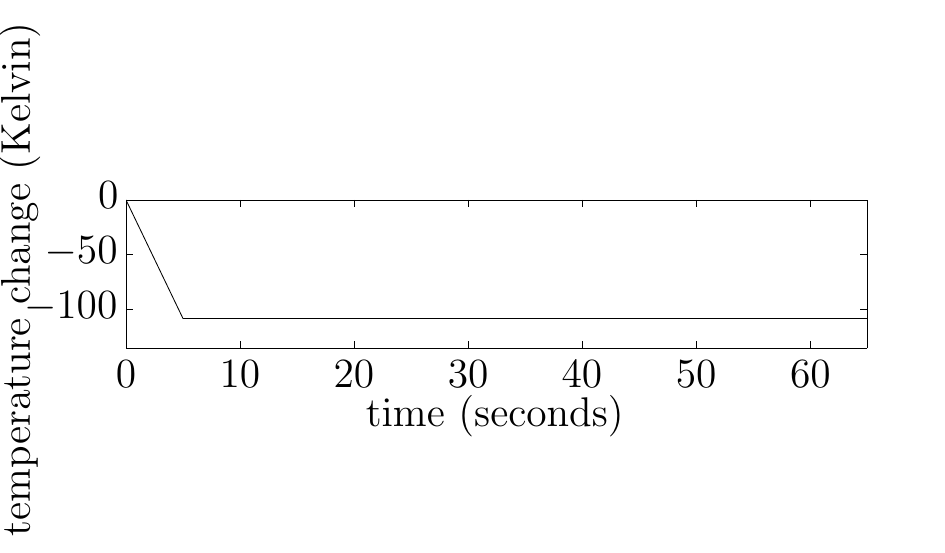} 
                \caption{Time-dependent temperature load applied to bolts.}
\label{fig:tempLoad}
        \end{subfigure}
				\caption{Time-dependent loadings for the pancake problem.}
\end{figure}

The $x$-, $y$-, and $z$-displacements of the rightmost surface (gray in Figure
\ref{fig:pancakeMesh}(b)) are set to zero. The $x$- and $y$- displacements of
the leftmost surface (red in Figure \ref{fig:pancakeMesh}(b)) are set to
zero; this surface is also subjected to the time-dependent pressure load depicted in
Figure \ref{fig:pressureLoad}.

The time-dependent thermal load depicted in Figure \ref{fig:tempLoad} is applied to the bolts (green components in
Figure \ref{fig:pancakeMesh}(b)) to emulate a pre-loading condition.
Contact surfaces are shown in blue in Figure \ref{fig:pancakeMesh}(b); they
are enforced by an augmented Lagrangian approach with DASH search using a
penalty factor of 1.25 and a friction coefficient of 0.5.  

The problem is discretized by the finite-element method 
using a mesh generated by the SIERRA toolkit \cite{edwards2010sierra}.
The mesh
consists of 9108 nodes and 4719 hexahedral elements. At each node, there
are three degrees of freedom (the $x$-, $y$-, and $z$-displacements), which
leads to $27,324$ total degrees of freedom for the finite-element model.

As the time scales of the load application are relatively large, we neglect
inertial effects and consider solving the quasi-static equations
\begin{equation}\label{eq:contactEOM}
\fint{\displaci} + \fcontact{\displaci}{\penaltyl}= \fext(\timei),
\quad l\innatseq{\npenalty},\quad i\innatseq{\ntimeSteps},
\end{equation}
where $\ntimeSteps$
denotes the number of time steps, $\fintNo:\RR{\ndof}\rightarrow\RR{\ndof}$ is a nonlinear operator
representing the internal
forces,
$\fcontactNo:\RR{\ndof}\times \RR{+}\rightarrow\RR{\ndof}$ is nonlinear in its
first argument and represents contact
forces, $\fext:
\left[0,65\right]\rightarrow\RR{\ndof}$ is the external-force vector, and
$\displaci\in\RR{\ndof}$ denotes the displacement at time $t_i$. As the
contact constraints are enforced using an augmented Lagrangian  
approach within
a continuation loop, 
$\penaltyl\in\RR{+}$ with $\penaltyArg{l}\leq\penaltyArg{l+1}$,
$l\innat{\npenalty}$ denotes the penalty factor at continuation iteration
$l$. Note
that these equations also include equality constraints arising from the
Dirichlet boundary conditions. 

Solving Eqs.~\eqref{eq:contactEOM} is mathematically equivalent to  solving
\begin{equation} \label{eq:optProb}
\underset{\dummyvec\in\RR{\ndof}}{\text{minimize}}\quad \objective{i}{l}{\dummyvec}
\end{equation} 
with $-\nabla\objectiveMap{i}{l}:{\dummyvec}\mapsto\fext(\timei) -
\fint{\dummyvec} -
\fcontact{\dummyvec}{\penaltyl} $ for a local minimum.
Sierra/SolidMechanics solves problem \eqref{eq:optProb} using the nonlinear conjugate gradient method with
a displacement-dependent
preconditioner  
$\precNonlinlNo: \displacDummy
\mapsto \nabla_{\displac}\fint{\displacDummy} +
\nabla_{\displac}\fcontact{\displacDummy}{\penaltyl}\in\spd{\ndof}$. This results in a sequence of
linear systems of the original form in Eqs.~\eqref{eq:sequence}: one at each nonlinear
conjugate-gradient iteration.
Here, $\Aj = \precNonlinl{\displacilk}$, $\bj =
\nabla\objective{i}{l}{\displacilk}$, $\displacilk\in\RR{\ndof}$ denotes the
displacement at time step $i$, continuation iteration $l$, and nonlinear
conjugate gradient iteration $k$. The mapping between
timestep $i$, continuation iteration $l$, nonlinear conjugate-gradient
iteration $k$ and the
index of the linear system $j$ is provided by 
$j:(i,l,k) \mapsto k + \sum_{\mathsf i=1}^i\sum_{\mathsf l=1}^lK(\mathsf
i,\mathsf l)$, where $K( i,  l)$ denotes the
number of conjugate-gradient iterations needed to solve problem
\eqref{eq:optProb}.

For this problem, the total number of linear systems we consider is $\nlhs =
47$. Each of these linear systems is preconditioned using a three-level
algebraic multigrid (AMG) preconditioner $\precj$ with incomplete Cholesky
smoothing for both pre-smoothing and post-smoothing. This preconditioner tends
to be expensive to apply, especially when compared with the cost of a
matrix--vector product
for this system.

\subsubsection{Problem 2: I-beam problem}

\begin{figure}[htb] 
\centering 
         \begin{subfigure}[b]{0.3\textwidth}
                \includegraphics[width=\textwidth]{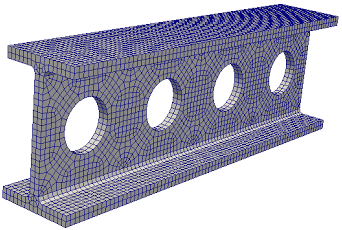} 
                \caption{Finite-element mesh}
	\label{fig:ibeamMeshDetail}
        \end{subfigure}
         \begin{subfigure}[b]{0.3\textwidth}
                \includegraphics[width=\textwidth]{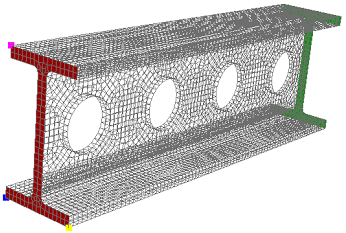} 
                \caption{Mesh with boundary-condition information.}
\label{fig:ibeamBC}
        \end{subfigure}
				\caption{I-beam problem.}
\label{fig:ibeamMesh} 
\end{figure} 

We next consider computing the quasistatic response (neglecting thermal
effects) of an
I-beam with holes in the web section; the domain is pictured in Figure \ref{fig:ibeamMeshDetail}. The material is
steel 304L, which is characterized by a Young's modulus of 
$E = 2.1\times 10^{8}\ \frac{\text{N}}{\text{mm}\ \cdot\ \text{s}^2}$,
Poisson's ratio of $\nu = 0.33$, and density of $\rho = 7.8\times 10^{-6}\
\text{kg}/\text{mm}^3$. 
The $x$-, $y$-, and $z$-displacement of the bottom-left point (blue point in Figure
\ref{fig:ibeamBC}) are set to zero. The $x$- and $y$- displacement of
the bottom-right point (yellow point in Figure \ref{fig:ibeamBC}) are set to
zero. The $x$-displacement of the top-left point (magenta point in Figure
\ref{fig:ibeamBC}) is set to zero. Finally, a torsional traction is
applied the end surfaces (red and green in Figure \ref{fig:ibeamMesh}); Ref.\
\cite[Eq.~(11)]{bishop2015direct} reports details on the tractional loading, where the scale
factor for the current problem is 0.01.  
The mesh 
consists of 13,137 nodes and 8,576 hexahedral elements. Each node is
characterized by 
three degrees of freedom (the $x$-, $y$-, and $z$-displacements), which
leads to $39,411$ total degrees of freedom in the finite-element model.

We again neglect inertial effects and consider solving the quasi-static
equations
\begin{equation}\label{eq:iBeamEOM}
\fint{\displaci} = \fext(\timei),
\quad i\innatseq{\ntimeSteps}.
\end{equation}
These equations also include equality constraints arising from the
Dirichlet boundary conditions; note that the continuation loop for computing contact forces does not appear
in this problem.
Solving Eqs.~\eqref{eq:iBeamEOM} is equivalent to solving
\begin{equation} \label{eq:optProbiBeam}
\underset{\dummyvec\in\RR{\ndof}}{\text{minimize}}\quad \objectiveibeam{i}{\dummyvec}
\end{equation} 
with $-\nabla\objectiveibeamMap{i}:\dummyvec\mapsto \fext(\timei) -
\fint{\dummyvec} $.
We again use the nonlinear conjugate gradient method to solve problem
\eqref{eq:optProbiBeam}; we also again employ
a displacement-dependent 
preconditioner $\precNonlinNo: \displac  \mapsto \nabla_{\displac}\fint{\displac} +
\nabla_{\displac}\fcontact{\displac}{\penaltyl}\in\spd{\ndof}$. This results in a sequence of
linear systems of the original form in Eqs.~\eqref{eq:sequence}: one at each nonlinear
conjugate-gradient iteration.
Here, $\Aj = \precNonlin{\displacik}$ and $\bj =
\nabla\objectiveibeam{i}{\displacik}$, where the mapping between time step and
nonlinear conjugate-gradient iteration is provided by 
$j:(i,k) \mapsto k + \sum_{\mathsf i=1}^iK(\mathsf i)$. Here, $K( i)$ denotes the
number of conjugate-gradient iterations needed to solve problem
\eqref{eq:optProbiBeam}.
For this problem, the total number of linear systems we consider is $\nlhs =
49$.
Each of these linear
systems is preconditioned using the same multigrid preconditioner $\precj$ as
described in Section \ref{sec:pancakeDescription}; the only modification is
that four levels of multigrid are used in this case due to the larger number
of degrees of freedom.

\subsection{Experimental setup}

We implemented our proposed method in Matlab and ran experiments on a Macbook
Pro with Intel 2.7 GHz i5 processor and 8 GB of RAM; the implementation
performs the linear-system solves after reading in the linear systems
generated by Sierra/SolidMechanics as described in Section
\ref{sec:problemDescription}.  

For all problems, we test our
framework using a full-orthogonalization method (FOM) rather than the
conjugate-gradient recurrence; this amounts to replacing Step
\ref{step:aug3part2} in Algorithm \ref{alg:augmentedCGgen} with the following
\cite{saad2003iterative}: 
\begin{algorithmic}
\renewcommand{\baselinestretch}{1}
\small\normalsize
\STATE $ \pArg{k+1} = \zArg{k+1}  - \augSpaceBasis\augComponentArg{k+1}$
\FOR {$i = 1,\ldots,k$}
\item 	\hspace{\algorithmicindent}$\beta^{(k+1),i} = \frac{(\rArg{k+1})^T
\zArg{k+1}}{(\rArg{i})^T \zArg{i}}$\\
	\hspace{\algorithmicindent}$ \pArg{k+1} = \pArg{k+1} + \beta^{(k+1),i} \pArg{i}$
\ENDFOR
\end{algorithmic}
In exact arithmetic, this modification does not affect the solution. However, in
finite precision, this modification ensures that that the basis $\kryVecj$ is
full rank; this is sometimes necessary to ensure nonsingular systems during
stage-1 solves.  Removing the effect of possible rank deficiency from the
numerical experiments also simplifies interpretation of the results. 

Iterative-solver performance can be measured in three primary ways: the number
of incurred matrix--vector products, the number of stage-3 iterations (which
is equal to the number of preconditioner applications), and the wall time
per linear-system solve. While the wall time is the most important metric in
practice, we report all three metrics (in terms of their averages over all
linear systems) to provide a more complete picture of performance, as the
specific linear system and choice of preconditioner can have a strong effect
on the relative cost of operations.

\subsection{Method comparison: Pancake problem} \label{sec:pancakeResults}

This section compares the following methods:
\begin{enumerate} 
\item \textit{FOM}. This approach solves each linear system
independently without recycling using the full orthogonalization method.
\item \textit{No truncation}. This approach does not perform truncation, and
employs $\augSpaceThreshold = \infty$, $\newvectorsdirect=1$, and
$\stagethreeoptionone=0$ in Algorithm \ref{alg:overall}. Because it employs
all Krylov vectors as the augmenting-subspace basis, it requires the fewest
stage-3 iterations and, thus, the fewest number of preconditioner
applications; however, the memory and orthogonalization costs are the
largest for this method.
\item \textit{DF(100,0)}.
This is the standard approach for deflation-based truncation, which places all
augmenting-subspace basis vectors in the stage-1 basis.  Here, Algorithm
\ref{alg:overall} parameters are $\augSpaceThreshold = 200$,
$\newvectorsdirect=1$,  $\stagethreeoptionone=0$, and $\stageonebasisArg{j+1} =
\left[\left[\augSpaceBasisjp\right]_1\ \cdots\
\left[\augSpaceBasisjp\right]_{100}\right]$ in step
\ref{step:compresslast}. In step \ref{step:compress1}, the augmenting space is
computed by solving Eq.~\eqref{eq:deflationEVP} for $j\leftarrow j+1$ and setting 
$\augSpaceBasisjp \leftarrow
\vectomat{\currentDataj\deflationEigenvector}{100}$.
\item \textit{POD(100,0)}. This approach employs POD truncation and places
augmenting-subspace basis 
vectors in the stage-1 basis. Algorithm \ref{alg:overall} parameters are $\augSpaceThreshold = 200$,
$\newvectorsdirect=1$,  $\stagethreeoptionone=0$, $ \augSpaceBasisjp  =
\podBasis{100}{\currentDataj}{\augSpaceWeightsVecRBF(j - \jcompress )}{\Aj} $
in step \ref{step:compress1}, and $\stageonebasisArg{j+1} =
\left[\left[\augSpaceBasisjp\right]_1\ \cdots\
\left[\augSpaceBasisjp\right]_{100}\right]$ in step
\ref{step:compresslast}.
\item \textit{POD(5,95)it stg1}. 
This approach employs POD truncation, places only the dominant POD modes in
the stage-1 basis, and places all post-truncation Krylov vectors in the
stage-1 basis.
Algorithm \ref{alg:overall} parameters are
$\augSpaceThreshold = 200$,
$\newvectorsdirect=1$,  $\stagethreeoptionone=1$, $ \augSpaceBasisjp  =
\podBasis{100}{\currentDataj}{\augSpaceWeightsVecRBF(j - \jcompress )}{\Aj} $
in step \ref{step:compress1}, and $\stageonebasisArg{j+1} =
\left[\left[\augSpaceBasisjp\right]_1\ \cdots\
\left[\augSpaceBasisjp\right]_{5}\right]$ in step
\ref{step:compresslast}.
Note that $\tolstagetwoj = 10^{-4}\tolj$ for $\tolj \ge 10^{-3}$, and
$\tolstagetwoj =
10^{-5}\tolj$ otherwise, while $\innertolstagetwoj = 10^{-2}\tolj$ for all
tolerances. 
\item \textit{POD(5,95)it mixed}. 
This approach employs POD truncation, places only the dominant POD modes in
the stage-1 basis, and places only the dominant post-truncation Krylov vectors
in the stage-1 basis. Algorithm \ref{alg:overall} parameters are
employs $\augSpaceThreshold = 200$,
$\newvectorsdirect=1\times 10^{-3}$,  $\stagethreeoptionone=1$, $ \augSpaceBasisjp  =
\podBasis{100}{\currentDataj}{\augSpaceWeightsVecRBF(j - \jcompress )}{\Aj} $
in step \ref{step:compress1}, and $\stageonebasisArg{j+1} =
\left[\left[\augSpaceBasisjp\right]_1\ \cdots\
\left[\augSpaceBasisjp\right]_{5}\right]$ in step
\ref{step:compresslast}.
Note $\tolstagetwoj = 10^{-4}\tolj$ for $\tolj \ge
10^{-2}$, and $\tolstagetwoj = 10^{-6}\tolj$ otherwise and $\innertolstagetwoj =
10^{-2}\tolj$ for all $\tolj$.
\item \textit{POD(5,95)it stg2}. 
This approach employs POD truncation, places only the dominant POD modes in
the stage-1 basis, and places none of the post-truncation Krylov vectors
in the stage-1 basis.  Algorithm \ref{alg:overall} parameters are
employs $\augSpaceThreshold = 200$,
$\newvectorsdirect=0$,  $\stagethreeoptionone=1$, $ \augSpaceBasisjp  =
\podBasis{100}{\currentDataj}{\augSpaceWeightsVecRBF(j - \jcompress )}{\Aj} $
in step \ref{step:compress1}, and $\stageonebasisArg{j+1} =
\left[\left[\augSpaceBasisjp\right]_1\ \cdots\
\left[\augSpaceBasisjp\right]_{5}\right]$ in step
\ref{step:compresslast}.
Note
$\tolstagetwoj = 10^{-4}\tolj$ for $\tolj \ge 10^{-2}$, and $\tolstagetwoj =
10^{-7}\tolj$ otherwise. $\innertolstagetwoj = 10^{-2}\tolj$ for $\tolj \ge
10^{-3}$ and $\innertolstagetwoj = 10^{-3}\tolj$ otherwise. 
\end{enumerate}

\begin{figure}[htbp] 
\centering 
         \begin{subfigure}[b]{0.48\textwidth}
                \includegraphics[width=\textwidth]{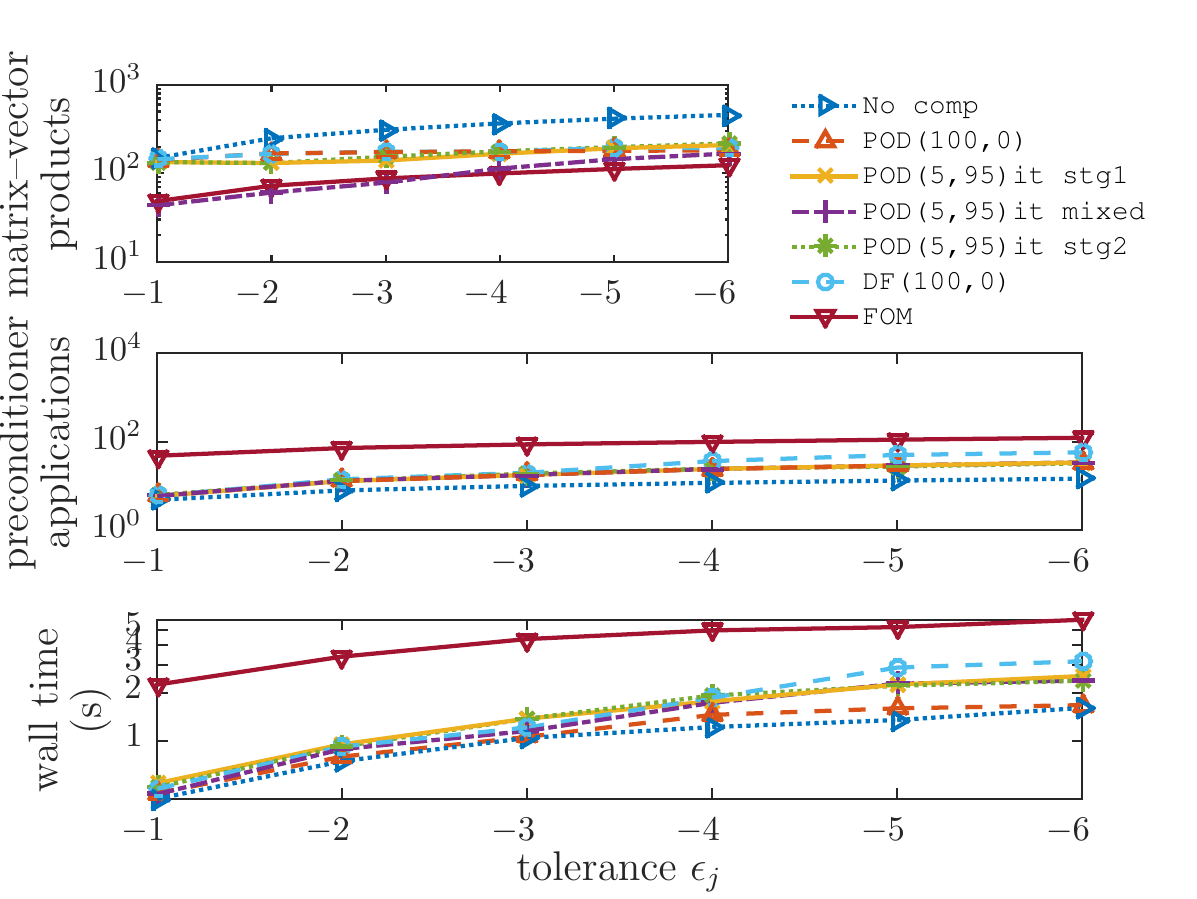} 
                \caption{All methods}
	\label{fig:original_all}
        \end{subfigure}
         \begin{subfigure}[b]{0.48\textwidth}
                \includegraphics[width=\textwidth]{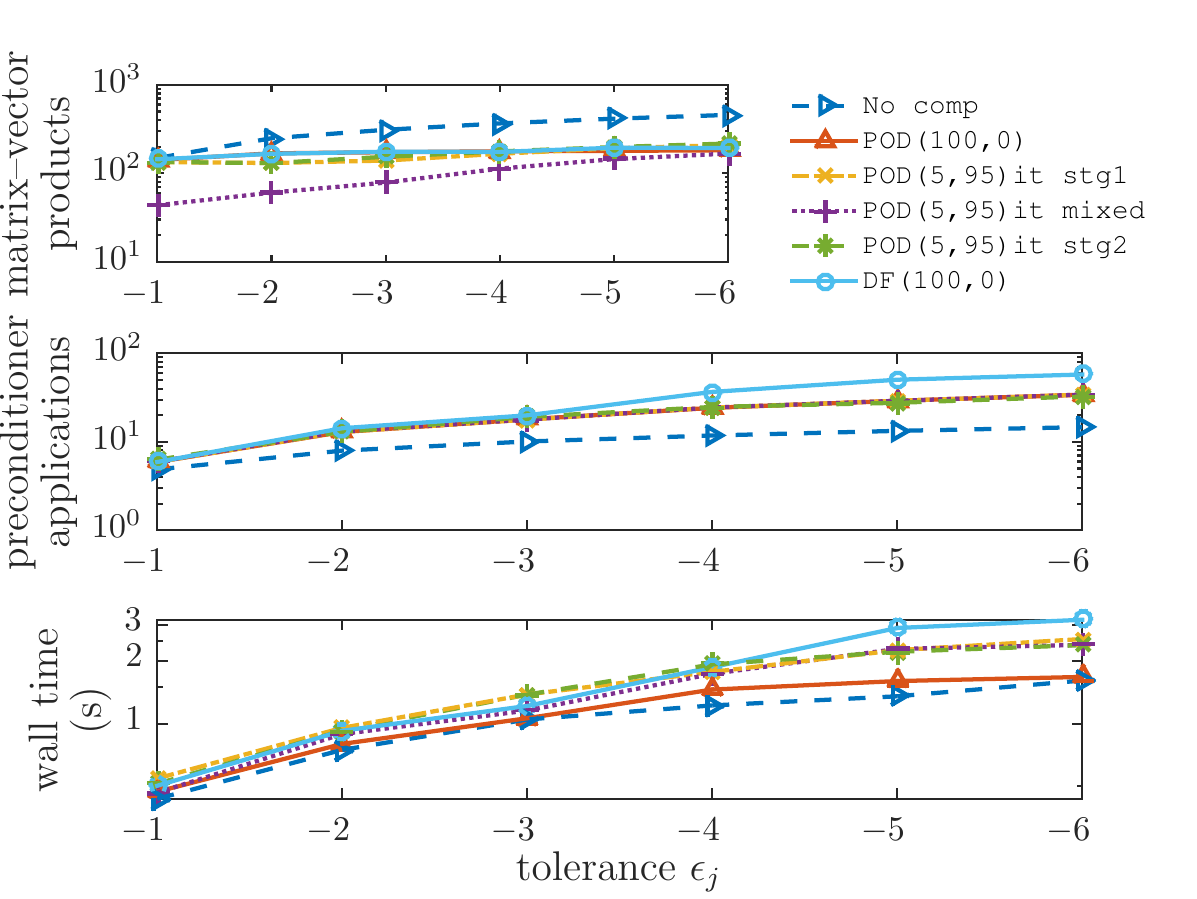} 
                \caption{Recycling methods only}
\label{fig:original_noCGFOM}
        \end{subfigure}
				\caption{Pancake problem: average number of matrix--vector products,
	preconditioner applications, and wall time to compute solutions within
	tolerances $\tolj = 10^{-1}$ through $\tolj = 10^{-6}$.}
\end{figure} 

Figure~\ref{fig:original_all} reports results for all tested methods.  First,
we note that applying the AMG preconditioner is computationally expensive for
this problem, especially relative to matrix--vector products.  Therefore,
there is a strong relationship between the wall-time performance and the
preconditions-application performance of iterative methods for this example.
Next, we note that recycling provides a significant benefit, as applying FOM
without recycling is the slowest method for all tested tolerances. To more
clearly distinguish the differences between recycling methods,
Figure~\ref{fig:original_noCGFOM} reports the same results with the FOM
performance removed.
Here, we see that the no-truncation case yields the best performance as
measured in both wall time and preconditioner applications; however, it
performs the worst in matrix--vector products. This arises from two primary
effects: 1) the preconditioner application is the dominant cost for this
problem, so minimizing the number of stage-3 iterations---which will always
occur by not truncating the augmenting subspace---yields the best wall
time performance, and 2) the problem is small scale, so there is not a
significant penalty to retaining all Krylov vectors.  We also note that the
POD methods (especially the POD(100,0) method) perform similarly to the no
truncation method; this suggests that POD truncation effectively captures the
most important subspace from the set of available vectors.

\begin{figure}[htbp] 
\centering 
\includegraphics[width=0.5\textwidth]{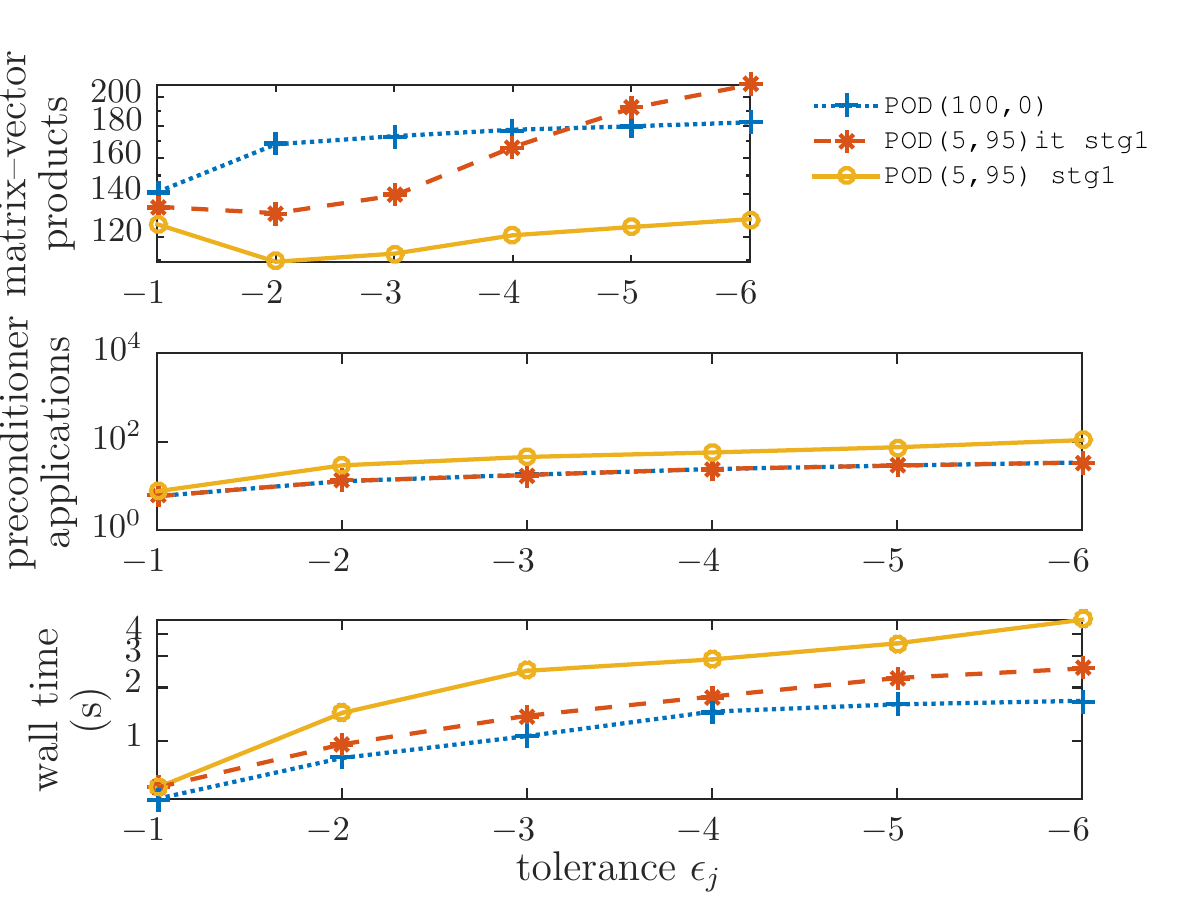}
\caption{Comparison of iterative and non-iterative Stage 3 for the problem 1. }
\label{fig:it_comp} 
\end{figure}

Next, we note that the `inner' iterative method described in Section
\ref{sec:stage3special}---which is used by the three POD(5,95)
methods---produces the same number of preconditioner applications as
POD(100,0). This illustrates that the inner iterative method has successfully
orthogonalized against the entire augmenting subspace.
Figure~\ref{fig:it_comp} illustrates this point further by comparing three
methods: the POD(5,95)it
stg1 method, the POD(5,95) stg1 method (which is identical to POD(5,95)it stg1 except
that it employs
$\stagethreeoptionone=0$ in Algorithm \ref{alg:overall}), and the POD(100,0)
method. While POD(5,95)it stg1 matches the preconditioning applications of
POD(100,0), this comes at the cost of additional matrix--vector products for
stricter tolerances; further these matrix--vector products are generally not applied as
a block as in POD(100,0), which makes them more expensive on average. 
In addition, the cost of repeatedly multiplying vectors by the full augmenting-subspace
basis $\augSpaceBasisj$ within the inner iterative method is not counted
toward matrix--vector products, even though this operation incurs a
non-negligible cost due to the density of $\augSpaceBasisj$.
As a
result, POD(100,0) yields the lowest wall time; we might expect POD(5,95)it
stg1 to produce a lower wall time 
when the dimension of the augmenting subspace is larger.

\begin{figure}[htbp] 
\centering 
\includegraphics[width=0.5\textwidth]{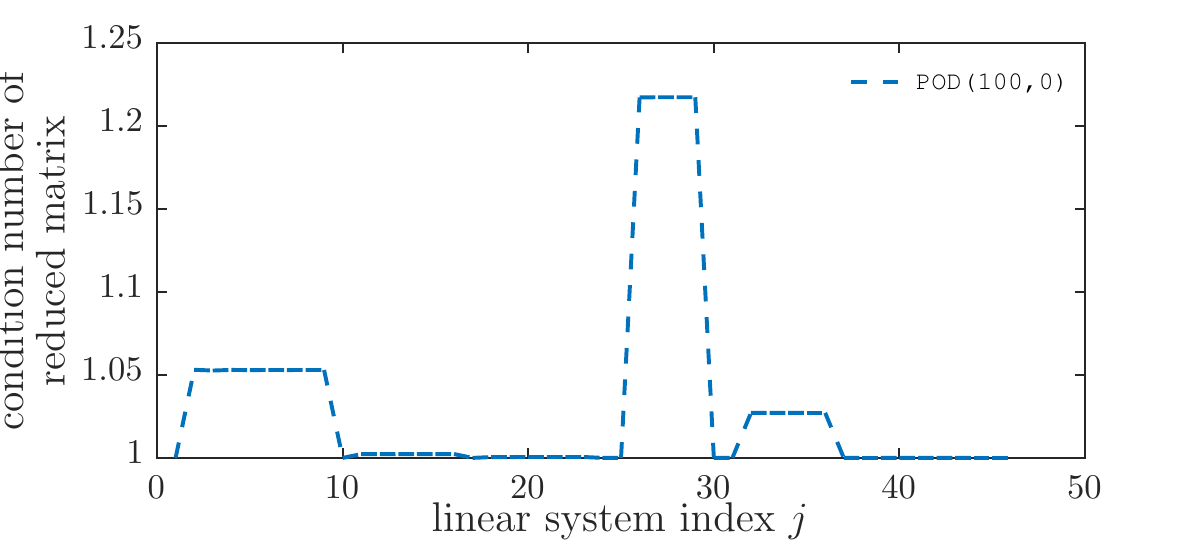}
\caption{Condition number for reduced linear systems.}
\label{fig:reduced_cond} 
\end{figure}

As was shown in Theorem \ref{wellConditioned}, the reduced matrix
$\augSpaceBasisj^T\Aj\augSpaceBasisj$ should be well conditioned if the system
matrices do not vary significantly. Figure~\ref{fig:reduced_cond} illustrates
this effect: the condition number of the reduced system is close to one for all
linear systems. This implies fast convergence of stage 2. Note that this
effect comes `for free' when employing a POD metric of $\A_{\jcompress}$, as the
basis is automatically $\A_{\jcompress}$-orthogonal; this precludes the need
for orthogonalization step \ref{step:enforceOrthog} in Algorithm
	\ref{alg:overall} as was noted in Section \ref{sec:integratePOD}.

\subsection{Method comparison: I-beam problem} \label{sec:methodComparison}

We compare the same methods as in Section \ref{sec:pancakeResults}; the only
modification is that 
all POD(5,95) methods employ parameters $\tolstagetwoj = 10^{-4}\tolj$,
$\innertolstagetwoj = 10^{-2} \tolj$ in Algorithm \ref{alg:overall}.

\begin{figure}[hbt]
\centering
 \begin{subfigure}[b]{0.48\textwidth}
				\includegraphics[width=1.0\textwidth]{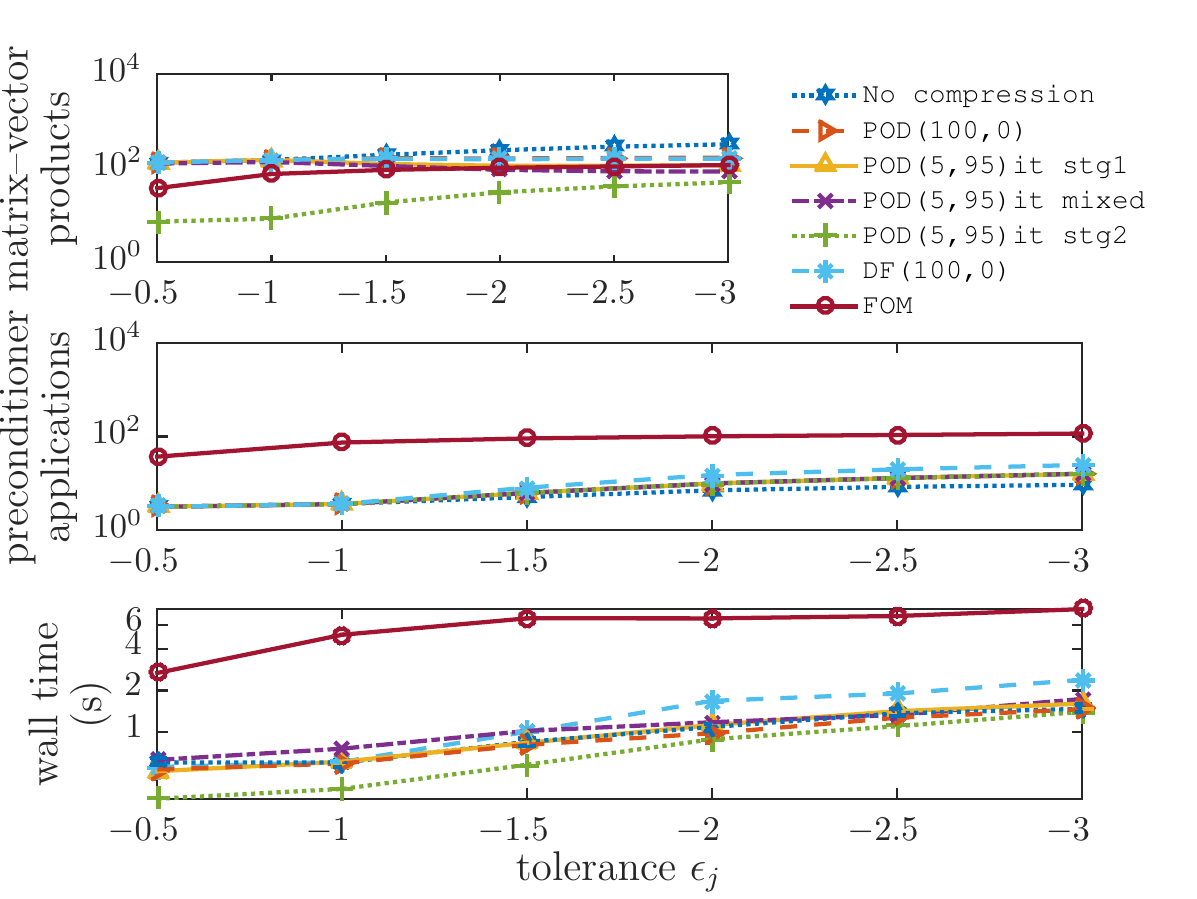} 
				\caption{All methods}
\label{fig:ibeam_zoom}
\end{subfigure}
 \begin{subfigure}[b]{0.48\textwidth}
				\includegraphics[width=1.0\textwidth]{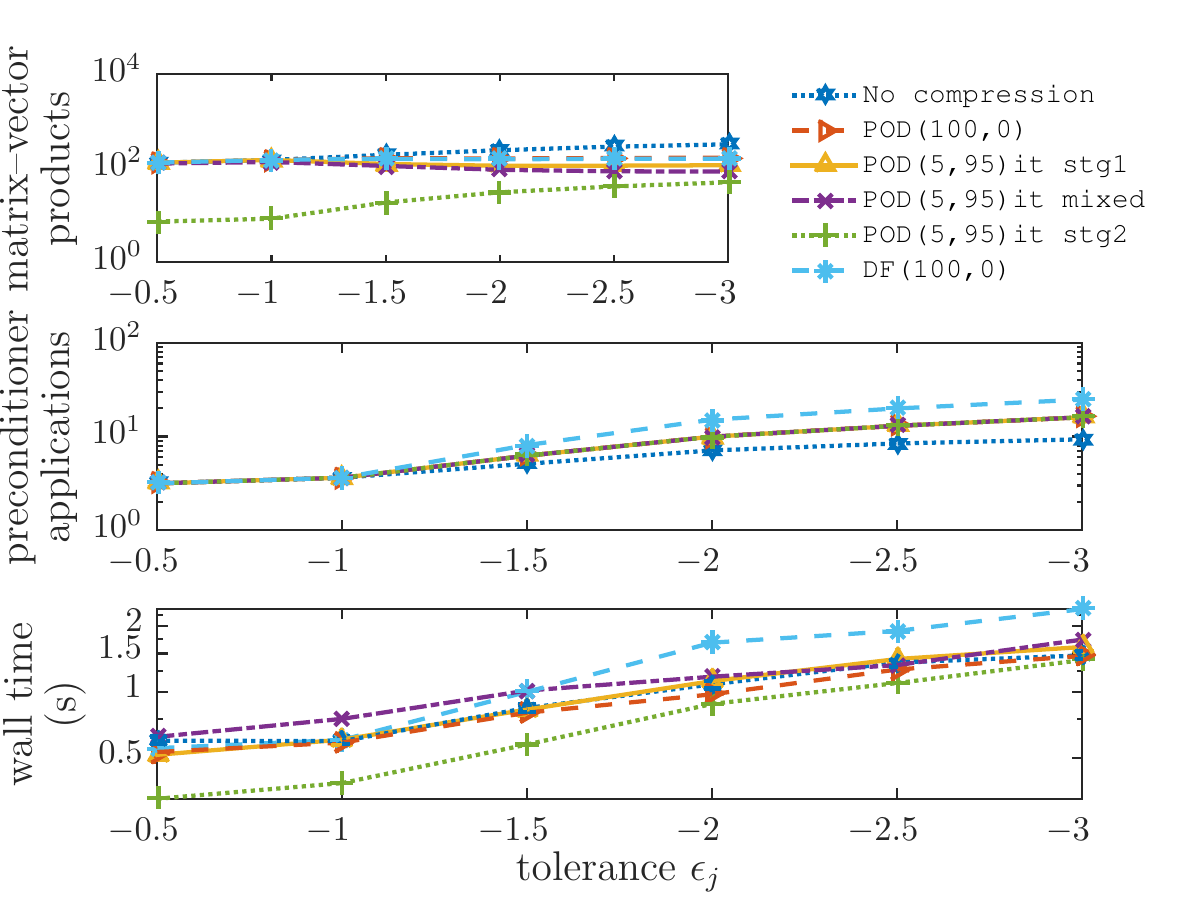}    
				\caption{Recycling methods only}
\label{fig:ibeam_zoomnoCGFOM}
\end{subfigure}
\caption{I-beam problem: average number of matrix--vector products,
	preconditioner applications, and wall time to compute solutions within
	tolerances $\tolj = 10^{-0.5}$ through $\tolj = 10^{-3}$.}
\end{figure}

Figure~\ref{fig:ibeam_zoom} reports results for all tested methods. Again, we
first note that the use of recycling leads to significant improvements, as the
FOM method (without recycling) produces the largest wall time and requires the
largest number of preconditioner applications.
To more easily distinguish the relative merits of the recycling methods,
Figure~\ref{fig:ibeam_zoomnoCGFOM} reports the results without FOM. This
figure illustrates the need for truncation within recycling. As for the pancake
problem, the `no truncation' case minimizes
the number of preconditioner applications; however, the 
lower matrix--vector multiplication cost and the lower overhead in the
stage-3 orthogonalization steps
lead to lower overall costs for several truncation methods.
Figure~\ref{fig:ibeam_zoomnoCGFOM} also illustrates the benefit of using the
hybrid direct/iterative approach to solve over the augmenting subspace, as 
POD(5,95)it stg2---which employs this approach---produces the lowest wall time
and number of matrix--vector products for all tested tolerances. We also note
that all POD-based truncation methods outperform deflation in terms of
preconditioner applications; all POD methods except for POD(5,95)it mixed also
outperform deflation in terms of matrix--vector multiplications and wall time.

\begin{figure}[htbp] 
\centering 
\includegraphics[width=0.5\textwidth]{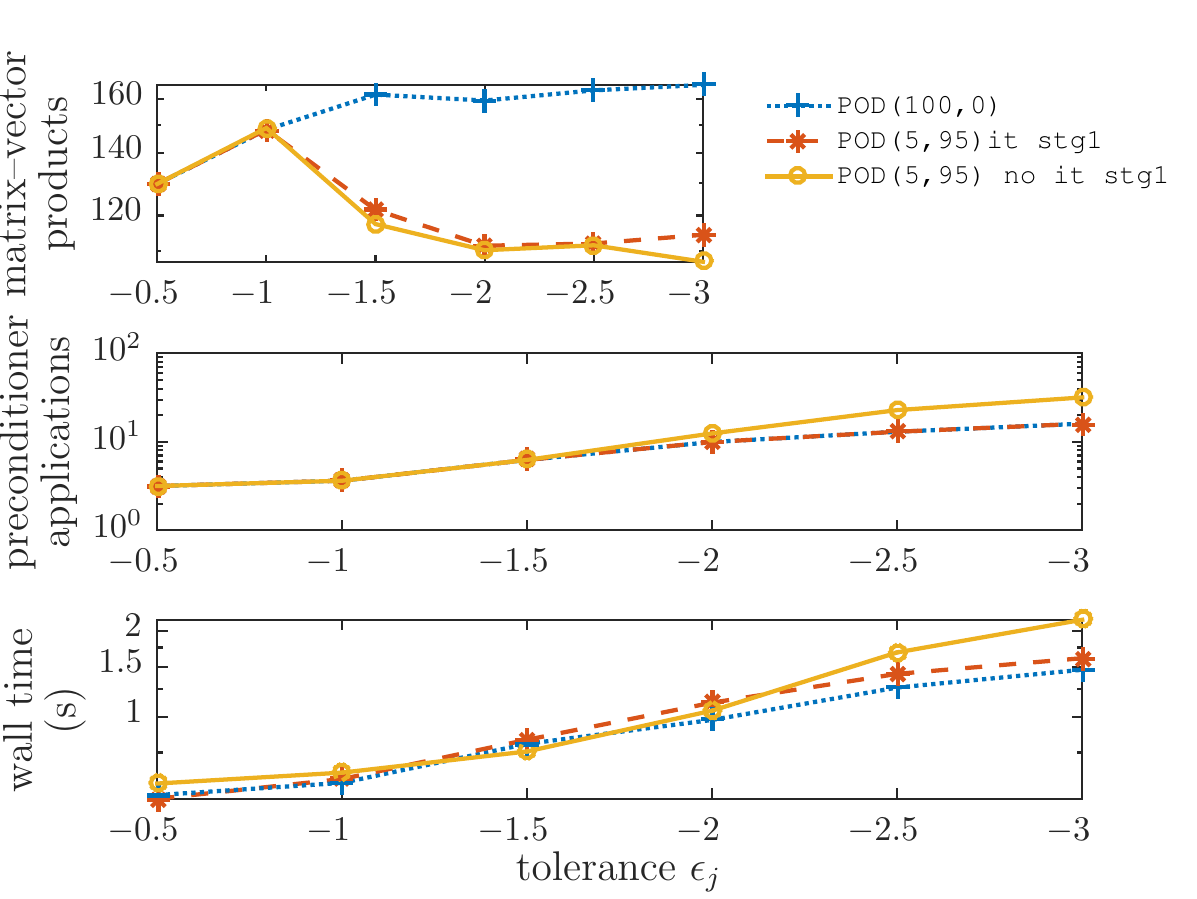}
\caption{I-beam problem: POD-method performance in terms of
average number of matrix--vector products, preconditioner applications, and
wall time to compute solutions within tolerances $\tolj = 10^{-0.5}$ through
$\tolj = 10^{-3}$.}
\label{fig:it_comp_ibeam} 
\end{figure}

Finally, Figure~\ref{fig:it_comp_ibeam} assess the performance of the `inner'
iterative method proposed in Section \ref{sec:stage3special}, which aims to
orthogonalize against the entire augmenting subspace within stage 3 using an
iterative method. These data show that the iterative method in POD(5,95)it
stg1 successfully matches the number of stage-3 iterations (i.e.,
preconditioner applications) as POD(100,0), yet it incurs far fewer
matrix--vector products. In fact, the number of matrix--vector products is
close to that realized by POD(5,95) stg1 without the inner iterative method. 

\subsection{POD-weights experiments}\label{sec:experimentWeights}

This section compares the performance of POD-based truncation methods when
various POD weights are employed (see Section \ref{sec:podweights}).
To assess
this, POD truncation is performed after the solution of the first ten linear
systems, and performance of the resulting truncated augmented subspace is
assessed for the eleventh linear system. That is, we propose computing the truncated augmenting subspace according
to
\begin{equation}\label{eq:augSpaceIsPODTest}
\augSpace_{11} =
\podSubspace{\augSpaceDim_{11}}{\currentDataArg{10}}{\weightVec}{\A_{10}}
\end{equation}
for various choices of $\weightVec$.
We compare the following POD weights proposed in Section
\ref{sec:podweights}:
\begin{enumerate} 
\item \textit{Ideal weights}. This case employs $\weightVec =
\augSpaceWeightsAVecArg{11}
$ as defined in 
Eq.~\eqref{eq:weightsIdealAj}.  Although this weighting scheme is not
practical, it illustrates the best possible choice.
\item \textit{Previous weights}. This case employs $\weightVec =
\augSpaceWeightsVecPrevArg{11}
$ as defined in 
Eq.~\eqref{eq:weightsPrevious}. 
\item \textit{Radial-basis-function weights}. This case employs $\weightVec =
\augSpaceWeightsVecRBFArg{11}(11)
$ as defined in 
Eq.~\eqref{eq:weightsRBF}.  This is the choice used by the previous
experiments in this section.
\end{enumerate}

\begin{figure}[htbp] 
\centering 
         \begin{subfigure}[b]{0.48\textwidth}
				 \captionsetup{width=.8\linewidth}
                \includegraphics[width=\textwidth]{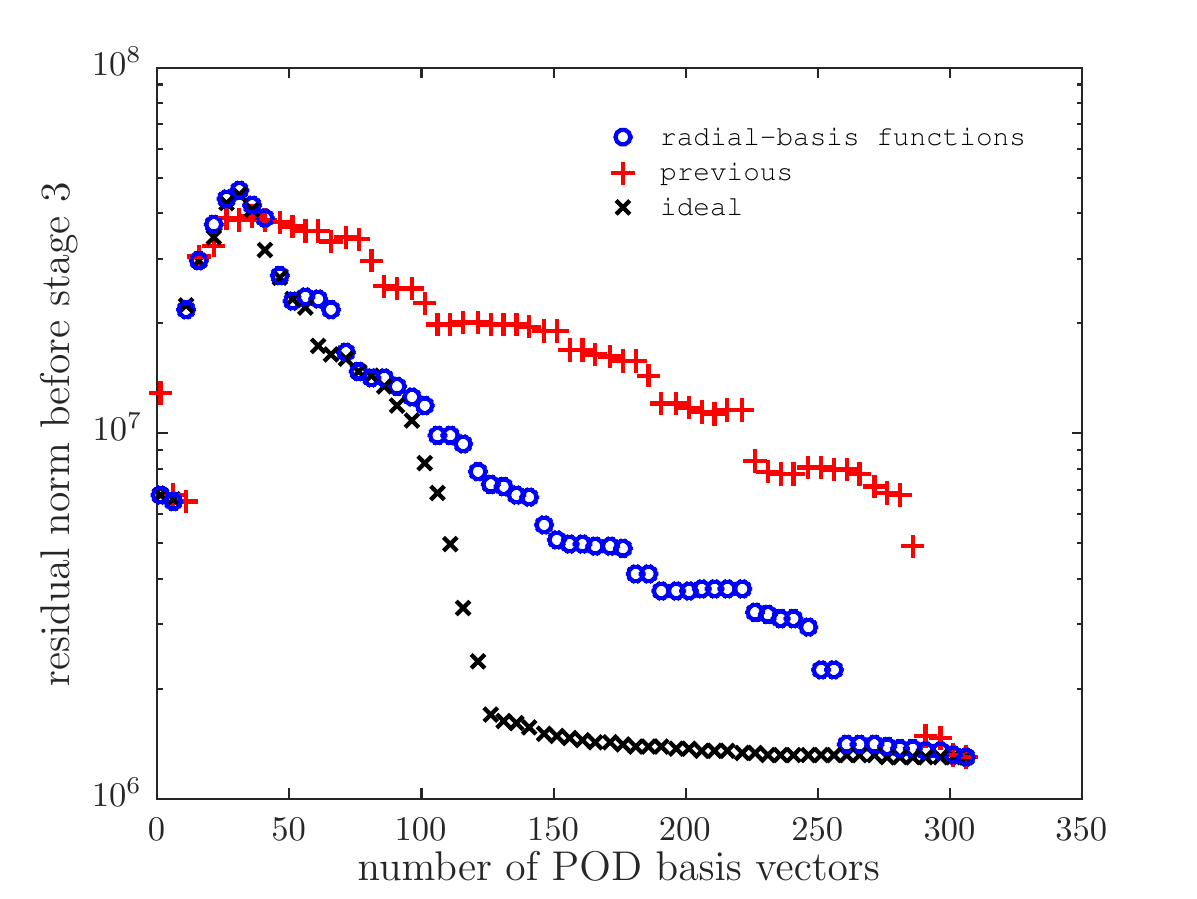} 
								\caption{The residual norm before stage 3 as a function of
								the number of POD vectors $k$.}
\label{fig:weights_comp_pancake11a}
        \end{subfigure}
         \begin{subfigure}[b]{0.48\textwidth}
				 \captionsetup{width=.8\linewidth}
                \includegraphics[width=\textwidth]{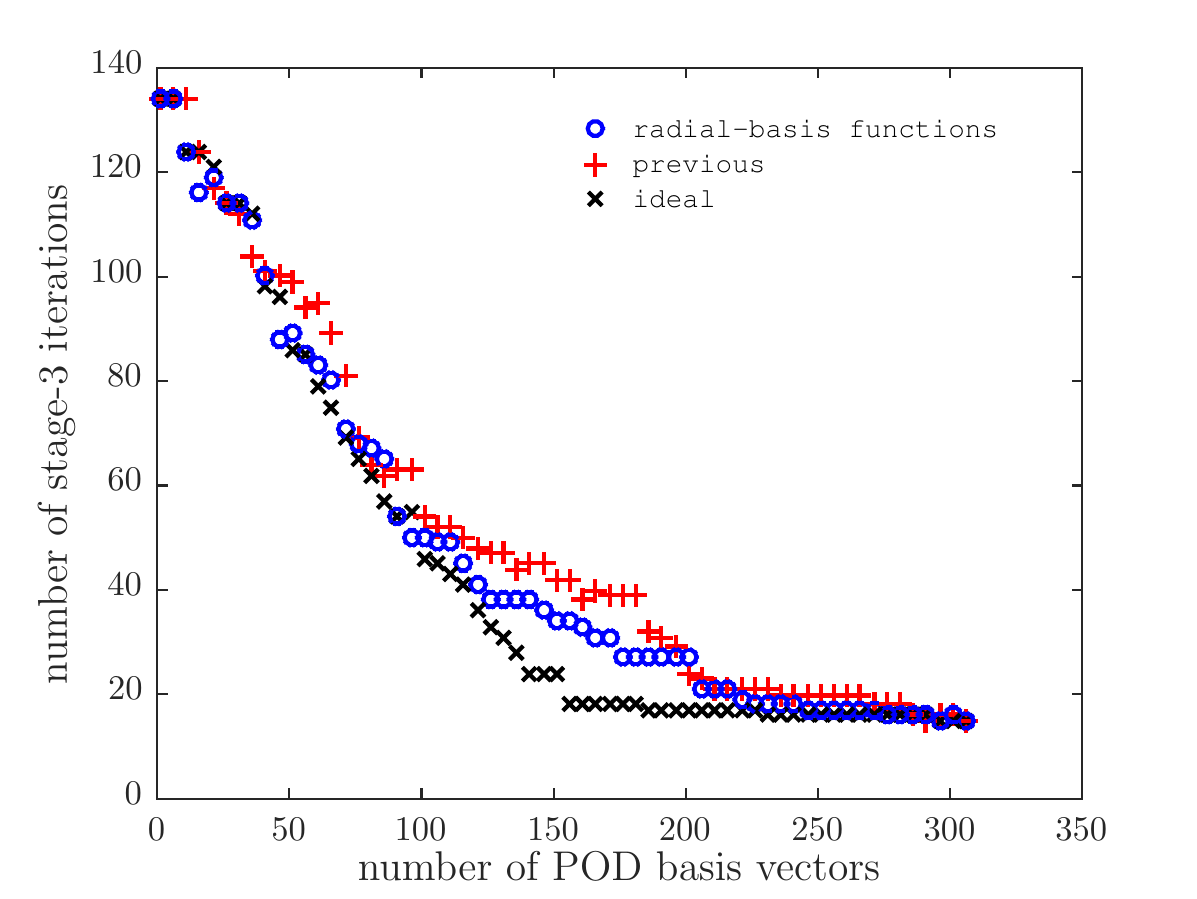} 
                \caption{Number of stage-3 iterations taken 
								to satisfy a tolerance of $\tol_{11} = 10^{-6}$
								as a function of the number of POD vectors $k$.}
\label{fig:weights_comp_pancake11b}
        \end{subfigure}
				\caption{Pancake problem: POD-weights experiments. Results
				correspond to solving linear system eleven after performing POD-based
				truncation for linear system ten using  radial-basis-function weights, previous weights, and ideal weights.}
\end{figure} 
Figure~\ref{fig:weights_comp_pancake11a} illustrates that the ideal weights
minimize the residual after the reduced system is solved, which implies that
the ideal weights lead to a better estimate of the solution in the augmenting
subspace.
Note
that the radial-basis-function weights (which was employed for POD results in
previous sections) yields results that are close to
the ideal case. Figure~\ref{fig:weights_comp_pancake11b} shows that the ideal weights
also minimize the number of stage-3 iterations, while the two other methods
produce a similar number of stage-3 iterations as the ideal-weights case. This
suggests that radial-basis-function weights provide a good approximation of
the ideal weights for both producing an accurate solution in the augmenting
subspace and yielding similar stage-3 convergence. 

\subsection{Output-oriented POD experiments}

This section assesses the performance of output-oriented POD, i.e., when the
metric is set to $\metric = \outputMat^T\outputMat$ as proposed in Section
\ref{sec:podmetric}. For this purpose, we consider a set of $\nout=100$ output
quantities of interest that are random linear functionals of the solution; as
such $\outputMat\in[0,1]^{100\times\ndof}$ with entries drawn from a uniform
distribution in the interval $[0,1]$. In addition to assessing the performance
of the No truncation, DF(100,0), and POD(100,0) as described in Section
\ref{sec:methodComparison}, we also compare the following methods:
\begin{enumerate} 
\item \textit{POD(5,95)}. 
This approach employs POD truncation, places only the dominant POD modes in
the stage-1 basis, and places all post-truncation Krylov vectors in the
stage-1 basis.
Algorithm \ref{alg:overall} parameters are
$\augSpaceThreshold = 200$,
$\newvectorsdirect=1$,  $\stagethreeoptionone=0$, $ \augSpaceBasisjp  =
\podBasis{100}{\currentDataj}{\augSpaceWeightsVecRBF(j - \jcompress )}{\Aj} $
in step \ref{step:compress1}, and $\stageonebasisArg{j+1} =
\left[\left[\augSpaceBasisjp\right]_1\ \cdots\
\left[\augSpaceBasisjp\right]_{5}\right]$ in step
\ref{step:compresslast}.
\item \textit{POD(5,95)it}. 
This approach employs POD truncation, places only the dominant POD modes in
the stage-1 basis, and places all post-truncation Krylov vectors in the
stage-1 basis.
Algorithm \ref{alg:overall} parameters are
$\augSpaceThreshold = 200$,
$\newvectorsdirect=1$,  $\stagethreeoptionone=1$, $ \augSpaceBasisjp  =
\podBasis{100}{\currentDataj}{\augSpaceWeightsVecRBF(j - \jcompress )}{\Aj} $
in step \ref{step:compress1}, and $\stageonebasisArg{j+1} =
\left[\left[\augSpaceBasisjp\right]_1\ \cdots\
\left[\augSpaceBasisjp\right]_{5}\right]$ in step
\ref{step:compresslast}.
Note that $\innertolstagetwoj = 10^{-2}\tolj$.
\item \textit{POD-$C^TC$(100,0)}. Identical to POD(100,0)
except for $ \augSpaceBasisjp  =
\podBasis{100}{\currentDataj}{\augSpaceWeightsVecRBF(j - \jcompress
)}{\outputMat^T\outputMat} $.
\item \textit{POD-$C^TC$(5,95)}. Identical to POD(5.95)
except for $ \augSpaceBasisjp  =
\podBasis{100}{\currentDataj}{\augSpaceWeightsVecRBF(j - \jcompress
)}{\outputMat^T\outputMat} $. 
\item \textit{POD-$C^TC$(5,95)it}. 
Identical to POD(5.95)
except for $ \augSpaceBasisjp  =
\podBasis{100}{\currentDataj}{\augSpaceWeightsVecRBF(j - \jcompress
)}{\outputMat^T\outputMat} $.
\end{enumerate}

To assess the performance of output-oriented POD truncation---which aims to
accurately represent the output quantity of interest---we monitor the error of
the solution in the output-oriented norm $|| \xExactj - \xjk||_{\outputMat^T
\outputMat}$. For the pancake problem, we employ $\tolj=10^{-6}$, while we use
$\tolj=10^{-3}$ for the I-beam problem.
During the execution of the stage-3 algorithm, we track the output-oriented
error and report the average number of matrix--vector products,
preconditioner applications, and wall time for the error to satisfy
$|| \xExactj - \xjk||_{\outputMat^T
\outputMat} < \tolOutput $ for a variety of output-oriented tolerances $\tolOutput$.

\begin{figure}[htbp] 
\centering 
         \begin{subfigure}[b]{0.48\textwidth}
                \includegraphics[width=\textwidth]{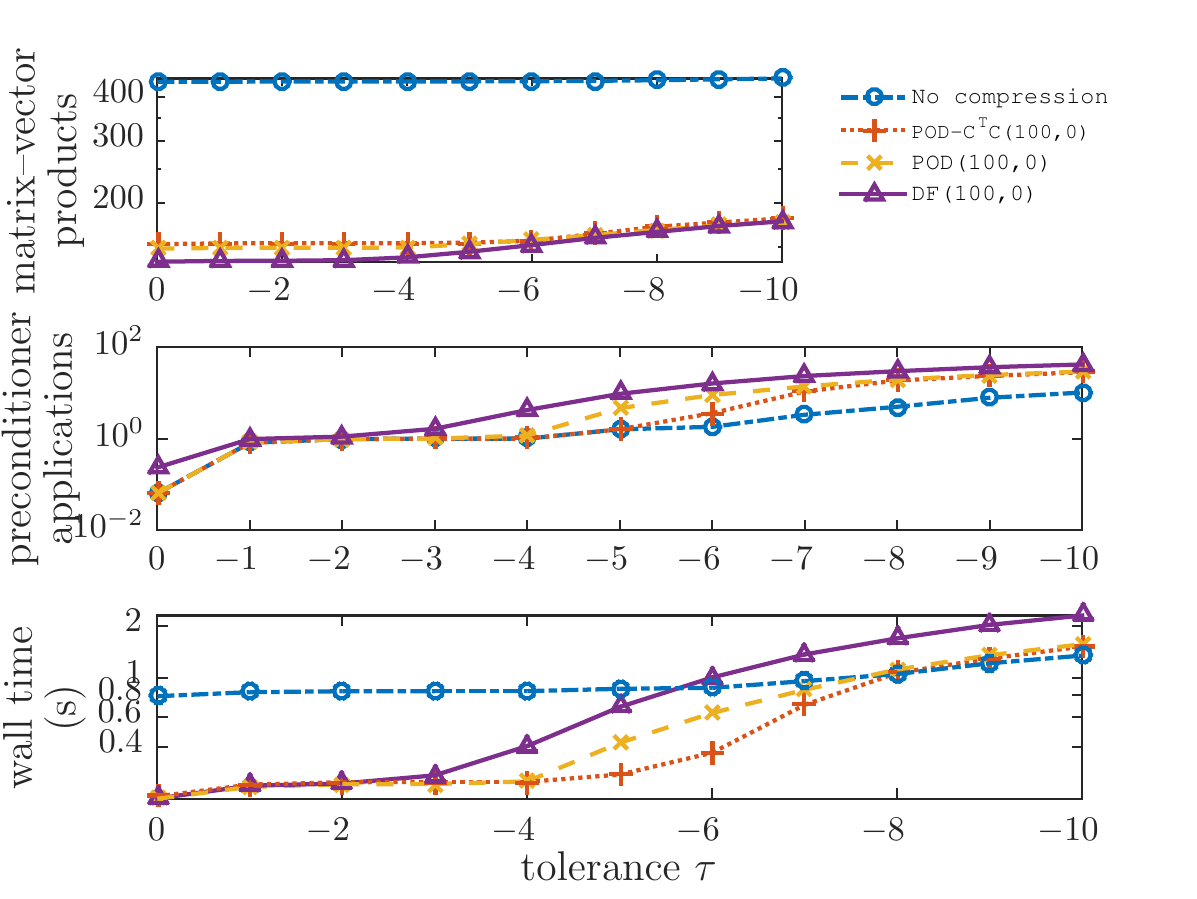} 
                \caption{Pancake problem: stage-1 methods only. }
	\label{fig:goal_pancakea}
        \end{subfigure}
         \begin{subfigure}[b]{0.48\textwidth}
                \includegraphics[width=\textwidth]{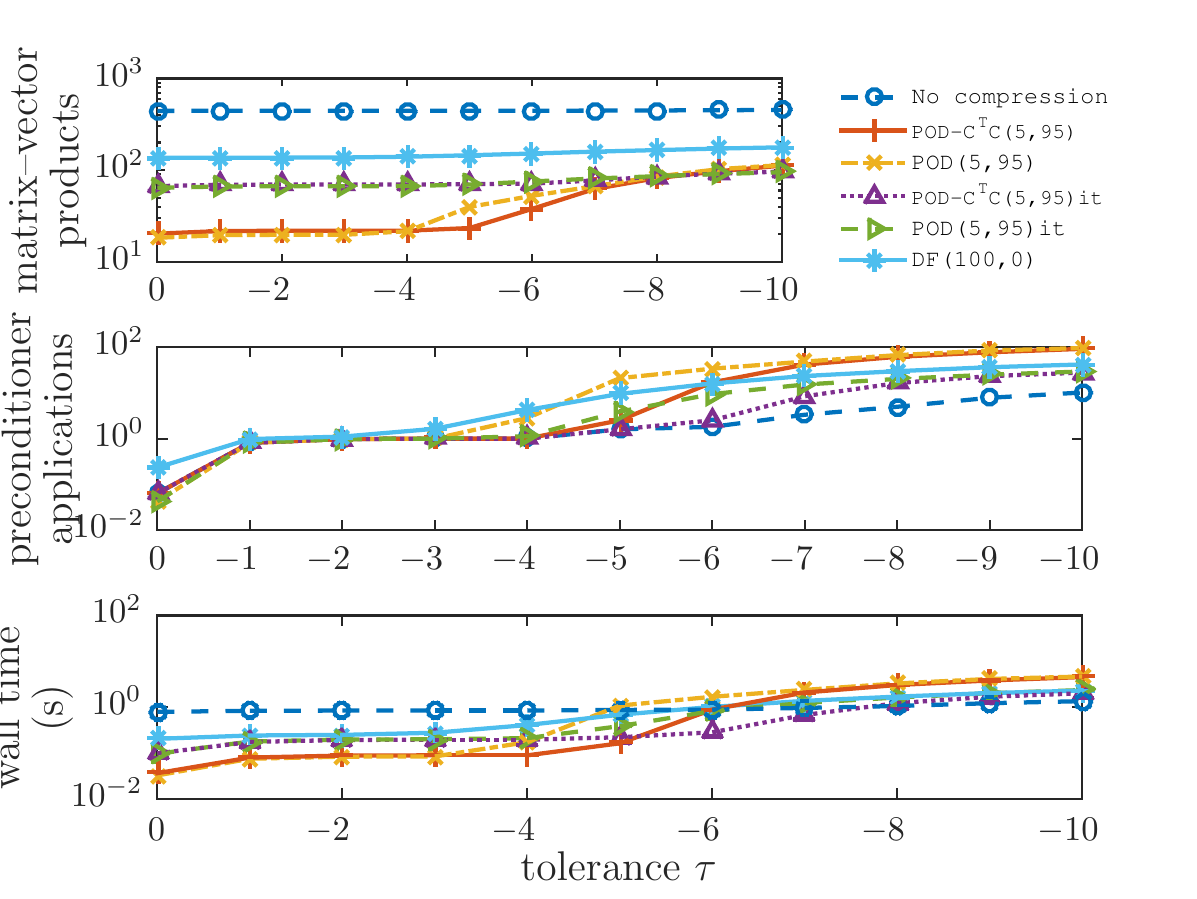} 
                \caption{Pancake problem: stage-1 and stage-2 methods. }
\label{fig:goal_pancakeb}
        \end{subfigure}
         \begin{subfigure}[b]{0.48\textwidth}
                \includegraphics[width=\textwidth]{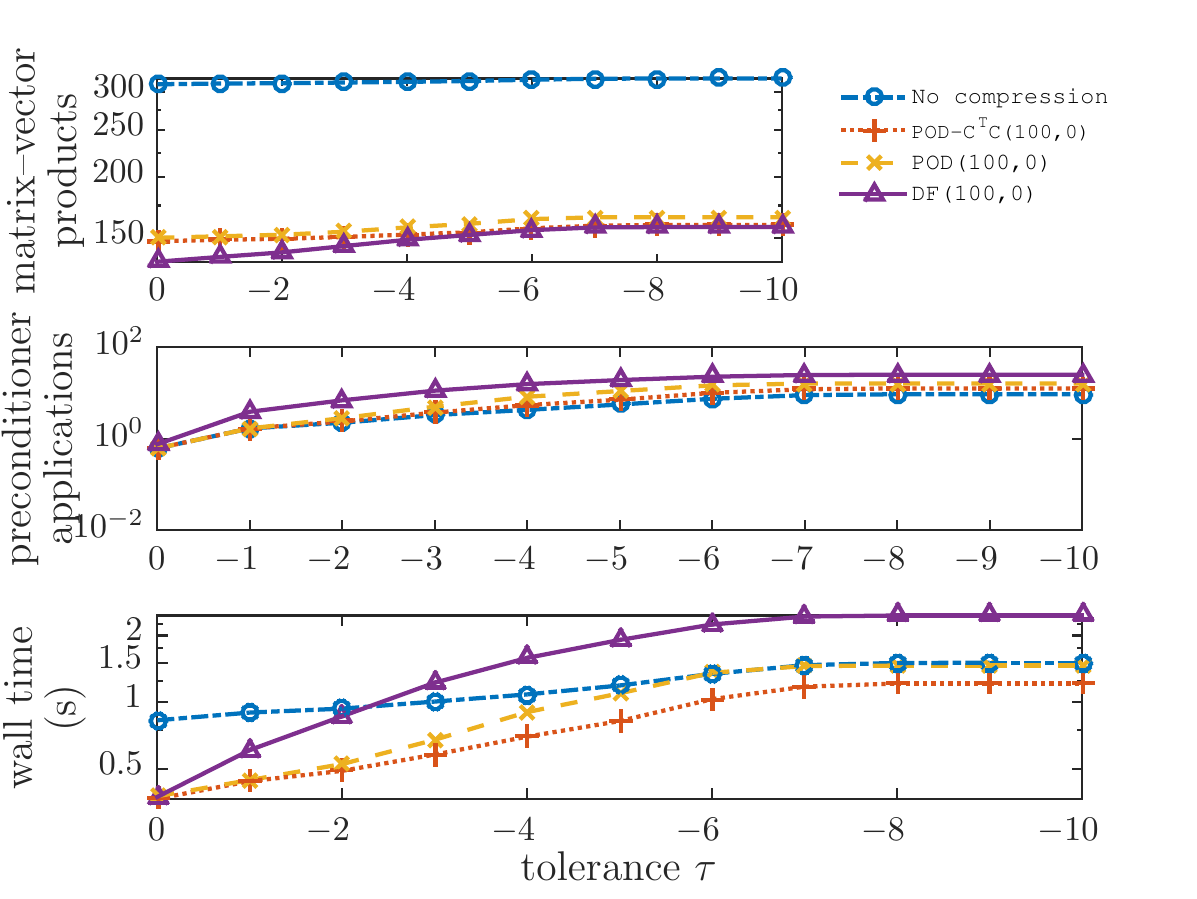} 
								\caption{I-beam problem: stage-1 methods only.}
	\label{fig:goal_ibeama}
        \end{subfigure}
         \begin{subfigure}[b]{0.48\textwidth}
                \includegraphics[width=\textwidth]{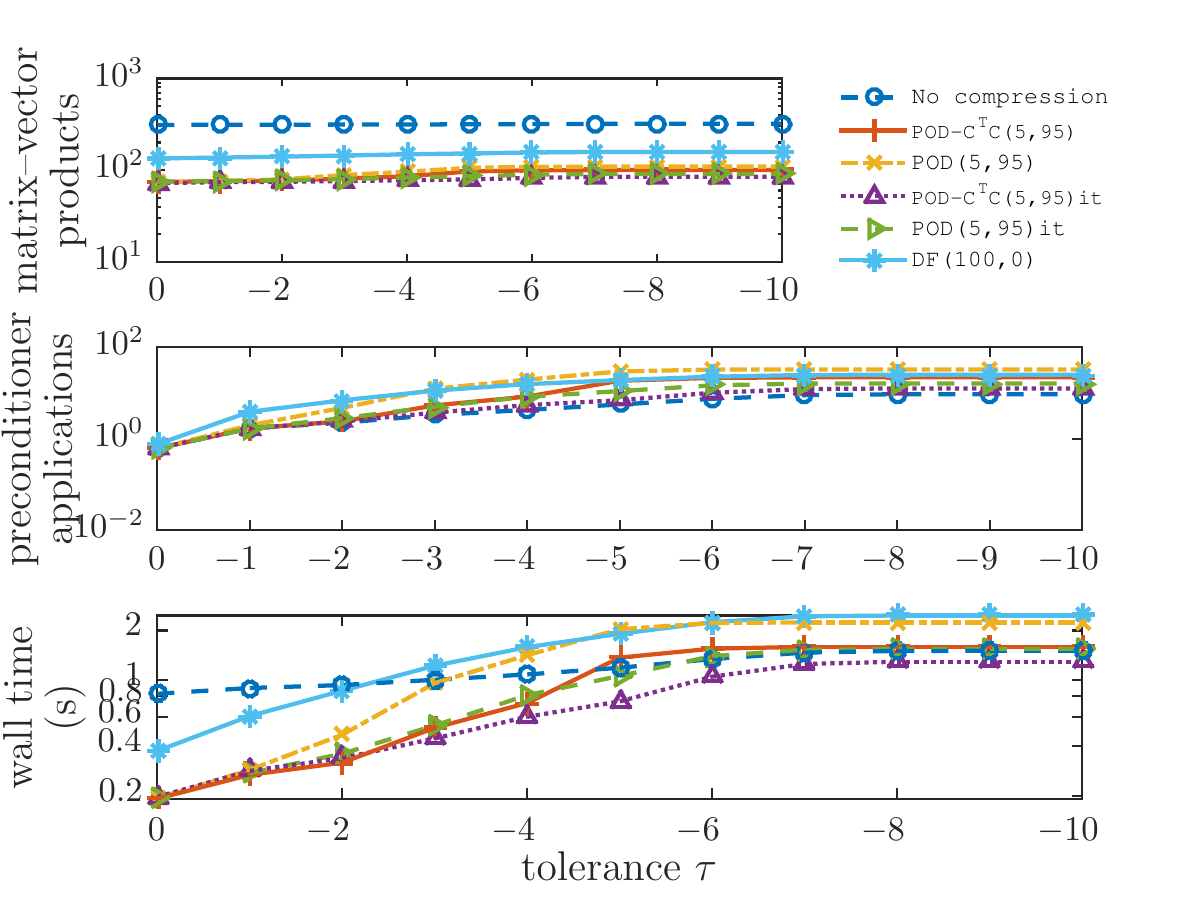} 
                \caption{I-beam problem: stage-1 and stage-2 methods. }
\label{fig:goal_ibeamb}		 
        \end{subfigure}
				\caption{Output-oriented truncation experiments: average number of matrix--vector products,
	preconditioner applications, and wall time to compute solutions as a
	function of the output-oriented tolerance $\tolOutput$. Note that (average) preconditioner
	applications below one indicate that the solution satisfied the specified
	tolerance before entering stage 3.}
\end{figure}

Figures~\ref{fig:goal_pancakea}--\ref{fig:goal_ibeamb} compare the
performances of the assessed methods as a function of the output-oriented
tolerance $\tolOutput$.  First, note that
Figures~\ref{fig:goal_pancakea}--\ref{fig:goal_pancakeb} show that POD-based
truncation methods produce the lowest wall times for \textit{inexact
output-oriented tolerances} for the pancake problem. These plots also
illustrate the benefit of output-oriented metric over the $\Aj$-metric:
Figure~\ref{fig:goal_pancakea} shows the superior performance of the former
for modest output-oriented tolerances, i.e., those between $\tolOutput = 10^{-4}$ and
	$10^{-7}$.  Figure~\ref{fig:goal_pancakeb} shows that the POD(5,95) methods
	without the inner iterative method described in
	Section~\ref{sec:stage3special} are faster for inexact output-oriented
	tolerances; this occurs because very few stage-3 iterations are required for
	convergence in the output-oriented norm to these tolerances. In contrast,
	for stricter tolerances $\tolOutput \le 10^{-6}$, the inner iterative modification
		(or no truncation) yields superior performance. 

Figures~\ref{fig:goal_ibeama}--\ref{fig:goal_ibeamb} present results for the
I-beam problem.  These results illustrate the advantage to using the
output-oriented metric  for most values of the output-oriented tolerance
$\tolOutput$,
as POD-$C^TC$(100,0) produces the lowest wall time in
Figure~\ref{fig:goal_ibeama} and the POD-$C^TC$(5,95)it
method produces the lowest wall time for $\tolOutput \le 10^{-3}$.

\begin{figure}
 \renewcommand{\baselinestretch}{1}
\small\normalsize
\centering
	\includegraphics[width=0.5\textwidth]{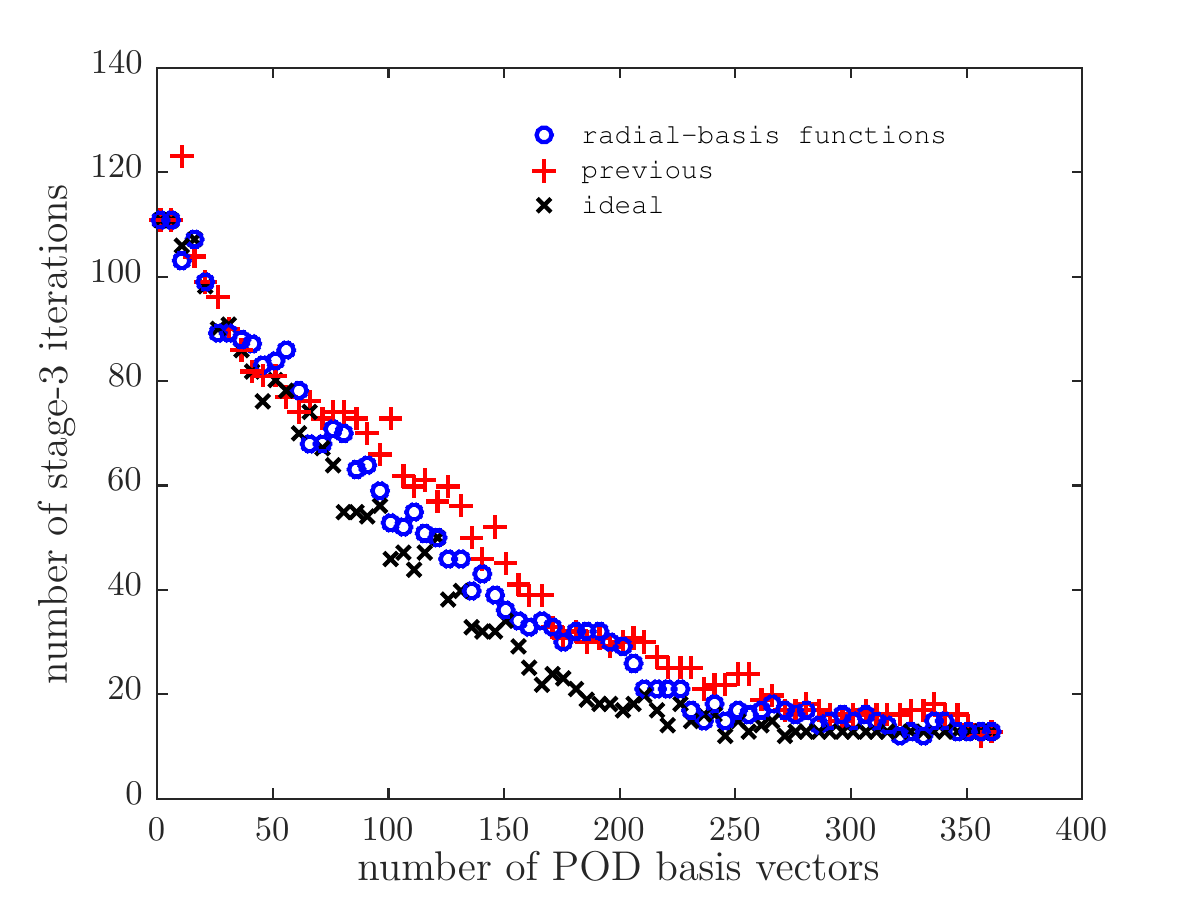}
\caption{Number of stage-3 iterations taken to satisfy an output-oriented
tolerance of $\tolOutput = 10^{-6}$ for the pancake problem as a function of
the augmenting-subspace dimension $k$ for radial-basis-function weights, previous weights, and ideal weights.}
\label{fig:goal_weights}
\end{figure}

Finally, we repeat the POD-weights experiments discussed in
Section~\ref{sec:experimentWeights} for the output-oriented POD-based
truncation method applied to the pancake problem. The experiments measure the
number of stage-3 iterations required to converge to a output-oriented
tolerance of $\tolOutput = 10^{-6}$ as a function of the dimension of the
augmenting subspace. We employ a stage-3 tolerance of $\tolj = 10^{-8}$ to
accrue the search directions for the first 10 linear systems and subsequently
perform truncation using the output-oriented POD metric with different POD
weights. The resulting stage-3 iterations are for the eleventh linear system.  As
we observed for the $\Aj$-metric, Figure~\ref{fig:goal_weights} shows that the
ideal weights yield the best performance, as they minimizing the number of
stage-3 iterations required for convergence in the output-oriented norm, while
the other weighting schemes closely follow the ideal weights. 


\section{Conclusions}\label{sec:conclusions}

This work has proposed a novel strategy for Krylov-subpace recycling inspired
by goal-oriented proper orthogonal decomposition (POD). We performed analyses
that expose the close connection between model reduction and Krylov-subspace
recycling, proposed specific goal-oriented POD ingredients for truncating
previous Krylov vectors, and developed a new `three-stage' algorithm that
employs a hybrid direct/iterative approach for efficiently solving over the
sum of augmenting and Krylov subspaces. Results on several solid-mechanics
problems highlighted the benefits of the new contributions, especially for
efficiently computing quantities of interest to modest tolerances.

\section*{Acknowledgments}\label{sec:acknowledgments}

We thank Paul Tsuji for his contributions to obtaining the linear systems used
in the numerical experiments, as well as Kendall Pierson for his assistance in
using Adagio. 
We also thank Howard Elman for both insightful conversations
and his support in establishing this collaboration. We are also grateful to Eric de
Sturler and Michael Parks for fruitful conversations an invaluable suggestions.
Sandia National Laboratories is a
multi-program laboratory managed and operated by Sandia Corporation, a wholly
owned subsidiary of Lockheed Martin Corporation, for the U.S.\ Department of
Energy's National Nuclear Security Administration under contract
DE-AC04-94AL85000. The content of this publication does not necessarily
reflect the position or policy of any of these institutions, and no official
endorsement should be inferred. 

\appendix

\section{POD computation}\label{app:podcomputation}

This section describes two techniques for computing a POD basis using given
snapshots, weights, and a pseudometric (see Ref.~\cite{carlbergCpodJour} for
additional details). Algorithm \ref{alg:PODEVD} describes the first technique,
which is based on the eigenvalue decomposition and is equivalent to the
well-known ``method of snapshots'' \cite{sirovich1987tad}. Algorithm
\ref{alg:PODSVD} reports the second case, which is based on the singular value
decomposition (SVD) and is more appropriate when the symmetric factorization $\metric =
(\metric^{1/2})^T\metric^{1/2}$
is available, where $\metric^{1/2}$ need not be upper triangular; this
approach leads to a more well-conditioned linear system.  Note that
Algorithms \ref{alg:PODEVD} and \ref{alg:PODSVD} produce equivalent POD bases
(in exact arithmetic) and differ only in their first two steps.
\begin{algorithm}[htbp]
\caption{\podevdalgorithmname}
\begin{algorithmic}[1]\label{alg:PODEVD}
\REQUIRE snapshot matrix $\snapshotMat\in\RR{\ndof\times \nsnapshots}$, weights
$\weightVec\in\RR{\nsnapshots}$, 
pseudometric $\metric\in\spsd{ \ndof}$, and energy criterion
$\energyCrit\in[0,1]$
\ENSURE  POD basis $\podBasis{\augSpaceDim}{\snapshotMat}{\totuple{\weight}{\nsnapshots}}{\metric}$
\STATE \label{step:redK}$\bar \metric
=\diag\totuple{\weight}{\nsnapshots} \snapshotMat^T \metric \snapshotMat\diag\totuple{\weight}{\nsnapshots}$
\STATE \label{step:EVP} Compute symmetric eigenvalue decomposition $\bar
\metric =\mathbf V\boldsymbol \Sigma^2 \mathbf V ^T$
\STATE \label{step:TruncateEVD}Choose dimension of truncated basis
$\augSpaceDim = \min_{i\in\feasibleSet}i$ with $ \feasibleSet\defeq
\{i\innat{\text{rank}(\bar\metric)}\ |\
\sum_{k=1}^i\sigma_k^2/\sum_{\ell=1}^\nsnapshots\sigma_\ell^2\geq \energyCrit\}$
\STATE \label{step:Truncate2EVD}$
\podBasis{\augSpaceDim}{\snapshotMat}{\weightVec}{\metric}
=\snapshotMat\diag\totuple{\weight}{\nsnapshots}\left[\frac{1}{\sigma_1}\mathbf v_1\ \ldots\
\frac{1}{\sigma_{\augSpaceDim}}\mathbf v_{\augSpaceDim}\right]$, where 
$\boldsymbol\Sigma\defeq\diag(\sigma_1,\ldots,\sigma_{\nsnapshots})$
\end{algorithmic}
\end{algorithm}
\begin{algorithm}[htbp]
\caption{\podsvdalgorithmname}
\begin{algorithmic}[1]\label{alg:PODSVD}
\REQUIRE snapshot matrix $\snapshotMat\in\RR{\ndof\times \nsnapshots}$, weights
$\weightVec\in\RR{\nsnapshots}$, 
pseudometric factor $\metric^{1/2}$ such that
$\metric=(\metric^{1/2})^T\metric^{1/2}\in\spsd{\ndof}$, and energy criterion
$\energyCrit\in[0,1]$

\ENSURE  POD basis $\podBasis{\augSpaceDim}{\snapshotMat}{\totuple{\weight}{\nsnapshots}}{\metric}$
\STATE \label{step:transform1}$\bar \snapshotMat=\metric^{1/2}\snapshotMat \diag\totuple{\weight}{\nsnapshots}$
\STATE\label{step:SVD} Compute thin singular value decomposition $\bar
\snapshotMat= \mathbf U
\boldsymbol\Sigma \mathbf V^T $
\STATE \label{step:Truncate}Choose dimension of truncated basis
$\augSpaceDim = \min_{i\in\feasibleSet}i$ with $ \feasibleSet\defeq
\{i\innat{\text{rank}(\bar\snapshotMat)}\ |\
\sum_{k=1}^i\sigma_k^2/\sum_{\ell=1}^\nsnapshots\sigma_\ell^2\geq \energyCrit\}$
\STATE \label{step:Truncate2}$
\podBasis{\augSpaceDim}{\snapshotMat}{\weightVec}{\metric}
=\snapshotMat\diag\totuple{\weight}{\nsnapshots}\left[\frac{1}{\sigma_1}\mathbf v_1\ \ldots\
\frac{1}{\sigma_{\augSpaceDim}}\mathbf v_{\augSpaceDim}\right]$, where 
$\boldsymbol\Sigma\defeq\diag(\sigma_1,\ldots,\sigma_{\nsnapshots})$
\end{algorithmic}
\end{algorithm}

\section{Proofs}\label{app:proofs}

This section provides proofs for the theoretical results presented in this
paper.

\begin{proofofthing}[Theorem \ref{Anorm}]
We begin by decomposing the centered exact solution 
as $\xExactj-\xGuessj=\xExactjIn + \xExactjOrth$, 
where $\xExactjIn = \projOrth{\Aj}{\xExactj-\xGuessj}{\currentDataSpacejm} =
\currentDatajm \augSpaceWeightsAjVec$ (from Eqs.~\eqref{eq:sequence} and
\eqref{eq:algebraicProjectorGen})
and $\xExactjOrth = \projError{\Aj}{\xExactj-\xGuessj}{\currentDataSpacejm}$.
We can then bound the error as
\begin{align} 
\normAj{\xExactj- \projOrth{\Aj}{\xExactj}{\xGuessj+\augSpacej}}
&=
\normAj{\xExactj-\xGuessj-
\projOrth{\Aj}{\xExactj-\xGuessj}{\augSpacej}}
=
\normAj{\xExactjIn + \xExactjOrth-
\projOrth{\Aj}{\xExactjIn}{\augSpacej}}\\
&=
\normAj{
\projError{\Aj}{\currentDatajm \augSpaceWeightsAjVec}{\augSpacej}}+\normAj{\xExactjOrth}\\
\label{eq:boundState}&\leq
\currentSpaceDimjmhalf\sqrt{\sum_{i=1}^\currentSpaceDimjm
\normAj{\projError{\Aj}{\augSpaceWeightsAj_i\currentDataveci{i}}{\augSpacej}}^2}+\underbrace{\normAj{\xExactjOrth}}_{\text{(I)}},
\end{align} 
where we have used $\projOrth{\Aj}{\xExactjOrth}{\augSpacej}=0$ (because
$\projOrth{\Aj}{\xExactjOrth}{\currentDataSpacejm}=0$ and
$\augSpacej\subseteq\currentDataSpacejm$), the triangle
inquality, and the norm-equivalence relation $\|\x\|_1\leq
n^{1/2}\|\x\|_2$.
By comparing Eqs.~\eqref{eq:podOpt} and \eqref{eq:boundState}, using
monotinicity of the square root function, and noting that term (I) is
independent of the subspace $\augSpacej$, it is clear
that the POD subspace minimizes (over all $\augSpaceDim$-dimensional subspaces of
$\currentDataSpacejm$) an 
upper bound for $\normAj{\xExactj- \projOrth{\Aj}{\xExactj}{\xGuessj+\augSpacej}}$ if $\metric = \Aj$,
$\snapshot_i = \currentDataveci{i}$ and  $\weight_i = \augSpaceWeightsAj_i$ for
$i\innatseq{\nsnapshots}$ with $\nsnapshots = \currentSpaceDimjm$.
\end{proofofthing}

\begin{proofofthing}[Theorem \ref{outNorm}]
We again decompose the centered exact solution $\xExactj-\xGuessj$ according
to
 \begin{equation} \label{eq:exactWeights}
 \xExactj-\xGuessj=
 \xExactjCIn + \xExactjCOrth
 =
 \xExactjCInWin + \xExactjCInWout + \xExactjCOrth.
  \end{equation} 
Here,
$\xExactjCIn = \projOrth{\I}{\xExactj-\xGuessj}{\range{\outputMat^T}}$,
$\xExactjCOrth = \projError{\I}{\xExactj-\xGuessj}{\range{\outputMat^T}}$,
$\xExactjCInWin =
\projOrth{\outputMat^T\outputMat}{\xExactjCIn}{\currentDataSpacejm} =
\currentDatajm\augSpaceWeightsCtCVec$ (from
Eq.~\eqref{eq:algebraicProjectorPseudoGen}),
and $\xExactjCInWout =
\projError{\outputMat^T\outputMat}{\xExactjCIn}{\currentDataSpacejm}$.
We can then bound the error as
\begin{align} 
\|\outputfun{\xExactj}- \outputfun{\projOrth{\Aj}{\xExactj}{\xGuessj+\augSpacej}}\|_2 &= 
\|\outputMat(\xExactj-\xGuessj) -
\outputMat(\projOrth{\Aj}{\xExactj-\xGuessj}{\augSpacej})\|_2
=
\normOutput{\xExactj-\xGuessj -
\projOrth{\Aj}{\xExactj-\xGuessj}{\augSpacej}}\\
&\leq \normOutput{\xExactj-\xGuessj -
\projOrth{\outputMat^T\outputMat}{\xExactj-\xGuessj}{\augSpacej}} +
\normOutput{(\projOrthNo{\Aj}{\augSpacej} -
\projOrthNo{\outputMat^T\outputMat}{\augSpacej})(\xExactj-\xGuessj)}\\
&\leq \normOutput{\xExactj-\xGuessj -
\projOrth{\outputMat^T\outputMat}{\xExactj-\xGuessj}{\augSpacej}} +
 \normOutput{\projOrthNo{\Aj}{\augSpacej} -
\projOrthNo{\outputMat^T\outputMat}{\augSpacej}}\normOutput{\xExactj-\xGuessj}\\
&= \normOutput{\xExactjCInWin + \xExactjCInWout 
 -
\projOrth{\outputMat^T\outputMat}{\xExactjCInWin}{\augSpacej}} +
\normOutput{\projOrthNo{\Aj}{\augSpacej} -
\projOrthNo{\outputMat^T\outputMat}{\augSpacej}}\normOutput{\xExactjCInWin + \xExactjCInWout}\\
&= \normOutput{
\projError{\outputMat^T\outputMat}{\xExactjCInWin}{\augSpacej}}  + \normOutput{\xExactjCInWout }+
\normOutput{\projOrthNo{\Aj}{\augSpacej} -
\projOrthNo{\outputMat^T\outputMat}{\augSpacej}}\normOutput{\xExactjCInWin + \xExactjCInWout}\\
\label{eq:boundOutput}
&\leq
\currentSpaceDimjmhalf\sqrt{\sum_{i=1}^\currentSpaceDimjm
\normOutput{\projError{\outputMat^T\outputMat}{[\augSpaceWeightsCtC]_i\currentDataveci{i}}{\augSpacej}}^2}
+ \underbrace{\normOutput{\xExactjCInWout }
+
\spaceBound\normOutput{\currentDatajm\augSpaceWeightsCtCVec}}_{\text{(I)}}
\end{align} 
Here, we have used $\outputMat\xExactjCOrth=0$ and
$\projOrth{\outputMat^T\outputMat}{\xExactjCInWout}{\augSpacej}=0$. Also,
$\normOutput{\A}=\sup_{\x\neq 0}\normOutput{\A\x}/\normOutput{\x}$ with
$\A\in\RR{\ndof\times\ndof}$ is an induced matrix norm and $\spaceBound \defeq
\sup_{\augSpace\in\grassman{\augSpaceDim}{\ndof}}\normOutput{\projOrthNo{\Aj}{\augSpace}
- \projOrthNo{\outputMat^T\outputMat}{\augSpace}}$.
Within inequality \eqref{eq:boundOutput}, term (I) is independent of the
subspace $\augSpace$; as such, a comparison of Eqs.~\eqref{eq:podOpt} and
\eqref{eq:boundOutput} reveals that the POD subspace minimizes (over all
subspaces $\augSpace$) an upper bound
for $\|\outputfun{\xExactj}- \outputfun{\projOrth{\Aj}{\xExactj}{\xGuessj+\augSpacej}}\|_2$ if $\metric =
	\outputMat^T\outputMat$, 
	$\snapshot_i = \currentDataveci{i}$ and $\weight_i =
	[\augSpaceWeightsCtC]_i$ for $i\innatseq{\nsnapshots}$ with
	$\nsnapshots = \currentSpaceDimjm$ under the stated
	assumptions. 
\end{proofofthing}

\begin{proofofthing}[Theorem \ref{weightsDifference}]
We prove the result in Eq.~\eqref{eq:weightsDiffAj};
Eq.~\eqref{eq:weightsDiffCtC} follows trivially. First note that
$\currentDatajm\left(\augSpaceWeightsAjVec - \augSpaceWeightsVecPrev\right) =
\projOrthNo{\Aj}{\currentDataSpacejm}{(\xExactj - \xGuessj)} -
\projOrthNo{\Ajm}{\currentDataSpacejm}{(\xExactjm - \xGuessjm)}$ from
Eqs.~\eqref{eq:algebraicProjectorGen}, \eqref{eq:weightsIdealAj}, and
\eqref{eq:weightsPrevious}. Then, we have
\begin{align*}
\|\currentDatajm\left(\augSpaceWeightsAjVec - \augSpaceWeightsVecPrev\right)\|
 \leq &
\|\projOrthNo{\Aj}{\currentDataSpacejm}{(\xExactj - \xGuessj)} - 
\projOrthNo{\Ajm}{\currentDataSpacejm}{(\xExactj - \xGuessj)}\| +
\|\projOrthNo{\Ajm}{\currentDataSpacejm}{\left(\left(\xExactj -
\xGuessj\right)- \left(\xExactjm - \xGuessjm\right)\right)}\|\\
\leq& \|\projOrthNo{\Aj}{\currentDataSpacejm} -
\projOrthNo{\Ajm}{\currentDataSpacejm}\|\|{\xExactj - \xGuessj}\| +
\projOrthAjmZjmSingVal\| (\xExactj - \xGuessj) - (\xExactjm - \xGuessjm)\|
\end{align*}
where we have applied the triangle inequality and used  $\projOrthAjmZjmSingVal= \max_{\x\neq
0}\|\projOrthNo{\Ajm}{\currentDataSpacejm}{\x}\|/\|\x\|$. Eq.~\eqref{eq:weightsDiffAj} follows
from applying $\currentDatajmSingVal= \min_{\x\neq
0}\|\currentDatajm\x\|/\|\x\|$.
\end{proofofthing}

\begin{proofofthing}[Theorem \ref{thm:podDistanceGen}]
We proceed by leveraging
perturbation bounds for invariant subspaces. 
From the algebraic definition of a POD basis provided by
Eq.~\eqref{eq:podAlgebraicDefinition}, we can write 
a basis for $\podSubspaceIdealShort$ as
\begin{equation} \label{ideal Y}
\myIdealY{j}\defeq\podBasisIdeal
=
\myZ{j-1} \myIdealEtaMatrix{j} \myIdealTruncatedX{j}
	[\myIdealTruncatedLambda{j}]^{-\frac{1}{2}}.
\end{equation} 
The matrices $\myIdealTruncatedX{j}
\in\stiefel{\myTruncDimension}{\mySnapDimension}$ and
$\myIdealTruncatedLambda{j}
\in\myRR{\myTruncDimension\times\myTruncDimension}$ are obtained by
truncation, i.e., by retaining the first $\myTruncDimension$ columns of 
eigenvector matrix $\myIdealX{j}
\in\stiefel{\mySnapDimension}{\mySnapDimension}$ and the first
$\myTruncDimension$ rows and columns of the eigenvalue matrix
 $\myIdealLambda{j} \in\myRR{\mySnapDimension\times\mySnapDimension}$, which
 are associated with 
the eigenvalue decomposition
\begin{equation} \label{eigen decomposition}
\myIdealEtaMatrix{j} \myZ{j-1}^T \myIdealMetric \myZ{j-1} \myIdealEtaMatrix{j}
= \myIdealX{j} \myIdealLambda{j} [\myIdealX{j}]^T.
\end{equation}
Similarly, the second POD subspace $\podSubspacePracticalShort$ is 
spanned by a basis
\begin{equation} \label{ideal Y 2}
\myPracticalY{j}\defeq\podBasisPractical
=
\myZ{j-1} \myPracticalEtaMatrix{j} \myPracticalTruncatedX{j}
	[\myPracticalTruncatedLambda{j}]^{-\frac{1}{2}}
\end{equation} 
 where again $\myPracticalTruncatedX{j}$ and
 $\myPracticalTruncatedLambda{j}\in\myRR{\myTruncDimension\times\myTruncDimension}$ are the truncated counterparts of the
 eigenvectors 
 $\myPracticalX{j}
\in\stiefel{\mySnapDimension}{\mySnapDimension}$
 and eigenvalues 
$\myPracticalLambda{j} \in\myRR{\mySnapDimension\times\mySnapDimension}$
 associated with eigenvalue problem
\begin{equation} \label{eigen decomposition 2}
\myPracticalEtaMatrix{j} \myZ{j-1}^T \myPracticalMetric \myZ{j-1} \myPracticalEtaMatrix{j}
= \myPracticalX{j} \myPracticalLambda{j} [\myPracticalX{j}]^T.
\end{equation}
Noting that a diagonal column scaling of a basis does not affect the
associated subspace,
we can write $$\podSubspaceIdealShort = \range{\myZ{j-1} \myIdealEtaMatrix{j}
\myIdealTruncatedX{j}}\quad \text{and}\quad
\podSubspacePracticalShort = \range{\myZ{j-1} \myPracticalEtaMatrix{j}
\myPracticalTruncatedX{j}}=\range{\myZ{j-1} \myIdealEtaMatrix{j} \myCombined{j}  \myPracticalTruncatedX{j}},$$
with
$ \myCombined{j} \defeq [\myIdealEtaMatrix{j}]^{-1} \myPracticalEtaMatrix{j} $.
Thus, we can compute the distance between $\podSubspaceIdealShort$ and $\podSubspacePracticalShort$  in two steps: 1) compute the 
distance between subspaces spanned by bases $\myIdealTruncatedX{j}$ and
$\myCombined{j}  \myPracticalTruncatedX{j}$, and 2) apply a change in
coordinates by pre-multiplying each basis by $\myZ{j-1}
\myIdealEtaMatrix{j}$.

\underline{Step 1}. We employ 
eigenvector perturbation 
theory~\cite{davis1969,davis1970,IpsenOne,IpsenTwo,IpsenThree} to bound the distance between
$\range{\myIdealTruncatedX{j}}$ and $\range{\myCombined{j}
\myPracticalTruncatedX{j}}$.
First, we recall a
general eigenvector perturbation bound that characterizes
the distance between $\myTruncDimension$-dimensional invariant subspaces
$\myGeneralIdealSubspace\in\grassman{\myTruncDimension}{\nsnapshots}$
and
$\myGeneralPracticalSubspace\in\grassman{\myTruncDimension}{\nsnapshots}$ 
associated with matrices 
$\myGeneralIdeal\in\RR{\nsnapshots\times\nsnapshots}$
and
$\myGeneralPractical\in\RR{\nsnapshots\times \nsnapshots}$, respectively.

From Ref.~\cite{IpsenThree}, we have the following bound when both
$\myGeneralPractical$
and 
$\myGeneralIdeal$ 
are diagonalizable:
\begin{equation} \label{general eigen perturb}
 \distanceMetric{\myGeneralIdealSubspace}{\myGeneralPracticalSubspace}  \le  \myIdealCondNumber \myPracticalCondNumber
          \| \E \|_2 /.
\text{abssep}(\myGeneralPerpIdealTruncatedLambda,\myGeneralPracticalTruncatedLambda)
\end{equation}
Here, the difference
matrix is denoted by 
$
\E \defeq \myGeneralPractical - \myGeneralIdeal, 
$
$\myPracticalCondNumber$ denote the condition number of
the matrix whose columns comprise the first $\myTruncDimension$ eigenvectors of 
$\myGeneralPractical$, and
$\myIdealCondNumber$ denote condition numbers of
the matrix whose columns comprise the last $\ndof-\myTruncDimension$ eigenvectors of 
$\myGeneralIdeal$; these condition numbers evaluate to
one when the associated matrix is normal. Also,
$\myGeneralPracticalTruncatedLambda$ denotes the diagonal matrix of
eigenvalues associated
with the invariant subspace spanned by the first $\augSpaceDim$ eigenvectors
of the matrix $\myGeneralPractical$, while $\myGeneralPerpIdealTruncatedLambda$ denotes
the diagonal matrix of eigenvalues associated with 
the orthogonal complement to the
invariant subspace spanned by the first $\augSpaceDim$ eigenvectors of
$\myGeneralIdeal$.

To apply bound \eqref{general eigen perturb} in the present context,
we take 
\begin{gather}
\myGeneralIdeal  = 
\myIdealEtaMatrix{j} \myZ{j-1}^T  
\myIdealMetric \myZ{j-1} \myIdealEtaMatrix{j}
= \myIdealEtaMatrix{j} \myZ{j-1}^T  
(\myPracticalMetric + \myMetricPerturb) \myZ{j-1} \myIdealEtaMatrix{j}= \myIdealX{j} \myIdealLambda{j} [\myIdealX{j}]^T\\
\myGeneralPractical = 
\myCombined{j} 
 \myPracticalEtaMatrix{j} \myZ{j-1}^T \myPracticalMetric \myZ{j-1}
 \myPracticalEtaMatrix{j} 
\myCombined{j}^{-1},
\end{gather}
where the latter satisfies
\begin{equation} 
\myGeneralPractical\myCombined{j}\myPracticalX{j} = 
\myCombined{j}\myPracticalX{j} \myPracticalLambda{j}
.
\end{equation} 
Here, we have used Eqs.~\eqref{eigen decomposition}, \eqref{eigen
decomposition 2}, and $[\myCombined{j}\myPracticalX{j}]^{-1} =
[\myPracticalX{j}]^T\myCombined{j}^{-1}$.
Thus, we have $\myGeneralPracticalTruncatedLambda=\myPracticalTruncatedLambda{j}$ because $\myGeneralPractical$ is
similar to $\myPracticalEtaMatrix{j} \myZ{j-1}^T \myPracticalMetric \myZ{j-1}
 \myPracticalEtaMatrix{j}$; we also can write
 $\myGeneralPerpIdealTruncatedLambda = \myPerpIdealTruncatedLambda{j}$.
We also have $\myGeneralIdealSubspace = \range{\myIdealTruncatedX{j}}$ and
$\myGeneralPracticalSubspace = \range{\myCombined{j}
\myPracticalTruncatedX{j}}$.
Further, note that $\myIdealCondNumber=1$ because $\myGeneralIdeal$ is normal,
while $\myPracticalCondNumber=\kappa(\myCombined{j}
\myPracticalTruncatedX{j})$ due to non-normality of $\myGeneralPractical$.

To complete bound \eqref{general eigen perturb}, we bound $\| \E\|_2$
using commutativity of diagonal matrices and the relations
$[\myPracticalEtaMatrix{j}]^2 [\myIdealEtaMatrix{j}]^{-1} = \myCombined{j}
\myPracticalEtaMatrix{j}$ and $\myIdealEtaMatrix{j} =
\myPracticalEtaMatrix{j} \myCombined{j}^{-1}$ as
\begin{align}
\| \E \|_2 \hskip -.045in &= \| \myGeneralPractical - \myGeneralIdeal \|_2 \\[5pt]
          &= \|  \myCombined{j} \myPracticalEtaMatrix{j} \myZ{j-1}^T \myPracticalMetric \myZ{j-1}  \myPracticalEtaMatrix{j} \myCombined{j}^{-1} ~~- 
                      \myIdealEtaMatrix{j} \myZ{j-1}^T \left( \myPracticalMetric + \myMetricPerturb \right ) \myZ{j-1} \myIdealEtaMatrix{j} \|_2 \\[5pt]
          &= \|  \myCombined{j} \myPracticalEtaMatrix{j} \myZ{j-1}^T \myPracticalMetric \myZ{j-1} \myIdealEtaMatrix{j} ~~- 
                      \myIdealEtaMatrix{j} \myZ{j-1}^T \left( \myPracticalMetric + \myMetricPerturb \right ) \myZ{j-1} \myIdealEtaMatrix{j} \|_2 \\[5pt]
          &= \| \left ( \myCombined{j} \myPracticalEtaMatrix{j} - \myIdealEtaMatrix{j} \right ) \myZ{j-1}^T \myPracticalMetric \myZ{j-1} \myIdealEtaMatrix{j} ~~- 
                      \myIdealEtaMatrix{j} \myZ{j-1}^T \myMetricPerturb \myZ{j-1} \myIdealEtaMatrix{j} \|_2 \\[5pt]
          &= \| \left ( [\myPracticalEtaMatrix{j}]^2 [\myIdealEtaMatrix{j}]^{-1} - \myIdealEtaMatrix{j} \right ) \myZ{j-1}^T \myPracticalMetric \myZ{j-1} \myIdealEtaMatrix{j} ~~- 
                      \myIdealEtaMatrix{j} \myZ{j-1}^T \myMetricPerturb \myZ{j-1} \myIdealEtaMatrix{j} \|_2 \\[5pt]
          &= \|  (\myPracticalEtaMatrix{j} +
					\myIdealEtaMatrix{j})(\myPracticalEtaMatrix{j} -
					\myIdealEtaMatrix{j})[\myIdealEtaMatrix{j}]^{-1}
\myZ{j-1}^T \myPracticalMetric \myZ{j-1} \myIdealEtaMatrix{j} - 
                      \myIdealEtaMatrix{j} \myZ{j-1}^T \myMetricPerturb \myZ{j-1} \myIdealEtaMatrix{j} \|_2 \\[5pt]
          &\label{eq:finalEinequality}\le \|  (\myPracticalEtaMatrix{j} +
					\myIdealEtaMatrix{j})(\myPracticalEtaMatrix{j} -
					\myIdealEtaMatrix{j})[\myIdealEtaMatrix{j}]^{-1}
 \myZ{j-1}^T \myPracticalMetric \myZ{j-1} \myIdealEtaMatrix{j} \|_2 + \|
                      \myMetricPerturb \|_2 \| \myZ{j-1} \myIdealEtaMatrix{j} \|^2_2 . 
\end{align}
Substituting inequality \eqref{eq:finalEinequality} and the values of the
condition numbers
in bound \eqref{general eigen perturb} leads to
\begin{align} \label{general eigen perturb final}
\begin{split}
 \distanceMetric{\myGeneralIdealSubspace}{\myGeneralPracticalSubspace}  \le
 \kappa(\myCombined{j} \myPracticalTruncatedX{j})
          \Bigl(&\|  \left(\myPracticalEtaMatrix{j} +
					\myIdealEtaMatrix{j}\right)\left(\myPracticalEtaMatrix{j} -
					\myIdealEtaMatrix{j}\right)[\myIdealEtaMatrix{j}]^{-1}
\myZ{j-1}^T \myPracticalMetric \myZ{j-1} \myIdealEtaMatrix{j} \|_2 +\\
& \|
                      \myMetricPerturb \|_2 \| \myZ{j-1} \myIdealEtaMatrix{j} \|^2_2\Bigr)/
\text{abssep}(\myPerpIdealTruncatedLambda{j},\myPracticalTruncatedLambda{j})
\end{split}
\end{align}

\underline{Step 2}. We next use the result in Reference \cite[Proposition
2.6.14]{watkins2007matrix} that allows for a change of coordinates using a
nonsingular square transformation matrix.
In particular, if 
$\left[\myZ{j-1}
\myIdealEtaMatrix{j}\right]$ has full column rank, then
$\left[\myZ{j-1}
\myIdealEtaMatrix{j},\ \myZComplement\right]$ 
with $\myZComplement\in\stiefel{\ndof-\mySnapDimension}{\ndof}$ and $\range{\myZComplement} =
\range{\myZ{j-1}
\myIdealEtaMatrix{j}} ^\perp$
is nonsingular. Then, we have 
\begin{align} \label{general eigen perturb final rotate}
 \distanceMetric{\podSubspaceIdealShort}{\podSubspacePracticalShort}&\le \kappa(\myZ{j-1}
\myIdealEtaMatrix{j}) \distanceMetric{\myGeneralIdealSubspace}{\myGeneralPracticalSubspace}.
\end{align}
Here, we have used $\kappa(\left[\myZ{j-1}
\myIdealEtaMatrix{j},\ \myZComplement\right]) = \kappa(\left[\myZ{j-1}
\myIdealEtaMatrix{j}\right])$. Combining bounds \eqref{general eigen perturb final} and \eqref{general eigen perturb final rotate} yields
the desired result.
\end{proofofthing}

 \begin{proofofthing}[Theorem \ref{wellConditioned}]
\begin{align}
\|\augSpaceBasisj^T\Aj\augSpaceBasisj - \I\| &=
\|\augSpaceBasisj^T\Ajm\augSpaceBasisj - \I + \augSpaceBasisj^T(\Aj -
\Ajm)\augSpaceBasisj\|\\
&\leq \|
\left[\augSpaceBasisjm,\ \kryVecArg{j-1}\sqrt{\TArg{j-1}}\right]^T\Ajm
\left[\augSpaceBasisjm,\ \kryVecArg{j-1}\sqrt{\TArg{j-1}}\right] - \I\| +
\|\augSpaceBasisj\|^2\|\Aj - \Ajm\|\\
& = \|
\left[
\begin{array}{cc}
\augSpaceBasisjm\Ajm\augSpaceBasisjm - \I & \zero\\
\zero & \zero
\end{array}
\right]
\| +
\|\augSpaceBasisj\|^2\|\Aj - \Ajm\|,\\
& = \|\augSpaceBasisjm\Ajm\augSpaceBasisjm - \I\| +
\|\augSpaceBasisj\|^2\|\Aj - \Ajm\|,\\
\end{align}
where we have used $\augSpaceBasisArg{k}^T\A_k\kryVecArg{k} = \zero$,
$k\innatseq{j-1}$ (which holds under the stated assumptions) and 
$\sqrt{\TArg{k}}^T\kryVecArg{k}^T\A_k\kryVecArg{k} \sqrt{\TArg{k}} = \I$, 
$k\innatseq{j-1}$. By induction, we arrive at inequality \eqref{eq:wellConditioned},
where we have used $\augSpaceBasisArg{\jcompress +
1}\A_{\jcompress}\augSpaceBasisArg{\jcompress + 1} = \I$.
 \end{proofofthing}

\bibliography{references,eigrefs}
\bibliographystyle{siam}
\end{document}